\documentclass{amsart}

\usepackage[utf8]{inputenc}
\usepackage[T1]{fontenc}

\usepackage{amsmath,amssymb,amscd,verbatim,amsthm}
\usepackage[usenames,svgnames]{xcolor}
\usepackage{tikz}
\usepackage[matrix,arrow,curve]{xy}

\usepackage{hyperref} 

\usepackage{ifthen}
\usepackage{verbatim}
\newboolean{draft}
\setboolean{draft}{false}
\newboolean{longversion}
\setboolean{longversion}{true}
\iflongversion\else
\renewenvironment{proof}{\comment}{\endcomment}
\fi

\ifdraft
\newcommand{\TODO}[2][To do: ]{\textcolor{red}{\textbf{#1#2}}}
\else
\newcommand{\TODO}[2][]{}
\fi
\newcommand{\FIXME}{\TODO[Fix-me: ]}

\newcommand{\realring}{{\mathbb{K}}}

\newcommand{\CC}{{\mathbb{C}}}
\newcommand{\ZZ}{{\mathbb{Z}}}
\newcommand{\QQ}{{\mathbb{Q}}}
\newcommand{\RR}{{\mathbb{R}}}
\newcommand{\s}{w}
\renewcommand{\t}{{w'}}
\newcommand{\End}{{\operatorname{End}}}

\newcommand{\geh}{\mathfrak{g}}
\newcommand{\wt}{\mathrm{wt}}
\newcommand{\D}{\Pi}

\newcommand{\W}{W}
\newcommand{\clW}{\mathring{W}}
\newcommand{\Wa}[1][\CC]{{#1[\W]}}
\newcommand{\clWa}[1][\CC]{{#1[\clW]}}
\newcommand{\kW}[1][\CC]{{#1\W}}
\newcommand{\kclW}[1][\CC]{{#1\clW}}
\newcommand{\heckeW         }[2][\W]{{           \operatorname{H} (#1)(#2)}}
\newcommand{\heckeWW        }[1][\W]{{           \operatorname{H}\!#1     }}

\newcommand{\Wmax}{{w_0}}

\newcommand{\sg}[1][n]{{\mathfrak{S}_{#1}}}

\newcommand{\hecke}[2][n]           {{\operatorname{H}_{#1}(#2)}}

\newcommand{\affinehecke}[2][n]     {{\widetilde{\operatorname{H}}_{#1}(#2)}}

\newcommand{\affheckcent}[2][n]     {{\mathcal H}_{#1}(#2)}
\newcommand{\Y}{{Y}}

\newcommand{\rootspace  }[1][]{\mathfrak{h}^*_{#1}}
\newcommand{\rootlattice}{\rootspace[\ZZ]}
\newcommand{\corootlattice}{{\rootspace[\ZZ]}}
\newcommand{\coweightspace}[1][]{\mathfrak{h}_{#1}}
\newcommand{\coweightlattice}{\coweightspace[\ZZ]}
\newcommand{\clcoweightspace}[1][]{\coweightspace[#1]^0}
\newcommand{\weightspace}{\mathcal{P}}

\newcommand{\clcorootlattice}{\mathring R^\vee}
\newcommand{\coroot}{\alpha^\vee}
\newcommand{\roots}{R}
\newcommand{\coroots}{\roots^\vee}
\newcommand{\clroots}{\mathring\roots}
\newcommand{\clcoroots}{{\clroots^\vee}}
\newcommand{\clpositivecoroots}{\clroots^{\vee +}}

\newcommand{\pisequence}[1][r]{\pi_{i_#1} \cdots \pi_{i_1}}

\newcommand{\Des}{{\operatorname{D}_R}}
\newcommand{\Rec}{{\operatorname{D}_L}}

\newcommand{\height}{\operatorname{ht}}
\newcommand{\clrho}{\mathring\rho}
\newcommand{\rhoc}{{\rho^\vee}}
\newcommand{\clrhoc}{\clrho^\vee}
\newcommand{\x}{x}
\newcommand{\xc}{\x^\vee}

\newcommand{\dc}{d^\vee}
\newcommand{\lc}{{\lambda^\vee}}
\newcommand{\Lambdac}{\Lambda^\vee}
\newcommand{\len}{\ell}
\newcommand{\cl}{{\operatorname{cl}}}

\newcommand{\e}{\varepsilon}

\newcommand{\opi}{{\overline{\pi}}}
\newcommand{\pim}{{\nabla}}

\newcommand{\oT} {{\overline T}}
\newcommand{\lavee}{{\lambda^\vee}}

\newcommand{\ip}[1]{\langle #1 \rangle}

\newcommand{\init}{{\operatorname{init}}}

\newcommand{\setofroots}{\mathcal S}

\newtheorem{theorem}{Theorem}[section]
\newtheorem{lemma}[theorem]{Lemma}
\newtheorem{proposition}[theorem]{Proposition}
\newtheorem{corollary}[theorem]{Corollary}

\newtheorem{example}[theorem]{Example}

\newtheorem{problem}[theorem]{Problem}

\newtheorem{remark}[theorem]{Remark}

\newcommand{\q}{-\frac{q_1}{q_2}}

\newcommand{\quotientmap}{\cl}

\newcommand{\suchthat}{{\ |\ }}

\makeatletter
\newskip\@bigflushglue \@bigflushglue = -100pt plus 1fil

\def\bigcentering{\let\\\@centercr\rightskip\@bigflushglue%
\leftskip\@bigflushglue
\parindent\z@\parfillskip\z@skip}

\makeatother

\newcommand{\dynkin}[1]{
  \begin{tikzpicture}[>=latex,join=bevel,baseline=(current bounding box.east)]
    #1
  \end{tikzpicture}}

\newcommand{\dynkinAIIa}{\dynkin{
    \node (N_0) at (21bp,22bp) [draw,draw=none] {$0$};
    \node (N_1) at (6bp,12bp) [draw,draw=none] {$1$};
    \node (N_2) at (36bp,12bp) [draw,draw=none] {$2$};
    \draw [] (N_0) -- (N_1);
    \draw [] (N_1) -- (N_2);
    \draw [] (N_2) -- (N_0);
  }}

\newcommand{\dynkinCIIa}{\dynkin{
    \node (N_0) at (6bp,6bp) [draw,draw=none] {$0$};
    \node (N_1) at (36bp,6bp) [draw,draw=none] {$1$};
    \node (N_2) at (66bp,6bp) [draw,draw=none] {$2$};
    \draw [->,double] (N_0) -- (N_1);
    \draw [<-,double] (N_1) -- (N_2);
  }}

\newcommand{\dynkinGIIa}{\dynkin{
    \node (N_1) at (10bp,6bp) [draw,draw=none] {$0$};
    \node (N_2) at (36bp,6bp) [draw,draw=none] {$1$};
    \node (N_3) at (66bp,6bp) [draw,draw=none] {$2$};
    \draw [] (N_1) -- (N_2);
    \draw [->] (N_2) -- node[above] {$3$} (N_3);
  }}

\newcommand{\dynkinAIVb}{\dynkin{
    \node (N_0) at (6bp,6bp) [draw,draw=none] {$0$};
    \node (N_1) at (36bp,6bp) [draw,draw=none] {$1$};
    \node (N_2) at (66bp,6bp) [draw,draw=none] {$2$};
    \draw [<-,double] (N_0) -- (N_1);
    \draw [<-,double] (N_1) -- (N_2);
  }}

\newcommand{\dynkinDIIb}{\dynkin{
    \node (N_0) at (6bp,6bp) [draw,draw=none] {$0$};
    \node (N_1) at (36bp,6bp) [draw,draw=none] {$1$};
    \node (N_2) at (66bp,6bp) [draw,draw=none] {$2$};
    \draw [<-,double] (N_0) -- (N_1);
    \draw [->,double] (N_1) -- (N_2);
  }}

\begin{document}
\title[Hecke group algebras as quotients of affine Hecke algebras]
{Hecke group algebras\\ as quotients of affine Hecke algebras at level $0$}
\iflongversion
\author{Florent Hivert}
\address{LITIS (EA 4108), Université de Rouen,  
  Avenue de l'Université  BP12
  76801 Saint-Etienne du Rouvray, France and
  Institut Gaspard Monge (UMR 8049), France}
\email{florent.hivert@univ-rouen.fr}

\author{Anne Schilling}
\address{Department of Mathematics, University of California, One
 Shields Avenue, Davis, CA 95616, U.S.A.}
\email{anne@math.ucdavis.edu}

\author{Nicolas M.~Thiéry}
\address{Univ Paris-Sud, Laboratoire de Mathématiques d'Orsay,
  Orsay, F-91405; CNRS, Orsay, F-91405, France}
\email{Nicolas.Thiery@u-psud.fr}
\else
\author{Florent Hivert\addressmark{1}, Anne Schilling\addressmark{2}, 
and Nicolas M.~Thiéry\addressmark{2,3}}

\address{\addressmark{1} LITIS (EA 4108), Université de Rouen,  
  Avenue de l'Université  BP12
  76801 Saint-Etienne du Rouvray, France and
  Institut Gaspard Monge (UMR 8049)
\addressmark{2}Department of Mathematics, University of California, One
 Shields Avenue, Davis, CA 95616, U.S.A.\\
\addressmark{3}Univ Paris-Sud, Laboratoire de Mathématiques d'Orsay,
  Orsay, F-91405; CNRS, Orsay, F-91405, France}
\fi


\iflongversion
\fi

\def\showrevision$Revision: #1${#1}
\def\showdate$Date: #1 #2${#1}

\iflongversion
\else
\maketitle
\fi
\TODO{Switch out of draft version!}
\begin{abstract}
  \ifdraft
  This is revision \showrevision$Revision: 94$ of this document
  (\showdate$Date: 2008-11-10 15:59:12 +0100 (lun, 10 nov 2008) $).

  \fi
  The Hecke group algebra $\heckeWW[\clW]$ of a finite Coxeter group
  $\clW$, as introduced by the first and last authors, is obtained
  from $\clW$ by gluing appropriately its $0$-Hecke algebra and its
  group algebra. In this paper, we give an equivalent alternative
  construction in the case when $\clW$ is the 
  finite Weyl group
  associated to an affine Weyl group $\W$. Namely, we prove that, for
  $q$ not a root of unity of small order, $\heckeWW[\clW]$ is the natural
  quotient of the affine Hecke algebra $\heckeW[\W]{q}$ through
  its level $0$ representation.
  
  The proof relies on the following core combinatorial result: at
  level $0$ the $0$-Hecke algebra $\heckeW[\W]{0}$ acts transitively on
  $\clW$. Equivalently, in type $A$, a word written on a circle can be
  both sorted and antisorted by elementary bubble sort operators.
  We further show that the level $0$
  representation is a calibrated principal series representation
  $M(t)$ for a suitable choice of character $t$, so that the quotient
  factors (non-trivially) through the principal central
  specialization. This explains in particular the similarities between
  the representation theory of the 
  $0$-Hecke algebra $\heckeW[\clW]{0}$ and
  that of the affine Hecke algebra $\heckeW[W]{q}$ at this specialization.

  \iflongversion
  \else
  L'algèbre de Hecke groupe $\heckeWW[\clW]$ d'un groupe de
  Coxeter fini $\clW$, introduite par le premier et le dernier
  auteur, est obtenue en recollant de manière appropriée son algèbre
  de Hecke dégénérée et son algèbre de groupe. Dans cet article, nous
  donnons une construction alternative dans le cas où $\clW$ est
  un groupe de Weyl associé à un groupe de Weyl affine $\W$. Plus
  précisément, nous montrons que quand $q$ n'est ni nul ni une racine
  de l'unité, $\heckeWW[\clW]$ est le quotient naturel de
  l'algèbre de Hecke affine $\heckeW[\clW]{q}$ dans sa
  représentation de niveau $0$. Nous montrons de plus que la
  représentation de niveau $0$ est une représentation de
  série principale calibrée $M(t)$ pour un certain caractère $t$, de sorte que
  le quotient se factorise par la spécialisation centrale
  principale. Ce fait explique en particulier les similarités entre
  les théories des représentations de l'algèbre de Hecke dégénérée et
  de l'algèbre de Hecke affine sous cette spécialisation. 
  \fi
\end{abstract}

\iflongversion
\maketitle
\fi

\section{Introduction}

The starting point of this research lies in the striking similarities
between the representation theories of the degenerate (Iwahori)-Hecke
algebras on one side and of the principal central specialization of
the affine Hecke algebras on the other. For the sake of simplicity, we
describe those similarities for type A in this introduction, but they carry over
straightforwardly to any affine Weyl group $\W$ and its associated
finite Weyl group $\clW$.

The representation theory of the degenerate Hecke algebras $\hecke{0}$
for general type has been worked out by Norton~\cite{Norton.1979} and
special combinatorial features of type A have been described by
Carter~\cite{Carter.1986}. In particular, the projective modules $P_I$
of the type $A$ degenerate Hecke algebra $\hecke{0}$ are indexed by
subsets $I$ of $\{1, \dots, n-1\}$, and the basis of each $P_I$ is
indexed by those permutations of $n$ whose descent set is $I$.

\smallskip

On the other hand, the classification of the irreducible
finite-dimensional representations of the affine Hecke algebra
$\affinehecke{q}$ is due to Zelevinsky~\cite{Zelevinsky.1980}. They
are indexed by simple combinatorial objects called
multisegments. However, in this work, we are interested in a
particular subcategory related to a central specialization for which
the multisegments are also in bijection with subsets of $\{1, \dots,
n-1\}$. This relation is as follows. It is well known from
Bernstein and Zelevinsky~\cite{Bernstein_Zelevinsky.1977} and
Lusztig~\cite{Lusztig.1983}, that the center of the affine Hecke
algebra is the ring of symmetric polynomials
$\CC[\Y_1,\dots,\Y_n]^{\sg}$ in some particular elements $\Y_1,\dots,
\Y_n$ such that as vector space,
\begin{equation}
  \affinehecke{q} \simeq \hecke{q} \otimes \CC[\Y_1,\dots,\Y_n]\,.
\end{equation}
As a center, it acts by scalar multiplication in all irreducible
representations, and one way to select a particular class of
representations is to specialize the center in the algebra
itself. Thus any ring morphism from $\CC[\Y_1,\dots,\Y_n]^{\sg}$ to
$\CC$, or in other words any scalar alphabet, defines a quotient of
the affine Hecke algebra of dimension
\begin{equation}
\dim\left(\hecke{q}\right)\, 
\dim\left(\CC[\Y_1,\dots,\Y_n]/\CC[\Y_1,\dots,\Y_n]^{\sg}\right)
 = n!^2\, .
\end{equation}
Let us denote by $\affheckcent{q}$ the quotient of
$\affinehecke{q}$ obtained by the principal specialization of its center
to the alphabet $\frac{1-q^n}{1-q\phantom{^n}} := \{1, q,\dots, q^{n-1}\}$, that is
\begin{equation}
  \affheckcent{q} :=
  \affinehecke{q} \ /\ \langle e_i(\Y_1,\dots,\Y_n) - e_i(1, q, \dots,
q^{n-1})\ |\ i=1,\dots, n\rangle\,,
\end{equation}
where $e_i$ denote the elementary symmetric polynomials.  Then, in
this particular case, the multisegments of Zelevinsky are in
bijection with subsets $I$ of $\{1,\dots,n-1\}$ and the irreducible
representations $S_I$ of $\affheckcent{q}$ have their bases indexed
by descent classes of permutations. Thus one expects a strong link
between $\hecke{0}$ and $\affheckcent{q}$.  \smallskip

The goal of this paper is to explain this relation by means of the
Hecke group algebra $\heckeWW[\clW]$ introduced by the first and
the last
authors~\cite{Hivert_Thiery.HeckeSg.2006,Hivert_Thiery.HeckeGroup.2007}. Indeed,
by definition, $\heckeWW[\clW]$ contains naturally the degenerated
Hecke algebra $\heckeW[\clW]{0}$ and it was shown that the simple
modules of $\heckeWW[\clW]$, when restricted to
$\heckeW[\clW]{0}$ form a complete family of projective ones.
The relation comes from the fact that there is a natural surjective
morphism from the affine Hecke algebra $\heckeW{q}$ to
$\heckeWW[\clW]$. As a consequence the simple modules of
$\heckeWW[\clW]$ are also simple modules of $\heckeW{q}$
elucidating the similarities. This can be restated as follows:
\begin{theorem}
  For $q$ not a root of unity, there is a particular
  finite-dimensional quotient $\heckeWW[\clW]$ of the affine Hecke
  algebra $\heckeW{q}$ which contains the $0$-Hecke algebra
  $\heckeW[\clW]{0}$ and such that any simple $\heckeWW[\clW]$ module
  is projective when restricted to $\heckeW[\clW]{0}$.
\end{theorem}

The remainder of this paper is structured as follows.

In Sections~\ref{section.coxeter} and~\ref{section.affineWeyl}, we
briefly review the required material on Coxeter groups, Hecke
algebras, and Hecke group algebras, as well as on the central theme of
this paper: the level $0$ action of an affine Weyl group $\W$ on the
associated 
finite Weyl group $\clW$ and the corresponding level $0$
representation of the affine Hecke algebra on $\kclW$.

In Section~\ref{section.transitivity}, we prove the core combinatorial
property (Theorem~\ref{theorem.connected}) which states that, at level
$0$, the affine $0$-Hecke algebra $\heckeW[\W]{0}$ acts transitively
on the chambers of $\clW$ (or equivalently on the finite Weyl
group). We first treat type $A$ where Theorem~\ref{theorem.connected}
states that a word written on a circle can be both sorted and
antisorted by elementary bubble sort operators (explicit (anti)sorting
algorithms are also provided for types $B$, $C$, and $D$).
We proceed with a type-free geometric proof of
Theorem~\ref{theorem.connected}. The ideas used in the proof are
inspired by private notes on finite-dimensional representations of
quantized affine algebras by Kashiwara~\cite{Kashiwara.2008}, albeit
reexpressed in terms of alcove walks.
We also mention connections with affine crystals.

In Section~\ref{section.quotient} we prove the main result of the
paper, namely that for $q$ not a root of unity of small order, the
Hecke group algebra is the natural quotient of the (extended) affine
Hecke algebra through its representation at level $0$
(Theorem~\ref{theorem.quotient}).  The proof relies on the results
from the subsequent sections, namely
Corollary~\ref{corollary.piGenerators} for $q=0$ and
Theorem~\ref{theorem.quotientQualibrated} for $q$ non-zero and not a
root of unity. Both yield a proof for generic $q$.

In Section~\ref{section.generators}, we derive new sets of generators
for the Hecke group algebra of a finite Weyl which, together with the
combinatorial results of Section~\ref{section.transitivity} give
Corollary~\ref{corollary.piGenerators}.

Unlike for the affine Weyl group $\W$, and interestingly enough, the
torus $Y$ does not degenerate trivially. In
Section~\ref{section.calibrated}, we describe precisely this
degeneracy, and show that, for a suitable choice of character on $Y$,
the level $0$ representation is a calibrated principal series
representation (Theorem~\ref{theorem.isomorphismCalibrated}). This
allows to us refine Theorem~\ref{theorem.quotient} to $q$ not a root
of unity.

Altogether, Theorems~\ref{theorem.quotient}
and~\ref{theorem.isomorphismCalibrated} can be interpreted as two new
equivalent alternative constructions of the Hecke group algebra, while
the latter provides a parametrized family of maximal decompositions of
its identity into idempotents
(Corollary~\ref{corollary.decompositionOfIdentity}).

\section{Coxeter groups, Hecke  algebras, and Hecke group algebras}
\label{section.coxeter}

In this and the next section, we briefly recall the notations and
properties of Coxeter groups, (affine) Weyl groups, their Hecke and
Hecke group algebras, as well as root systems and alcove walks that we
need in the sequel. For further reading on those topics, we refer the
reader
to~\cite{Humphreys.1990,kac.1990,Macdonald.2003,Bjorner_Brenti.2005,Ram.2006}.

\subsection{Coxeter groups and their geometric representations}
\label{subsection.coxeter.groups}

Let $\W$ be a Coxeter group and $I$ the index set of its Dynkin
diagram.  Denote by $(s_i)_{i\in I}$ its simple reflections and by
$\Wmax$ its maximal element (when $\W$ is finite). A presentation of
$W$ is given by the generators $s_i$ together with their quadratic and
braid-like relations:
\begin{equation}
	s_i^2=1 \qquad \text{and} \qquad  
	\underbrace{s_is_j\cdots}_{m(i,j)}= \underbrace{s_js_i\cdots}_{m(i,j)} \quad  \text{for $i\ne j$,}
\end{equation}
where the $m(i,j)$'s are integers depending on $\W$.

For $J\subset I$, write $W_J$ for the parabolic subgroup
generated by $(s_i)_{i\in J}$. The left and right descent sets of an
element $w\in\W$ are respectively 
\begin{equation*}
	\Rec(w) := \{ i \in I \mid s_i w < w \} \qquad \text{and} \qquad 
	\Des(w) := \{ i \in I \mid w s_i < w\}\,.
\end{equation*}

The Coxeter group $\W$ can be realized geometrically as follows. Take the
module $\rootspace:=\rootspace[\realring]:=\bigoplus_{i\in I}
\realring \alpha_i$ and its $\realring$-dual
$\coweightspace:=\coweightspace[\realring] := \bigoplus_{i\in I}
\realring \Lambdac_i$, with the natural pairing $\langle
\Lambdac_i,\alpha_j\rangle=\delta_{ij}$.  The $\alpha_i$ are the
\emph{simple roots}, and the $\Lambdac_i$ the \emph{fundamental
coweights}. The \emph{simple coroots} are given by
$\coroot_i:=\sum_j a_{i,j} \Lambdac_i$, where $M=(a_{i,j})_{i,j\in I}$ with
$a_{i,j} = \langle \coroot_i , \alpha_i \rangle$ is the
(generalized) Cartan matrix for $\W$ with coefficients in a ring $\realring\subset\RR$.
The Coxeter group acts on $\coweightspace$ by the number game:
\begin{equation}
  s_i(\xc) := \xc - \ip{\xc,\alpha_i} \coroot_i
  \quad \text{for $\xc\in \coweightspace$},
\end{equation}
and on $\rootspace$ by the dual number game:
\begin{equation}
  s_i(\x) := \x - \ip{\coroot_i,\x} \alpha_i
  \quad \text{for $\x\in \rootspace$}.
\end{equation}

Denote by $R:=\{w(\alpha_i) \suchthat w\in W, i\in I\}$ the set of
\emph{roots}, and by $R^\vee:=\{w(\coroot_i) \suchthat w\in W, i\in
I\}$ the set of \emph{coroots}.
To each root $\alpha$ corresponds the reflection $s_\alpha$ across the
associated coroot $\coroot$ and along the hyperplane $H_\alpha$ which
splits $\coweightspace$ into a  positive $H_\alpha^+$  and a negative
half-space $H_\alpha^-$:
\iflongversion
\begin{equation}
\begin{split}
  H_{\alpha}   &:= \left\{ \xc\in\coweightspace \suchthat \ip{\xc,\alpha}=0 \right\},\\
  H_{\alpha}^+ &:= \left\{ \xc\in\coweightspace \suchthat \ip{\xc,\alpha}>0 \right\}, \\
  H_{\alpha}^- &:= \left\{ \xc\in\coweightspace \suchthat \ip{\xc,\alpha}<0 \right\}.
\end{split}
\end{equation}
\else
\begin{equation*}
  H_{\alpha}   := \left\{ \xc\in\coweightspace \suchthat \ip{\xc,\alpha}=0 \right\} \quad
  H_{\alpha}^+ := \left\{ \xc\in\coweightspace \suchthat \ip{\xc,\alpha}>0 \right\} \quad
  H_{\alpha}^- := \left\{ \xc\in\coweightspace \suchthat \ip{\xc,\alpha}<0 \right\}\,.
\end{equation*}
\fi

Take now $\realring=\RR$. Define the \emph{fundamental chamber} as the
open simplicial cone $C:=\{\xc \suchthat \ip{\xc,\alpha_i} > 0,
\forall i\in I \}$. For each root $\alpha$, the fundamental chamber
$C$ lies either entirely in $H_\alpha^+$ or in $H_\alpha^-$; $R$
splits accordingly into the sets of \emph{positive roots} $R^+:=\{\alpha\in R \mid C \subseteq 
H_\alpha^+\}$ and of \emph{negative roots} $R^- := \{ \alpha \in R \mid C \subseteq H_\alpha^-\} = 
- R^+$.

The closure $\overline C$ of $C$ is a fundamental domain for the
action of $\W$ on the \emph{Tits cone} $U:=\bigcup_{w\in W} w(\overline
C)$, and the elements $w$ of $\W$ are in bijection with the
\emph{chambers} $w(C)$.
This bijection induces both a left and a right actions of $\W$ on the
chambers. The right action is particularly nice as the chambers $w(C)$
and $w(C).s_i = ws_i(C)$ share a common wall. Any sequence
$i_1,\ldots,i_r$ gives therefore rise to a sequence of adjacent
chambers $C,\ s_{i_1}(C),\ s_{i_1}s_{i_2}(C),\ \dots,\
(s_{i_1}\!\cdots s_{i_r}) (C)$ from $C$ to $w(C)$ (where
$w=s_{i_1}\cdots s_{i_r}$), called a \emph{gallery}. For short, we
often denote this gallery by just $i_1,\ldots,i_r$.

\subsection{(Iwahori)-Hecke algebras}

Let $\W$ be a Coxeter group and $q_1$ and $q_2$ two complex numbers.
When defined, set $q=:\q$.  The (generic, Iwahori)
$(q_1,q_2)$-Hecke algebra $\heckeW{q_1,q_2}$ of $\W$ is the
$\CC$-algebra generated by the operators $T_i$ subject to the
quadratic and braid-like relations:
\begin{equation}
  (T_i-q_1)(T_i-q_2)=0 \qquad \text{ and } \qquad
  \underbrace{T_iT_j\cdots}_{m(i,j)}= \underbrace{T_jT_i\cdots}_{m(i,j)} \quad  \text{for $i\ne j$.}
\end{equation}
Its dimension is $|W|$, and a basis is given by the elements
$T_w:=T_{i_1}\cdots T_{i_r}$ where $w\in W$ and $i_1,\dots,i_r$ is a
reduced word for $w$. The right regular representation of
$\heckeW{q_1,q_2}$ is given by
\begin{equation}
  T_w T_i =
  \begin{cases}
    (q_1+q_2) T_w - q_1 q_2 T_{ws_i} & \text{ if $i$ descent of $w$,}\\
    T_{ws_i}                 & \text{ otherwise.}
    \end{cases}
\end{equation}
Define the unique operators $\oT_i$ such that $T_i + \oT_i = q_1 +
q_2$. They satisfy the same relations as the $T_i$, and further $T_i
\oT_i = \oT_i T_i = q_1q_2$.

At $q_1=1, q_2=-1$ (so $q=1$), we recover the usual group algebra
$\Wa$ of $\W$; in general, when $q_1+q_2=0$ one still recovers $\Wa$
up to a scaling of the generators: $s_i=\frac1{q_1} T_i$. 
Note that when $q_1$ and $q_2$ are non-zero and $q$ is not a root of
unity $\heckeW[\W]{q_1,q_2}$ is still isomorphic to $\Wa$, but in a
non-trivial way. 
On the
opposite side, taking $q_1=0$ and $q_2\ne 0$ (so $q=0$) yields the
\emph{$0$-Hecke algebra} $\heckeW[\W]{0}$; it is also a monoid algebra for the
\emph{$0$-Hecke monoid} $\{\pi_w
\suchthat w\in W \}$ generated by the idempotents
$\pi_i:=\frac1{q_2}T_i$. At $q_1=q_2=0$, one obtains the nilCoxeter
algebra. Traditionally, and depending on the application in mind,
different authors choose different specializations of $q_1$ and $q_2$,
typically $q_1=q$ and $q_2=-1$ (cf.~\cite{Wikipedia.HeckeAlgebra}), or
$q_1=t^{\frac12}$ and $q_2=t^{-\frac12}$ (cf. for
example~\cite{Ram.2008}). For our needs, keeping the two eigenvalues
generic yields more symmetrical formulas which are also easier to
specialize to other conventions. 
There also exists a more general definition of the Hecke algebra by
allowing a different pair of parameters $(q_1,q_2)$ for each conjugacy
class of reflections in $W$. For the sake of simplicity, we did not
try to extend the results presented in this paper to this larger
setting, but would not expect specific difficulties either.

We may realize the $0$-Hecke monoid geometrically on $\coweightspace$
as follows. For each $i\in I$, define the (half-linear) idempotent
$\pi_i$ (resp.  $\opi_i$) which projects onto the negative
(resp. positive) half space with respect to the root $\alpha_i$:
\begin{equation} \label{eq:pi}
  \pi_i(\xc) :=
  \begin{cases}
    s_i(\xc) & \text{ if $\ip{\xc,\alpha_i} > 0$},\\
    \xc & \text{ otherwise;}
  \end{cases}
  \quad
  \opi_i(\xc) :=
  \begin{cases}
    s_i(\xc) & \text{ if $\ip{\xc,\alpha_i} < 0$},\\
    \xc & \text{ otherwise.}
  \end{cases}
\end{equation}
As with the reflection $s_i$, these projections map chambers to chambers.
None of the projections $\pi_1,\dots,\pi_n$ fix the fundamental
chamber, and (when $\W$ is finite) all of them fix the negative
chamber. The correspondence between chambers and Weyl group elements
induces an action on the group $\W$ itself: this is the usual right
regular actions of the $0$-Hecke monoid, where $\pi_i$ adds a left
descent at position $i$ if it is not readily there, and $\opi_i$ does
the converse.  The action of
the $\pi_i$'s can be depicted by a graph on $\W$, with an $i$-arrow
from $w$ to $w'$ if $\pi_i(w)=w'$. Examples of such graphs are given
in Figure~\ref{figure.antisorting} (ignoring the $0$-arrows).

Let $\kW$ be the vector space of dimension $|\W|$ spanned by $W$.
Except for the nilCoxeter algebra ($q_1=q_2=0$), the Hecke algebra
$\heckeW{q_1,q_2}$ can be realized as acting on $\kW$ by
interpolation, mapping $T_i$ to $(q_1+q_2) \pi_i - q_1 s_i$. This
amounts to identify $\kW$ with the right regular representation of
$\heckeW{q_1,q_2}$ via $w\mapsto q_2^{-\ell(w)} T_w$, where $\ell(w)$
is the length of $w$.
Through this mapping, $\oT_i = (q_1+q_2) \opi_i - q_2 s_i$.

\subsection{Hecke group algebras}
\label{subsection.heckeGroupAlgebras}

Let now $\W$ be a finite Coxeter group. As we have just seen, we may
embed simultaneously the Hecke algebra $\heckeW{0}$ and the group
algebra $\Wa$ in $\End(\kW)$, via their right regular representations.
The \emph{Hecke group algebra} $\heckeWW$ of $\W$ is the smallest
subalgebra of $\End(\kW)$ containing them both
(see~\cite{Hivert_Thiery.HeckeGroup.2007}). It is therefore generated
by $(\pi_i)_{i\in I}$ and $(s_i)_{i\in I}$, and by interpolation it
contains all $q_1,q_2$-Hecke algebras where $(q_1,q_2) \neq
(0,0)$\footnote{However, the nilCoxeter algebra does not embed
  naturally. More precisely, up to a scalar there is a single nilpotent
  element $d_i := 1 + s_i -2\pi_i$ in the algebraic span of $s_i$ and
  $\pi_i$. A direct calculation shows that, for example, $d_1$ and
  $d_2$ do not satisfy the braid relations.}.

\FIXME{why the nilCoxeter algebra does \emph{not} embed} A basis for
$\heckeWW$ is given by $\{ \s\pi_\t \suchthat \Des(\s) \cap \Rec(\t) =
\emptyset\}$.  A more conceptual characterization is as follows: call
a vector $v$ in $\kW$ \emph{$i$-left antisymmetric} if $s_iv=-v$;
then, $\heckeWW$ is the subalgebra of $\End(\kW)$ of those operators
which preserve all $i$-left antisymmetries~\cite{Hivert_Thiery.HeckeGroup.2007}.

\section{Affine Weyl groups, Hecke algebras, and their level $0$ actions}
\label{section.affineWeyl}

Now let $\W$ be an affine Weyl group, with index set
$I:=\{0,\dots,n\}$ and Cartan matrix $M$. We always assume that
$\W$ is irreducible. We denote respectively by $a_i$ and
$a_i^\vee$ the coefficients of the canonical linear combination
annihilating the columns and rows of $M$, respectively.

In the sequel, we stick to the number game / dual number game
geometric setting of Section~\ref{subsection.coxeter.groups}.
(see also Figure~\ref{figure.alcoves})
This differs slightly from the usual setting for affine or Kac-Moody
Lie algebras~\cite{kac.1990}; it turns out to be simpler yet
sufficient for our purpose. Note first that $R:=\{w(\alpha_i)
\suchthat w\in W, i \in I\}$ is the set of \emph{real roots}; by
abuse, we call them roots, as the imaginary roots do not play a role
for our purposes. The geometric representations $\rootlattice$ and
$\coweightspace$ defined in Section~\ref{subsection.coxeter.groups}
correspond to the root lattice and the coweight space, respectively; we
do not use the central extension by $c:=\sum_{i=0}^n a_i^\vee
\coroot_i$. As a consequence, the coroot lattice $\bigoplus_{i\in I}
\ZZ \coroot_i$does not embed faithfully in $\coweightlattice$ (since
$c=0$ in $\coweightlattice$).  In particular, the set of coroots
$\coroots$ is finite, and (essentially) coincides with the set
$\clcoroots$ of coroots of $\clW$.
We also keep separate the dual lattices, without embedding them in
a single ambient space endowed with an inner product.

\subsection{Affine Weyl groups and alcove walks}
\label{subsection.affine.weyl.groups}
\iflongversion
\begin{figure}
  \begin{tikzpicture}[x={(-1cm,1cm)}, y={(1cm,1cm)}]
  \tikzstyle{point}=[circle,draw,fill=black,inner sep=0mm, minimum size=1mm]
  \tikzstyle{ref}=[inner sep=0mm, minimum size=0mm]
  \tikzstyle{alcove}=[DarkRed]
  \tikzstyle{s0}=[black]
  \tikzstyle{s1}=[DarkBlue]
  
  \node[ref] (tl) at (4,-2)  {};
  \node (tr) at (-2,4)  {};
  \node (bl) at (2,-4)  {};
  \node (br) at (-4,2)  {};

  \node[ref] (L0) at ( 1, 0) {};
  \node[ref] (L1) at ( 0, 1) {};
  \node[ref] (A0) at ( 2,-2) {};
  \node[ref] (A1) at (-2, 2) {};

  \draw[-] (intersection of A0--A1 and bl--tl)[very thin] --
           (intersection of A0--A1 and br--tr) node [below]{$\coweightspace^0$};
  \draw[-,alcove] (intersection of L0--L1 and bl--tl)[very thin] --
           (intersection of L0--L1 and br--tr) node [below]{$\coweightspace^1$};

  \node[point] (L) at (.5,.5)  [label=below:$\rhoc$] {};

  \node[alcove] at (1.15,0.15) {$\phantom{{}^1}0^1$};

  \draw[-] (intersection of 0,0--L0 and bl--br)[very thin,s1] -- 
           (intersection of 0,0--L0 and tl--tr)
           node [above right]{$H_{\alpha_{1,0}}=H_{\alpha_1}$};
  \draw[-] (intersection of 0,0--L1 and bl--br)[very thin,s0] -- 
           (intersection of 0,0--L1 and tl--tr)
           node [above right]{$H_{\alpha_{1,1}}=H_{\alpha_0}$};
  \foreach \x/\c in { -2/s1, -1/s0, 2/s1, 3/s0} {
    \node (start) at (intersection cs: first line  = {(0,0)--(1-\x,\x)},
                     second line = {(tl)--(bl)}) {};
    \node (end)   at (intersection cs: first line  = {(0,0)--(1-\x,\x)},
                     second line = {(tr)--(br)}) {};
    \draw[-,\c] (start) -- (end)[very thin] node [right]{$H_{\alpha_1,\x}$};
  }

  \foreach \x / \l in { -1/$s_0(C)$, 0/$C$, 1/$s_1(C)$ }
    \node at ( .9-1.6*\x,.9+1.6*\x) {\l};

  \foreach \x/\c in { -2/s0, -1/s1, 0/s0, 1/s1, 2/s0, 3/s1} {
    \draw[-,\c] (-.05+\x,0.95-\x) -- (.05+\x,1.05-\x);
  }
  \foreach \x / \l in { -2/$s_0s_1(A)$, -1/$s_0(A)$, 0/$A$,
    1/$s_1(A)$, 2/$s_1s_0(A)$ } {
    \node[alcove] at ( .6-1.15*\x,.6+1.15*\x) {\l};
    \node[point,alcove] at (.5-\x,.5+\x) {};
  }

  \foreach \x / \l in { -1.5, -0.5, 0.5, 1.5 }{
    \node at ( 0.8-1.83*\x,1  +1.77*\x) {+};
    \node at ( 1  -1.77*\x,0.8+1.83*\x) {-};
}

  \draw[->,thick,black] (0,0) -- (L0) node [below left ]{$\Lambdac_0$};
  \draw[->,thick,DarkBlue] (0,0) -- (L1) node [below right]{$\Lambdac_1$};

  \draw[->,thick] (0,0) -- (A0) node [below]{$\coroot_0$}; 
  \draw[->,thick,DarkBlue] (0,0) -- (A1) node [below]{$\coroot_1$}; 
\end{tikzpicture}

  \caption{Realization of the alcove picture at the level $1$
    hyperplane $\coweightspace^1$ of the coweight space
    $\coweightspace$ in type $A_1^{(1)}$.}
  \label{figure.alcoves}
\end{figure}
\fi

Let $\delta:=\sum_{i\in I} a_i\alpha_i$ be the so-called \emph{null root}%
\footnote{Beware that this is not a root in the current setting!}.
The level of an element $\xc$ of $\coweightspace$ is given by
$\ell(\xc) = \ip{\xc, \delta}$; in particular, and by construction,
all the coroots are of level $0$. Since $\delta$ is fixed by
$\W$, the affine hyperplanes $\coweightspace^\ell:= \{\xc
\suchthat \ip{\xc,\delta}=\ell\}$ are stabilized by $\W$.

At level $0$, the action $\cl$ of the affine Weyl group $\W$ on
$\clcoweightspace$ reduces to that of a finite Weyl group
$\clW:=\cl(W)$; in fact $\clW = \langle s_1,\dots,s_n\rangle$,
assuming an appropriate labeling of the Dynkin diagram. This induces a
right action of $\W$ on $\clW$: for $w$ in $\clW$ and $s_i\in W$,
$w.s_i:=w\cl(s_i)$, where $\cl:\,W \to \clW$ denotes the canonical
quotient map.  We denote respectively by $\clroots:=\{w(\alpha_i)
\suchthat w\in W, i=1,\dots,n\}$ and $\clcoroots:=\{w(\coroot_i)
\suchthat w\in W, i=1,\dots,n\}$ the sets of
roots and coroots of $\clW$. 
The coroot $\coroot_0$ is of the form $\coroot_0=\epsilon\coroot$ where $\coroot \in \clpositivecoroots$
and $\epsilon<0$. In the untwisted case, $\epsilon=-1$ so that  $\coroots=\clcoroots$. 
In the other cases $\coroots$ and $\clcoroots$ may differ by the orbit of $\coroot_0$.

The reflections in $\W$ are given by
\begin{equation}
  \{ s_{\alpha,m} := s_{\alpha-m\delta} \suchthat \alpha\in \clroots^+ \text{ and } m\in c_\alpha\ZZ\}\,.
\end{equation}
Here $s_{\alpha,m}$ is the reflection across the hyperplane
$H_{\alpha,m}:=H_{\alpha-m\delta}$ along the 
coroot $\alpha^\vee$ of $\clW$, and $c_\alpha\in \QQ$ 
($c_\alpha=1$ always in the untwisted case; for the twisted case see
Kac~\cite[Proposition 6.5]{kac.1990}). 

At level $\ell$, each positive root $\alpha$ of $\clW$ gives rise to a
family $(H^\ell_{\alpha,m})_{m\in c_\alpha\ZZ}$ of parallel reflection
hyperplanes (which all collapse to $H_{\alpha}^0$ at level $0$):
\begin{equation}
  H_{\alpha,m}^\ell := H_{\alpha-m\delta} \cap \coweightspace^\ell
  = \{ \xc\in\coweightspace^\ell \suchthat \ip{\xc,\alpha}=\ell m \}\,.
\end{equation}
The Tits cone is $\{\xc \suchthat \ip{\xc,\delta}>0\}$, and slicing it
at level $\ell>0$ gives rise to the \emph{alcove picture} (see Figure~\ref{figure.alcoves}). The
\emph{fundamental alcove} $A:=C\cap \coweightspace^\ell$ is a simplex,
and the \emph{alcoves} $w(A)$ in its orbit form a tessellation of
$\coweightspace^\ell$.  
Each gallery $C,s_{i_1}(C),\dots,(s_{i_1}\! \cdots s_{i_r})(C)$
induces an \emph{alcove walk} $A,s_{i_1}(A),$ $\dots,
(s_{i_1}\! \cdots s_{i_r})(A)$. As for galleries, we often denote this
alcove walk by just $i_1,\ldots,i_r$. 

For a simple coroot $\coroot_i$, let $c_i=c_{\alpha_i}$ and define 
$t_{\coroot_i}=s_{\alpha_i,c_i}s_{\alpha_i,0}$; at level $\ell$,
$t_{\coroot_i}$ is the composition of two reflections along parallel
hyperplanes, and acts as a translation by $c_i \ell\coroot_i$. For any
$\lc=\sum_{i=1}^n \lambda_i \coroot_i$ in the coroot lattice $\clcorootlattice$ of $\clW$,
set $c(\lc)=\sum_{i=1}^n c_i \lambda_i \coroot_i$. Then, in general,
$t_\lc:\,\coweightspace\to \coweightspace$ defined by
\begin{equation}
  t_\lc(\xc) = \xc + \ell(\xc) c(\lc)
\end{equation}
belongs to $\W$. More specifically, $t_\lc = s_{i_1}\cdots s_{i_r}$,
where $i_1,\dots,i_r$ is an alcove walk from $A$ to the translated
alcove $t_\lc{A}$. By abuse, we call $t_\lc$ a \emph{translation} of $W$.
This gives the usual semi-direct product decomposition $W = \clW
\rtimes \clcorootlattice$. In particular, $\cl: W\mapsto \clW$ is the
group morphism which kills the translations $t_\lc$,
$\lc\in\clcorootlattice$.

The \emph{fundamental chamber} for $\clW$ is the open simplicial cone
\begin{displaymath}
\{\xc\in \coweightspace^\ell \suchthat \ip{\xc,\alpha_i} > 0, \forall i=1,\dots,n \}\,.
\end{displaymath}
We denote by $0^\ell$ the intersection point of its walls
$(H^\ell_{\alpha_i})_{i=1,\dots,n}$. 
The orientation of the alcove walls is the periodic orientation where
only points infinitely deep inside the fundamental chamber for $\clW$
is on the positive side of all walls.
Consider an $i$-\emph{crossing} for $i\in\{0,\dots,n\}$ from an alcove
$w(A)$ to the adjacent alcove $ws_i(A)$, and let $H_{\alpha,m}$ the
crossed affine wall.  The crossing is \emph{positive} if $ws_i(A)$ is
on the positive side of $H_{\alpha,m}$, and \emph{negative} otherwise.
For an alcove walk $i_1,\dots,i_k$, define
$\epsilon_1,\dots,\epsilon_r$ by $\epsilon_k=1$ if the $k$th crossing
is positive and $-1$ otherwise.

The \emph{height} of an alcove $w(A)$ is given by
$\height(w(A))=\frac12(\epsilon_1+\dots+\epsilon_k)$, for any alcove
walk $i_1,\dots,i_k$ from $A$ to $w(A)$. This is well-defined, since
$\epsilon_1+\dots+\epsilon_k$ counts the number of hyperplanes
$H_{\alpha,m}$ separating $A$ from $w(A)$, where those with $w(A)$ on
the positive side are counted positively, and the others negatively.

\begin{remark}
  \label{remark.height}
  The height of the alcove $t_\lc(A)$ coincides with the \emph{height}
  of the coroot $\lc$ of $\clW$, $\height(\lc) := \ip{\lc,\clrho}$,
  where $\clrho:=\frac12\sum_{\alpha\in \clroots^+} \alpha$. In particular, a
  coroot is of height one if and only if it is a simple coroot
  ($\clrho$ is also the sum of the fundamental weights of $\clW$).
\end{remark}
\begin{proof}
  For each positive root $\alpha$ of $\clW$, the family of parallel
  hyperplanes $(H_{\alpha,m})_{m\in c_\alpha\ZZ}$ contributes to
  $\epsilon_1+\dots+\epsilon_k$
  the (relative) number of those separating $\frac \ell{n+1}\rho^\vee$
  and $\frac \ell{n+1}\rho^\vee + \ell c(\lc)$; this is given by
  $\ip{\lc,\alpha}$. The result follows by summing up over all
  positive roots.
\end{proof}

\subsection{Affine Hecke algebras}

The affine Hecke algebra of $\W$ is $\heckeW{q_1,q_2}$. In particular,
it is isomorphic to $\heckeW[\clW]{q_1,q_2} \otimes \CC[Y]$, where 
\begin{equation}
  \CC[Y]:=\CC.
  \{Y^\lc\suchthat \lc\in \clcorootlattice \}
\end{equation}
is the group algebra of the
coroot lattice.  The $Y^\lc$'s have an expression in terms of the $T_i$'s
which generalizes that for translations $t_\lc$ in the affine Weyl
group~\cite[Equation (3.2.10)]{Macdonald.2003}:
\begin{equation}
\label{equation.Y}
  Y^\lc := (\frac{1}{\sqrt{-q_1q_2}} T_{i_1})^{\epsilon_1}\cdots 
  (\frac{1}{\sqrt{-q_1q_2}} T_{i_r})^{\epsilon_r}
  = (-q_1q_2)^{-\height(\lc)} T_{i_1}^{\epsilon_1}\cdots T_{i_r}^{\epsilon_r} \;,
\end{equation}
where $i_1,\dots,i_r$ is an alcove walk from $A$ to $t_\lc(A)$.
The center of $\heckeW{q_1,q_2}$ is the subring of invariants $Y^W :=
\{ p\in Y \suchthat p.w =p\}$. In type $A$, this is the ring of
symmetric functions.

As for $\W$, the geometric realization at level $0$ induces an action
$\cl$ of the $0$-Hecke monoid $\langle \pi_i \suchthat i \in I
\rangle$ on the chambers of $\clW$, and therefore on $\clW$ itself:
\begin{equation}
  \label{equation.level0heckeaction}
  w.\cl(\pi_i) :=
  \begin{cases}
    w s_i & \text{ if $\pi_i(w^{-1}(\clrhoc))=w^{-1}(\clrhoc)$, that
      is $\ip{w^{-1}(\clrhoc), \alpha_i} > 0$,}\\
    w & \text{otherwise,}
  \end{cases}
\end{equation}
where $\clrhoc= \frac{1}{2}\sum_{\coroot \in \clcoroots} \coroot$ is the
canonical representative of the fundamental chamber of $\clcoroots$.
Geometrically, it can be interpreted as a quotient of the action at
level $\ell$ by identifying a point in a chamber at level $0$ with a
point infinitely deep inside the corresponding chamber for $\clW$ at
level $\ell$.
We recognize the usual action of $\pi_1,\dots,\pi_n$, where
$w.\pi_i=w$ if $i$ is a (right) descent of $w$ and $w.\pi_i=ws_i$
otherwise. By extension $0$ is called an \emph{(affine) descent} if
$w.\pi_0=w$. Since there is no ambiguity, we write $w.\pi_i$ for
$w.\cl(\pi_i)$.  Let us relate affine descents and positivity of
crossings.
\begin{remark}
  \label{remark.positiveCrossing}
  Consider an $i$-\emph{crossing} for $i\in\{0,\dots,n\}$ from an
  alcove $w(A)$ to the adjacent alcove $ws_i(A)$. Let $H_{\alpha,m}$
  be the wall separating $w(A)$ and $ws_i(A)$. Then $w(\alpha_i)$ can
  be written as $w(\alpha_i)=\epsilon (\alpha -m\delta)$, where
  $\epsilon\in \RR$ (in fact $\epsilon=\pm1$ in the untwisted
  case). Furthermore, the following conditions are equivalent:
  \begin{itemize}
  \item[(i)] The $i$-crossing is positive;
  \item[(ii)] $i$ is an (affine) descent of $\cl(w)$;
  \item[(iii)] $\epsilon < 0$.
  \end{itemize}
  Condition (iii) is to be interpreted as $\cl(w)$ maps $\alpha_i$
  (resp. $\coroot_i$) to a negative root (resp. coroot) for $\clW$
  (possibly up to a positive scalar factor for $i=0$ in the twisted
  case).
\end{remark}
\begin{proof}
  Note that $ws_i(A) = w s_i w^{-1} w(A) = s_{w(\alpha_i)} w(A)$, so
  $s_{w(\alpha_i)}=s_{\alpha,m}$. The form for $w(\alpha_i)$
  follows. It remains to prove the equivalence between the three
  conditions.

  (i) $\Longleftrightarrow$ (ii): Let $\rhoc = \frac \ell{n+1}
  (\Lambdac_0+\cdots+\Lambdac_n)$ be the canonical representative of
  the fundamental alcove at level $\ell$: for $i$ in $I$, $\ip{\rhoc,
    \alpha_i} = \frac \ell{n+1}>0$. We compute how the representative
  $w(\rhoc)$ of $w(A)$ is moved in the crossing:
  \begin{equation}
    \begin{split}
      ws_i(\rhoc) - w(\rhoc) &= s_{w(\alpha_i)} w(\rhoc) - w(\rhoc)
      = - \ip{w(\rhoc), w(\alpha_i)} w(\coroot_i)\\
      &= - \ip{\rhoc, \alpha_i} w(\coroot_i)
      = - \frac \ell{n+1} w(\coroot_i)\,.
    \end{split}
  \end{equation}
  The crossing is positive if $\ip{ws_i(\rhoc) -
    w(\rhoc),\alpha}>0$, or equivalently
  \begin{equation}
    0 
    > \ip{ w(\coroot_i), \alpha }
    = \ip{ w(\coroot_i), \frac1\epsilon w(\alpha_i) + m\delta }
    = \frac1\epsilon \ip{ w(\coroot_i), w(\alpha_i) }
    = \frac2\epsilon\,,
  \end{equation}
  that is $\epsilon <0$.

  (i) $\Longleftrightarrow$ (iii): Using
  \eqref{equation.level0heckeaction}, $i$ is a descent of $\cl(w)$ if
  and only if:
  \begin{equation}
    0
    > \ip {w^{-1}(\clrhoc), \alpha_i}
    = \ip {       \clrhoc, w(\alpha_i)}
    = \ip {       \clrhoc, \epsilon (\alpha - m\delta)}
    = \epsilon \ip { \clrhoc, \alpha}\,,
  \end{equation}
  or equivalently $\epsilon<0$.
\end{proof}

By using the interpolation formula $T_i = (q_1+q_2)\pi_i - q_1 s_i$,
the level $0$ actions $\cl$ of the Weyl group $W$ and of the $0$-Hecke
monoid $\langle \pi_i \suchthat i \in I \rangle$ on $\clW$ can be
extended for any $(q_1,q_2)\ne(0,0)$ to a representation $\cl$ of the
affine Hecke algebra $\heckeW{q_1,q_2}$ on $\kclW$.

\emph{Interestingly enough, and this is the central topic of this
  paper, the algebra $\cl(\heckeW{q_1,q_2})=\langle \cl(T_0),\dots,
  \cl(T_n)\rangle$ turns out \emph{not} to be the Hecke algebra
  $\heckeW[\clW]{q_1,q_2}$, except at $q=1$ and certain roots of
  unity.}

\iflongversion

\subsection{Cartan matrix independence}

\label{section.independence}

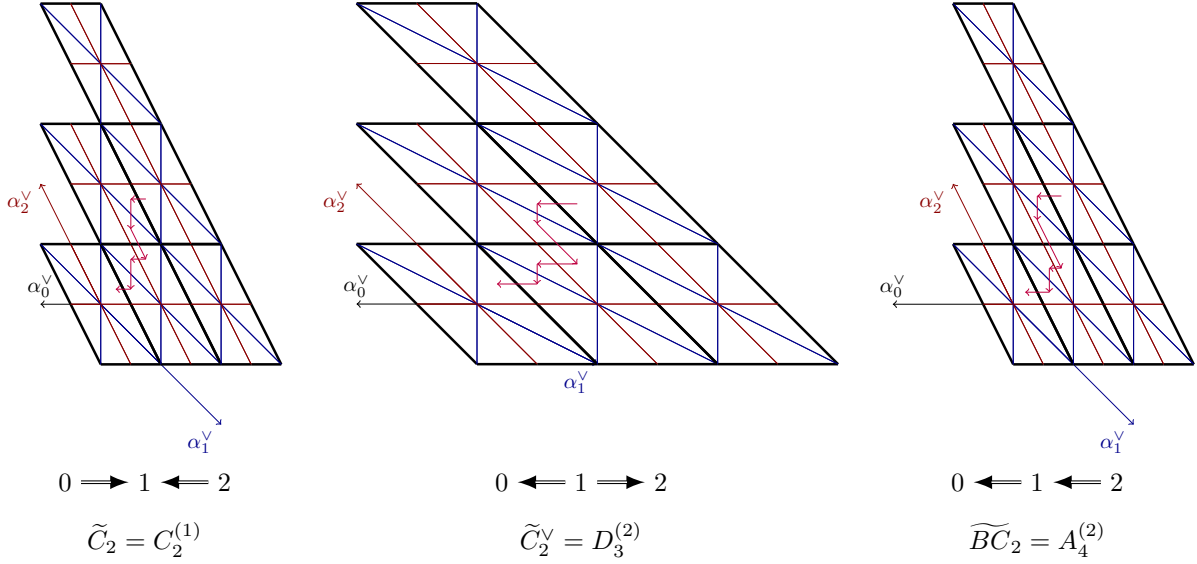
\begin{figure}
  \begin{bigcenter}
  $\begin{array}{ccc}
    \scalebox{.8}{\begin{tikzpicture}[baseline=(current bounding box.east)]
\draw[->,  color = black,] (0.0,0.0) -- (-1.0,0.0) node[at end, auto=right] {$\alpha^\vee_{0}$};
\draw[->,  color = DarkBlue,] (0.0,0.0) -- (2.0,-2.0) node[at end, auto=right] {$\alpha^\vee_{1}$};
\draw[->,  color = DarkRed,] (0.0,0.0) -- (-1.0,2.0) node[at end, auto=left] {$\alpha^\vee_{2}$};
\draw[ color = black, very thick,](0.5,0.0) -- (0.0,1.0);
\draw[ color = DarkBlue, ,](0.0,0.0) -- (0.0,1.0);
\draw[ color = DarkRed, very thin,](0.0,0.0) -- (0.5,0.0);
\draw[ color = black, very thick,](0.5,0.0) -- (1.0,-1.0);
\draw[ color = DarkBlue, ,](0.0,0.0) -- (1.0,-1.0);
\draw[ color = DarkRed, very thin,](0.0,0.0) -- (0.5,0.0);
\draw[ color = black, very thick,](-0.5,0.0) -- (-1.0,1.0);
\draw[ color = DarkBlue, ,](0.0,0.0) -- (-1.0,1.0);
\draw[ color = DarkRed, very thin,](0.0,0.0) -- (-0.5,0.0);
\draw[ color = black, very thick,](-0.5,0.0) -- (0.0,-1.0);
\draw[ color = DarkBlue, ,](0.0,0.0) -- (0.0,-1.0);
\draw[ color = DarkRed, very thin,](0.0,0.0) -- (-0.5,0.0);
\draw[ color = black, very thick,](0.5,-1.0) -- (1.0,-1.0);
\draw[ color = DarkBlue, ,](0.0,0.0) -- (1.0,-1.0);
\draw[ color = DarkRed, very thin,](0.0,0.0) -- (0.5,-1.0);
\draw[ color = black, very thick,](0.5,-1.0) -- (0.0,-1.0);
\draw[ color = DarkBlue, ,](0.0,0.0) -- (0.0,-1.0);
\draw[ color = DarkRed, very thin,](0.0,0.0) -- (0.5,-1.0);
\draw[ color = black, very thick,](-0.5,1.0) -- (0.0,1.0);
\draw[ color = DarkBlue, ,](0.0,0.0) -- (0.0,1.0);
\draw[ color = DarkRed, very thin,](0.0,0.0) -- (-0.5,1.0);
\draw[ color = black, very thick,](-0.5,1.0) -- (-1.0,1.0);
\draw[ color = DarkBlue, ,](0.0,0.0) -- (-1.0,1.0);
\draw[ color = DarkRed, very thin,](0.0,0.0) -- (-0.5,1.0);
\draw[ color = black, very thick,](1.5,0.0) -- (1.0,1.0);
\draw[ color = DarkBlue, ,](1.0,0.0) -- (1.0,1.0);
\draw[ color = DarkRed, very thin,](1.0,0.0) -- (1.5,0.0);
\draw[ color = black, very thick,](1.5,0.0) -- (2.0,-1.0);
\draw[ color = DarkBlue, ,](1.0,0.0) -- (2.0,-1.0);
\draw[ color = DarkRed, very thin,](1.0,0.0) -- (1.5,0.0);
\draw[ color = black, very thick,](0.5,0.0) -- (0.0,1.0);
\draw[ color = DarkBlue, ,](1.0,0.0) -- (0.0,1.0);
\draw[ color = DarkRed, very thin,](1.0,0.0) -- (0.5,0.0);
\draw[ color = black, very thick,](0.5,0.0) -- (1.0,-1.0);
\draw[ color = DarkBlue, ,](1.0,0.0) -- (1.0,-1.0);
\draw[ color = DarkRed, very thin,](1.0,0.0) -- (0.5,0.0);
\draw[ color = black, very thick,](1.5,-1.0) -- (2.0,-1.0);
\draw[ color = DarkBlue, ,](1.0,0.0) -- (2.0,-1.0);
\draw[ color = DarkRed, very thin,](1.0,0.0) -- (1.5,-1.0);
\draw[ color = black, very thick,](1.5,-1.0) -- (1.0,-1.0);
\draw[ color = DarkBlue, ,](1.0,0.0) -- (1.0,-1.0);
\draw[ color = DarkRed, very thin,](1.0,0.0) -- (1.5,-1.0);
\draw[ color = black, very thick,](0.5,1.0) -- (1.0,1.0);
\draw[ color = DarkBlue, ,](1.0,0.0) -- (1.0,1.0);
\draw[ color = DarkRed, very thin,](1.0,0.0) -- (0.5,1.0);
\draw[ color = black, very thick,](0.5,1.0) -- (0.0,1.0);
\draw[ color = DarkBlue, ,](1.0,0.0) -- (0.0,1.0);
\draw[ color = DarkRed, very thin,](1.0,0.0) -- (0.5,1.0);
\draw[ color = black, very thick,](0.5,2.0) -- (0.0,3.0);
\draw[ color = DarkBlue, ,](0.0,2.0) -- (0.0,3.0);
\draw[ color = DarkRed, very thin,](0.0,2.0) -- (0.5,2.0);
\draw[ color = black, very thick,](0.5,2.0) -- (1.0,1.0);
\draw[ color = DarkBlue, ,](0.0,2.0) -- (1.0,1.0);
\draw[ color = DarkRed, very thin,](0.0,2.0) -- (0.5,2.0);
\draw[ color = black, very thick,](-0.5,2.0) -- (-1.0,3.0);
\draw[ color = DarkBlue, ,](0.0,2.0) -- (-1.0,3.0);
\draw[ color = DarkRed, very thin,](0.0,2.0) -- (-0.5,2.0);
\draw[ color = black, very thick,](-0.5,2.0) -- (0.0,1.0);
\draw[ color = DarkBlue, ,](0.0,2.0) -- (0.0,1.0);
\draw[ color = DarkRed, very thin,](0.0,2.0) -- (-0.5,2.0);
\draw[ color = black, very thick,](0.5,1.0) -- (1.0,1.0);
\draw[ color = DarkBlue, ,](0.0,2.0) -- (1.0,1.0);
\draw[ color = DarkRed, very thin,](0.0,2.0) -- (0.5,1.0);
\draw[ color = black, very thick,](0.5,1.0) -- (0.0,1.0);
\draw[ color = DarkBlue, ,](0.0,2.0) -- (0.0,1.0);
\draw[ color = DarkRed, very thin,](0.0,2.0) -- (0.5,1.0);
\draw[ color = black, very thick,](-0.5,3.0) -- (0.0,3.0);
\draw[ color = DarkBlue, ,](0.0,2.0) -- (0.0,3.0);
\draw[ color = DarkRed, very thin,](0.0,2.0) -- (-0.5,3.0);
\draw[ color = black, very thick,](-0.5,3.0) -- (-1.0,3.0);
\draw[ color = DarkBlue, ,](0.0,2.0) -- (-1.0,3.0);
\draw[ color = DarkRed, very thin,](0.0,2.0) -- (-0.5,3.0);
\draw[ color = black, very thick,](2.5,0.0) -- (2.0,1.0);
\draw[ color = DarkBlue, ,](2.0,0.0) -- (2.0,1.0);
\draw[ color = DarkRed, very thin,](2.0,0.0) -- (2.5,0.0);
\draw[ color = black, very thick,](2.5,0.0) -- (3.0,-1.0);
\draw[ color = DarkBlue, ,](2.0,0.0) -- (3.0,-1.0);
\draw[ color = DarkRed, very thin,](2.0,0.0) -- (2.5,0.0);
\draw[ color = black, very thick,](1.5,0.0) -- (1.0,1.0);
\draw[ color = DarkBlue, ,](2.0,0.0) -- (1.0,1.0);
\draw[ color = DarkRed, very thin,](2.0,0.0) -- (1.5,0.0);
\draw[ color = black, very thick,](1.5,0.0) -- (2.0,-1.0);
\draw[ color = DarkBlue, ,](2.0,0.0) -- (2.0,-1.0);
\draw[ color = DarkRed, very thin,](2.0,0.0) -- (1.5,0.0);
\draw[ color = black, very thick,](2.5,-1.0) -- (3.0,-1.0);
\draw[ color = DarkBlue, ,](2.0,0.0) -- (3.0,-1.0);
\draw[ color = DarkRed, very thin,](2.0,0.0) -- (2.5,-1.0);
\draw[ color = black, very thick,](2.5,-1.0) -- (2.0,-1.0);
\draw[ color = DarkBlue, ,](2.0,0.0) -- (2.0,-1.0);
\draw[ color = DarkRed, very thin,](2.0,0.0) -- (2.5,-1.0);
\draw[ color = black, very thick,](1.5,1.0) -- (2.0,1.0);
\draw[ color = DarkBlue, ,](2.0,0.0) -- (2.0,1.0);
\draw[ color = DarkRed, very thin,](2.0,0.0) -- (1.5,1.0);
\draw[ color = black, very thick,](1.5,1.0) -- (1.0,1.0);
\draw[ color = DarkBlue, ,](2.0,0.0) -- (1.0,1.0);
\draw[ color = DarkRed, very thin,](2.0,0.0) -- (1.5,1.0);
\draw[ color = black, very thick,](1.5,2.0) -- (1.0,3.0);
\draw[ color = DarkBlue, ,](1.0,2.0) -- (1.0,3.0);
\draw[ color = DarkRed, very thin,](1.0,2.0) -- (1.5,2.0);
\draw[ color = black, very thick,](1.5,2.0) -- (2.0,1.0);
\draw[ color = DarkBlue, ,](1.0,2.0) -- (2.0,1.0);
\draw[ color = DarkRed, very thin,](1.0,2.0) -- (1.5,2.0);
\draw[ color = black, very thick,](0.5,2.0) -- (0.0,3.0);
\draw[ color = DarkBlue, ,](1.0,2.0) -- (0.0,3.0);
\draw[ color = DarkRed, very thin,](1.0,2.0) -- (0.5,2.0);
\draw[ color = black, very thick,](0.5,2.0) -- (1.0,1.0);
\draw[ color = DarkBlue, ,](1.0,2.0) -- (1.0,1.0);
\draw[ color = DarkRed, very thin,](1.0,2.0) -- (0.5,2.0);
\draw[ color = black, very thick,](1.5,1.0) -- (2.0,1.0);
\draw[ color = DarkBlue, ,](1.0,2.0) -- (2.0,1.0);
\draw[ color = DarkRed, very thin,](1.0,2.0) -- (1.5,1.0);
\draw[ color = black, very thick,](1.5,1.0) -- (1.0,1.0);
\draw[ color = DarkBlue, ,](1.0,2.0) -- (1.0,1.0);
\draw[ color = DarkRed, very thin,](1.0,2.0) -- (1.5,1.0);
\draw[ color = black, very thick,](0.5,3.0) -- (1.0,3.0);
\draw[ color = DarkBlue, ,](1.0,2.0) -- (1.0,3.0);
\draw[ color = DarkRed, very thin,](1.0,2.0) -- (0.5,3.0);
\draw[ color = black, very thick,](0.5,3.0) -- (0.0,3.0);
\draw[ color = DarkBlue, ,](1.0,2.0) -- (0.0,3.0);
\draw[ color = DarkRed, very thin,](1.0,2.0) -- (0.5,3.0);
\draw[ color = black, very thick,](0.5,4.0) -- (0.0,5.0);
\draw[ color = DarkBlue, ,](0.0,4.0) -- (0.0,5.0);
\draw[ color = DarkRed, very thin,](0.0,4.0) -- (0.5,4.0);
\draw[ color = black, very thick,](0.5,4.0) -- (1.0,3.0);
\draw[ color = DarkBlue, ,](0.0,4.0) -- (1.0,3.0);
\draw[ color = DarkRed, very thin,](0.0,4.0) -- (0.5,4.0);
\draw[ color = black, very thick,](-0.5,4.0) -- (-1.0,5.0);
\draw[ color = DarkBlue, ,](0.0,4.0) -- (-1.0,5.0);
\draw[ color = DarkRed, very thin,](0.0,4.0) -- (-0.5,4.0);
\draw[ color = black, very thick,](-0.5,4.0) -- (0.0,3.0);
\draw[ color = DarkBlue, ,](0.0,4.0) -- (0.0,3.0);
\draw[ color = DarkRed, very thin,](0.0,4.0) -- (-0.5,4.0);
\draw[ color = black, very thick,](0.5,3.0) -- (1.0,3.0);
\draw[ color = DarkBlue, ,](0.0,4.0) -- (1.0,3.0);
\draw[ color = DarkRed, very thin,](0.0,4.0) -- (0.5,3.0);
\draw[ color = black, very thick,](0.5,3.0) -- (0.0,3.0);
\draw[ color = DarkBlue, ,](0.0,4.0) -- (0.0,3.0);
\draw[ color = DarkRed, very thin,](0.0,4.0) -- (0.5,3.0);
\draw[ color = black, very thick,](-0.5,5.0) -- (0.0,5.0);
\draw[ color = DarkBlue, ,](0.0,4.0) -- (0.0,5.0);
\draw[ color = DarkRed, very thin,](0.0,4.0) -- (-0.5,5.0);
\draw[ color = black, very thick,](-0.5,5.0) -- (-1.0,5.0);
\draw[ color = DarkBlue, ,](0.0,4.0) -- (-1.0,5.0);
\draw[ color = DarkRed, very thin,](0.0,4.0) -- (-0.5,5.0);
\draw[->,  color = purple,] (0.75,1.75) -- (0.5,1.75);
\draw[->,  color = purple,] (0.5,1.75) -- (0.5,1.25);
\draw[->,  color = purple,] (0.5,1.25) -- (0.75,0.75);
\draw[->,  color = purple,] (0.75,0.75) -- (0.5,0.75);
\draw[->,  color = purple,] (0.5,0.75) -- (0.5,0.25);
\draw[->,  color = purple,] (0.5,0.25) -- (0.25,0.25);
\end{tikzpicture}}&
    \raisebox{2.65ex}{\scalebox{.8}{\begin{tikzpicture}[baseline=(current bounding box.east)]
\draw[->,  color = black,] (0.0,0.0) -- (-2.0,0.0) node[at end, auto=right] {$\alpha^\vee_{0}$};
\draw[->,  color = DarkBlue,] (0.0,0.0) -- (2.0,-1.0) node[at end, auto=right] {$\alpha^\vee_{1}$};
\draw[->,  color = DarkRed,] (0.0,0.0) -- (-2.0,2.0) node[at end, auto=left] {$\alpha^\vee_{2}$};
\draw[ color = black, very thick,](1.0,0.0) -- (0.0,1.0);
\draw[ color = DarkBlue, ,](0.0,0.0) -- (0.0,1.0);
\draw[ color = DarkRed, very thin,](0.0,0.0) -- (1.0,0.0);
\draw[ color = black, very thick,](1.0,0.0) -- (2.0,-1.0);
\draw[ color = DarkBlue, ,](0.0,0.0) -- (2.0,-1.0);
\draw[ color = DarkRed, very thin,](0.0,0.0) -- (1.0,0.0);
\draw[ color = black, very thick,](-1.0,0.0) -- (-2.0,1.0);
\draw[ color = DarkBlue, ,](0.0,0.0) -- (-2.0,1.0);
\draw[ color = DarkRed, very thin,](0.0,0.0) -- (-1.0,0.0);
\draw[ color = black, very thick,](-1.0,0.0) -- (0.0,-1.0);
\draw[ color = DarkBlue, ,](0.0,0.0) -- (0.0,-1.0);
\draw[ color = DarkRed, very thin,](0.0,0.0) -- (-1.0,0.0);
\draw[ color = black, very thick,](1.0,-1.0) -- (2.0,-1.0);
\draw[ color = DarkBlue, ,](0.0,0.0) -- (2.0,-1.0);
\draw[ color = DarkRed, very thin,](0.0,0.0) -- (1.0,-1.0);
\draw[ color = black, very thick,](1.0,-1.0) -- (0.0,-1.0);
\draw[ color = DarkBlue, ,](0.0,0.0) -- (0.0,-1.0);
\draw[ color = DarkRed, very thin,](0.0,0.0) -- (1.0,-1.0);
\draw[ color = black, very thick,](-1.0,1.0) -- (0.0,1.0);
\draw[ color = DarkBlue, ,](0.0,0.0) -- (0.0,1.0);
\draw[ color = DarkRed, very thin,](0.0,0.0) -- (-1.0,1.0);
\draw[ color = black, very thick,](-1.0,1.0) -- (-2.0,1.0);
\draw[ color = DarkBlue, ,](0.0,0.0) -- (-2.0,1.0);
\draw[ color = DarkRed, very thin,](0.0,0.0) -- (-1.0,1.0);
\draw[ color = black, very thick,](3.0,0.0) -- (2.0,1.0);
\draw[ color = DarkBlue, ,](2.0,0.0) -- (2.0,1.0);
\draw[ color = DarkRed, very thin,](2.0,0.0) -- (3.0,0.0);
\draw[ color = black, very thick,](3.0,0.0) -- (4.0,-1.0);
\draw[ color = DarkBlue, ,](2.0,0.0) -- (4.0,-1.0);
\draw[ color = DarkRed, very thin,](2.0,0.0) -- (3.0,0.0);
\draw[ color = black, very thick,](1.0,0.0) -- (0.0,1.0);
\draw[ color = DarkBlue, ,](2.0,0.0) -- (0.0,1.0);
\draw[ color = DarkRed, very thin,](2.0,0.0) -- (1.0,0.0);
\draw[ color = black, very thick,](1.0,0.0) -- (2.0,-1.0);
\draw[ color = DarkBlue, ,](2.0,0.0) -- (2.0,-1.0);
\draw[ color = DarkRed, very thin,](2.0,0.0) -- (1.0,0.0);
\draw[ color = black, very thick,](3.0,-1.0) -- (4.0,-1.0);
\draw[ color = DarkBlue, ,](2.0,0.0) -- (4.0,-1.0);
\draw[ color = DarkRed, very thin,](2.0,0.0) -- (3.0,-1.0);
\draw[ color = black, very thick,](3.0,-1.0) -- (2.0,-1.0);
\draw[ color = DarkBlue, ,](2.0,0.0) -- (2.0,-1.0);
\draw[ color = DarkRed, very thin,](2.0,0.0) -- (3.0,-1.0);
\draw[ color = black, very thick,](1.0,1.0) -- (2.0,1.0);
\draw[ color = DarkBlue, ,](2.0,0.0) -- (2.0,1.0);
\draw[ color = DarkRed, very thin,](2.0,0.0) -- (1.0,1.0);
\draw[ color = black, very thick,](1.0,1.0) -- (0.0,1.0);
\draw[ color = DarkBlue, ,](2.0,0.0) -- (0.0,1.0);
\draw[ color = DarkRed, very thin,](2.0,0.0) -- (1.0,1.0);
\draw[ color = black, very thick,](1.0,2.0) -- (0.0,3.0);
\draw[ color = DarkBlue, ,](0.0,2.0) -- (0.0,3.0);
\draw[ color = DarkRed, very thin,](0.0,2.0) -- (1.0,2.0);
\draw[ color = black, very thick,](1.0,2.0) -- (2.0,1.0);
\draw[ color = DarkBlue, ,](0.0,2.0) -- (2.0,1.0);
\draw[ color = DarkRed, very thin,](0.0,2.0) -- (1.0,2.0);
\draw[ color = black, very thick,](-1.0,2.0) -- (-2.0,3.0);
\draw[ color = DarkBlue, ,](0.0,2.0) -- (-2.0,3.0);
\draw[ color = DarkRed, very thin,](0.0,2.0) -- (-1.0,2.0);
\draw[ color = black, very thick,](-1.0,2.0) -- (0.0,1.0);
\draw[ color = DarkBlue, ,](0.0,2.0) -- (0.0,1.0);
\draw[ color = DarkRed, very thin,](0.0,2.0) -- (-1.0,2.0);
\draw[ color = black, very thick,](1.0,1.0) -- (2.0,1.0);
\draw[ color = DarkBlue, ,](0.0,2.0) -- (2.0,1.0);
\draw[ color = DarkRed, very thin,](0.0,2.0) -- (1.0,1.0);
\draw[ color = black, very thick,](1.0,1.0) -- (0.0,1.0);
\draw[ color = DarkBlue, ,](0.0,2.0) -- (0.0,1.0);
\draw[ color = DarkRed, very thin,](0.0,2.0) -- (1.0,1.0);
\draw[ color = black, very thick,](-1.0,3.0) -- (0.0,3.0);
\draw[ color = DarkBlue, ,](0.0,2.0) -- (0.0,3.0);
\draw[ color = DarkRed, very thin,](0.0,2.0) -- (-1.0,3.0);
\draw[ color = black, very thick,](-1.0,3.0) -- (-2.0,3.0);
\draw[ color = DarkBlue, ,](0.0,2.0) -- (-2.0,3.0);
\draw[ color = DarkRed, very thin,](0.0,2.0) -- (-1.0,3.0);
\draw[ color = black, very thick,](5.0,0.0) -- (4.0,1.0);
\draw[ color = DarkBlue, ,](4.0,0.0) -- (4.0,1.0);
\draw[ color = DarkRed, very thin,](4.0,0.0) -- (5.0,0.0);
\draw[ color = black, very thick,](5.0,0.0) -- (6.0,-1.0);
\draw[ color = DarkBlue, ,](4.0,0.0) -- (6.0,-1.0);
\draw[ color = DarkRed, very thin,](4.0,0.0) -- (5.0,0.0);
\draw[ color = black, very thick,](3.0,0.0) -- (2.0,1.0);
\draw[ color = DarkBlue, ,](4.0,0.0) -- (2.0,1.0);
\draw[ color = DarkRed, very thin,](4.0,0.0) -- (3.0,0.0);
\draw[ color = black, very thick,](3.0,0.0) -- (4.0,-1.0);
\draw[ color = DarkBlue, ,](4.0,0.0) -- (4.0,-1.0);
\draw[ color = DarkRed, very thin,](4.0,0.0) -- (3.0,0.0);
\draw[ color = black, very thick,](5.0,-1.0) -- (6.0,-1.0);
\draw[ color = DarkBlue, ,](4.0,0.0) -- (6.0,-1.0);
\draw[ color = DarkRed, very thin,](4.0,0.0) -- (5.0,-1.0);
\draw[ color = black, very thick,](5.0,-1.0) -- (4.0,-1.0);
\draw[ color = DarkBlue, ,](4.0,0.0) -- (4.0,-1.0);
\draw[ color = DarkRed, very thin,](4.0,0.0) -- (5.0,-1.0);
\draw[ color = black, very thick,](3.0,1.0) -- (4.0,1.0);
\draw[ color = DarkBlue, ,](4.0,0.0) -- (4.0,1.0);
\draw[ color = DarkRed, very thin,](4.0,0.0) -- (3.0,1.0);
\draw[ color = black, very thick,](3.0,1.0) -- (2.0,1.0);
\draw[ color = DarkBlue, ,](4.0,0.0) -- (2.0,1.0);
\draw[ color = DarkRed, very thin,](4.0,0.0) -- (3.0,1.0);
\draw[ color = black, very thick,](3.0,2.0) -- (2.0,3.0);
\draw[ color = DarkBlue, ,](2.0,2.0) -- (2.0,3.0);
\draw[ color = DarkRed, very thin,](2.0,2.0) -- (3.0,2.0);
\draw[ color = black, very thick,](3.0,2.0) -- (4.0,1.0);
\draw[ color = DarkBlue, ,](2.0,2.0) -- (4.0,1.0);
\draw[ color = DarkRed, very thin,](2.0,2.0) -- (3.0,2.0);
\draw[ color = black, very thick,](1.0,2.0) -- (0.0,3.0);
\draw[ color = DarkBlue, ,](2.0,2.0) -- (0.0,3.0);
\draw[ color = DarkRed, very thin,](2.0,2.0) -- (1.0,2.0);
\draw[ color = black, very thick,](1.0,2.0) -- (2.0,1.0);
\draw[ color = DarkBlue, ,](2.0,2.0) -- (2.0,1.0);
\draw[ color = DarkRed, very thin,](2.0,2.0) -- (1.0,2.0);
\draw[ color = black, very thick,](3.0,1.0) -- (4.0,1.0);
\draw[ color = DarkBlue, ,](2.0,2.0) -- (4.0,1.0);
\draw[ color = DarkRed, very thin,](2.0,2.0) -- (3.0,1.0);
\draw[ color = black, very thick,](3.0,1.0) -- (2.0,1.0);
\draw[ color = DarkBlue, ,](2.0,2.0) -- (2.0,1.0);
\draw[ color = DarkRed, very thin,](2.0,2.0) -- (3.0,1.0);
\draw[ color = black, very thick,](1.0,3.0) -- (2.0,3.0);
\draw[ color = DarkBlue, ,](2.0,2.0) -- (2.0,3.0);
\draw[ color = DarkRed, very thin,](2.0,2.0) -- (1.0,3.0);
\draw[ color = black, very thick,](1.0,3.0) -- (0.0,3.0);
\draw[ color = DarkBlue, ,](2.0,2.0) -- (0.0,3.0);
\draw[ color = DarkRed, very thin,](2.0,2.0) -- (1.0,3.0);
\draw[ color = black, very thick,](1.0,4.0) -- (0.0,5.0);
\draw[ color = DarkBlue, ,](0.0,4.0) -- (0.0,5.0);
\draw[ color = DarkRed, very thin,](0.0,4.0) -- (1.0,4.0);
\draw[ color = black, very thick,](1.0,4.0) -- (2.0,3.0);
\draw[ color = DarkBlue, ,](0.0,4.0) -- (2.0,3.0);
\draw[ color = DarkRed, very thin,](0.0,4.0) -- (1.0,4.0);
\draw[ color = black, very thick,](-1.0,4.0) -- (-2.0,5.0);
\draw[ color = DarkBlue, ,](0.0,4.0) -- (-2.0,5.0);
\draw[ color = DarkRed, very thin,](0.0,4.0) -- (-1.0,4.0);
\draw[ color = black, very thick,](-1.0,4.0) -- (0.0,3.0);
\draw[ color = DarkBlue, ,](0.0,4.0) -- (0.0,3.0);
\draw[ color = DarkRed, very thin,](0.0,4.0) -- (-1.0,4.0);
\draw[ color = black, very thick,](1.0,3.0) -- (2.0,3.0);
\draw[ color = DarkBlue, ,](0.0,4.0) -- (2.0,3.0);
\draw[ color = DarkRed, very thin,](0.0,4.0) -- (1.0,3.0);
\draw[ color = black, very thick,](1.0,3.0) -- (0.0,3.0);
\draw[ color = DarkBlue, ,](0.0,4.0) -- (0.0,3.0);
\draw[ color = DarkRed, very thin,](0.0,4.0) -- (1.0,3.0);
\draw[ color = black, very thick,](-1.0,5.0) -- (0.0,5.0);
\draw[ color = DarkBlue, ,](0.0,4.0) -- (0.0,5.0);
\draw[ color = DarkRed, very thin,](0.0,4.0) -- (-1.0,5.0);
\draw[ color = black, very thick,](-1.0,5.0) -- (-2.0,5.0);
\draw[ color = DarkBlue, ,](0.0,4.0) -- (-2.0,5.0);
\draw[ color = DarkRed, very thin,](0.0,4.0) -- (-1.0,5.0);
\draw[->,  color = purple,] (1.666666667,1.666666667) -- (1.0,1.666666667);
\draw[->,  color = purple,] (1.0,1.666666667) -- (1.0,1.333333333);
\draw[->,  color = purple,] (1.0,1.333333333) -- (1.666666667,0.6666666667);
\draw[->,  color = purple,] (1.666666667,0.6666666667) -- (1.0,0.6666666667);
\draw[->,  color = purple,] (1.0,0.6666666667) -- (1.0,0.3333333333);
\draw[->,  color = purple,] (1.0,0.3333333333) -- (0.3333333333,0.3333333333);
\end{tikzpicture}}} &
    \scalebox{.8}{\begin{tikzpicture}[baseline=(current bounding box.east)]
\draw[->,  color = black,] (0.0,0.0) -- (-2.0,0.0) node[at end, auto=right] {$\alpha^\vee_{0}$};
\draw[->,  color = DarkBlue,] (0.0,0.0) -- (2.0,-2.0) node[at end, auto=right] {$\alpha^\vee_{1}$};
\draw[->,  color = DarkRed,] (0.0,0.0) -- (-1.0,2.0) node[at end, auto=left] {$\alpha^\vee_{2}$};
\draw[ color = black, very thick,](0.5,0.0) -- (0.0,1.0);
\draw[ color = DarkBlue, ,](0.0,0.0) -- (0.0,1.0);
\draw[ color = DarkRed, very thin,](0.0,0.0) -- (0.5,0.0);
\draw[ color = black, very thick,](0.5,0.0) -- (1.0,-1.0);
\draw[ color = DarkBlue, ,](0.0,0.0) -- (1.0,-1.0);
\draw[ color = DarkRed, very thin,](0.0,0.0) -- (0.5,0.0);
\draw[ color = black, very thick,](-0.5,0.0) -- (-1.0,1.0);
\draw[ color = DarkBlue, ,](0.0,0.0) -- (-1.0,1.0);
\draw[ color = DarkRed, very thin,](0.0,0.0) -- (-0.5,0.0);
\draw[ color = black, very thick,](-0.5,0.0) -- (0.0,-1.0);
\draw[ color = DarkBlue, ,](0.0,0.0) -- (0.0,-1.0);
\draw[ color = DarkRed, very thin,](0.0,0.0) -- (-0.5,0.0);
\draw[ color = black, very thick,](0.5,-1.0) -- (1.0,-1.0);
\draw[ color = DarkBlue, ,](0.0,0.0) -- (1.0,-1.0);
\draw[ color = DarkRed, very thin,](0.0,0.0) -- (0.5,-1.0);
\draw[ color = black, very thick,](0.5,-1.0) -- (0.0,-1.0);
\draw[ color = DarkBlue, ,](0.0,0.0) -- (0.0,-1.0);
\draw[ color = DarkRed, very thin,](0.0,0.0) -- (0.5,-1.0);
\draw[ color = black, very thick,](-0.5,1.0) -- (0.0,1.0);
\draw[ color = DarkBlue, ,](0.0,0.0) -- (0.0,1.0);
\draw[ color = DarkRed, very thin,](0.0,0.0) -- (-0.5,1.0);
\draw[ color = black, very thick,](-0.5,1.0) -- (-1.0,1.0);
\draw[ color = DarkBlue, ,](0.0,0.0) -- (-1.0,1.0);
\draw[ color = DarkRed, very thin,](0.0,0.0) -- (-0.5,1.0);
\draw[ color = black, very thick,](1.5,0.0) -- (1.0,1.0);
\draw[ color = DarkBlue, ,](1.0,0.0) -- (1.0,1.0);
\draw[ color = DarkRed, very thin,](1.0,0.0) -- (1.5,0.0);
\draw[ color = black, very thick,](1.5,0.0) -- (2.0,-1.0);
\draw[ color = DarkBlue, ,](1.0,0.0) -- (2.0,-1.0);
\draw[ color = DarkRed, very thin,](1.0,0.0) -- (1.5,0.0);
\draw[ color = black, very thick,](0.5,0.0) -- (0.0,1.0);
\draw[ color = DarkBlue, ,](1.0,0.0) -- (0.0,1.0);
\draw[ color = DarkRed, very thin,](1.0,0.0) -- (0.5,0.0);
\draw[ color = black, very thick,](0.5,0.0) -- (1.0,-1.0);
\draw[ color = DarkBlue, ,](1.0,0.0) -- (1.0,-1.0);
\draw[ color = DarkRed, very thin,](1.0,0.0) -- (0.5,0.0);
\draw[ color = black, very thick,](1.5,-1.0) -- (2.0,-1.0);
\draw[ color = DarkBlue, ,](1.0,0.0) -- (2.0,-1.0);
\draw[ color = DarkRed, very thin,](1.0,0.0) -- (1.5,-1.0);
\draw[ color = black, very thick,](1.5,-1.0) -- (1.0,-1.0);
\draw[ color = DarkBlue, ,](1.0,0.0) -- (1.0,-1.0);
\draw[ color = DarkRed, very thin,](1.0,0.0) -- (1.5,-1.0);
\draw[ color = black, very thick,](0.5,1.0) -- (1.0,1.0);
\draw[ color = DarkBlue, ,](1.0,0.0) -- (1.0,1.0);
\draw[ color = DarkRed, very thin,](1.0,0.0) -- (0.5,1.0);
\draw[ color = black, very thick,](0.5,1.0) -- (0.0,1.0);
\draw[ color = DarkBlue, ,](1.0,0.0) -- (0.0,1.0);
\draw[ color = DarkRed, very thin,](1.0,0.0) -- (0.5,1.0);
\draw[ color = black, very thick,](0.5,2.0) -- (0.0,3.0);
\draw[ color = DarkBlue, ,](0.0,2.0) -- (0.0,3.0);
\draw[ color = DarkRed, very thin,](0.0,2.0) -- (0.5,2.0);
\draw[ color = black, very thick,](0.5,2.0) -- (1.0,1.0);
\draw[ color = DarkBlue, ,](0.0,2.0) -- (1.0,1.0);
\draw[ color = DarkRed, very thin,](0.0,2.0) -- (0.5,2.0);
\draw[ color = black, very thick,](-0.5,2.0) -- (-1.0,3.0);
\draw[ color = DarkBlue, ,](0.0,2.0) -- (-1.0,3.0);
\draw[ color = DarkRed, very thin,](0.0,2.0) -- (-0.5,2.0);
\draw[ color = black, very thick,](-0.5,2.0) -- (0.0,1.0);
\draw[ color = DarkBlue, ,](0.0,2.0) -- (0.0,1.0);
\draw[ color = DarkRed, very thin,](0.0,2.0) -- (-0.5,2.0);
\draw[ color = black, very thick,](0.5,1.0) -- (1.0,1.0);
\draw[ color = DarkBlue, ,](0.0,2.0) -- (1.0,1.0);
\draw[ color = DarkRed, very thin,](0.0,2.0) -- (0.5,1.0);
\draw[ color = black, very thick,](0.5,1.0) -- (0.0,1.0);
\draw[ color = DarkBlue, ,](0.0,2.0) -- (0.0,1.0);
\draw[ color = DarkRed, very thin,](0.0,2.0) -- (0.5,1.0);
\draw[ color = black, very thick,](-0.5,3.0) -- (0.0,3.0);
\draw[ color = DarkBlue, ,](0.0,2.0) -- (0.0,3.0);
\draw[ color = DarkRed, very thin,](0.0,2.0) -- (-0.5,3.0);
\draw[ color = black, very thick,](-0.5,3.0) -- (-1.0,3.0);
\draw[ color = DarkBlue, ,](0.0,2.0) -- (-1.0,3.0);
\draw[ color = DarkRed, very thin,](0.0,2.0) -- (-0.5,3.0);
\draw[ color = black, very thick,](2.5,0.0) -- (2.0,1.0);
\draw[ color = DarkBlue, ,](2.0,0.0) -- (2.0,1.0);
\draw[ color = DarkRed, very thin,](2.0,0.0) -- (2.5,0.0);
\draw[ color = black, very thick,](2.5,0.0) -- (3.0,-1.0);
\draw[ color = DarkBlue, ,](2.0,0.0) -- (3.0,-1.0);
\draw[ color = DarkRed, very thin,](2.0,0.0) -- (2.5,0.0);
\draw[ color = black, very thick,](1.5,0.0) -- (1.0,1.0);
\draw[ color = DarkBlue, ,](2.0,0.0) -- (1.0,1.0);
\draw[ color = DarkRed, very thin,](2.0,0.0) -- (1.5,0.0);
\draw[ color = black, very thick,](1.5,0.0) -- (2.0,-1.0);
\draw[ color = DarkBlue, ,](2.0,0.0) -- (2.0,-1.0);
\draw[ color = DarkRed, very thin,](2.0,0.0) -- (1.5,0.0);
\draw[ color = black, very thick,](2.5,-1.0) -- (3.0,-1.0);
\draw[ color = DarkBlue, ,](2.0,0.0) -- (3.0,-1.0);
\draw[ color = DarkRed, very thin,](2.0,0.0) -- (2.5,-1.0);
\draw[ color = black, very thick,](2.5,-1.0) -- (2.0,-1.0);
\draw[ color = DarkBlue, ,](2.0,0.0) -- (2.0,-1.0);
\draw[ color = DarkRed, very thin,](2.0,0.0) -- (2.5,-1.0);
\draw[ color = black, very thick,](1.5,1.0) -- (2.0,1.0);
\draw[ color = DarkBlue, ,](2.0,0.0) -- (2.0,1.0);
\draw[ color = DarkRed, very thin,](2.0,0.0) -- (1.5,1.0);
\draw[ color = black, very thick,](1.5,1.0) -- (1.0,1.0);
\draw[ color = DarkBlue, ,](2.0,0.0) -- (1.0,1.0);
\draw[ color = DarkRed, very thin,](2.0,0.0) -- (1.5,1.0);
\draw[ color = black, very thick,](1.5,2.0) -- (1.0,3.0);
\draw[ color = DarkBlue, ,](1.0,2.0) -- (1.0,3.0);
\draw[ color = DarkRed, very thin,](1.0,2.0) -- (1.5,2.0);
\draw[ color = black, very thick,](1.5,2.0) -- (2.0,1.0);
\draw[ color = DarkBlue, ,](1.0,2.0) -- (2.0,1.0);
\draw[ color = DarkRed, very thin,](1.0,2.0) -- (1.5,2.0);
\draw[ color = black, very thick,](0.5,2.0) -- (0.0,3.0);
\draw[ color = DarkBlue, ,](1.0,2.0) -- (0.0,3.0);
\draw[ color = DarkRed, very thin,](1.0,2.0) -- (0.5,2.0);
\draw[ color = black, very thick,](0.5,2.0) -- (1.0,1.0);
\draw[ color = DarkBlue, ,](1.0,2.0) -- (1.0,1.0);
\draw[ color = DarkRed, very thin,](1.0,2.0) -- (0.5,2.0);
\draw[ color = black, very thick,](1.5,1.0) -- (2.0,1.0);
\draw[ color = DarkBlue, ,](1.0,2.0) -- (2.0,1.0);
\draw[ color = DarkRed, very thin,](1.0,2.0) -- (1.5,1.0);
\draw[ color = black, very thick,](1.5,1.0) -- (1.0,1.0);
\draw[ color = DarkBlue, ,](1.0,2.0) -- (1.0,1.0);
\draw[ color = DarkRed, very thin,](1.0,2.0) -- (1.5,1.0);
\draw[ color = black, very thick,](0.5,3.0) -- (1.0,3.0);
\draw[ color = DarkBlue, ,](1.0,2.0) -- (1.0,3.0);
\draw[ color = DarkRed, very thin,](1.0,2.0) -- (0.5,3.0);
\draw[ color = black, very thick,](0.5,3.0) -- (0.0,3.0);
\draw[ color = DarkBlue, ,](1.0,2.0) -- (0.0,3.0);
\draw[ color = DarkRed, very thin,](1.0,2.0) -- (0.5,3.0);
\draw[ color = black, very thick,](0.5,4.0) -- (0.0,5.0);
\draw[ color = DarkBlue, ,](0.0,4.0) -- (0.0,5.0);
\draw[ color = DarkRed, very thin,](0.0,4.0) -- (0.5,4.0);
\draw[ color = black, very thick,](0.5,4.0) -- (1.0,3.0);
\draw[ color = DarkBlue, ,](0.0,4.0) -- (1.0,3.0);
\draw[ color = DarkRed, very thin,](0.0,4.0) -- (0.5,4.0);
\draw[ color = black, very thick,](-0.5,4.0) -- (-1.0,5.0);
\draw[ color = DarkBlue, ,](0.0,4.0) -- (-1.0,5.0);
\draw[ color = DarkRed, very thin,](0.0,4.0) -- (-0.5,4.0);
\draw[ color = black, very thick,](-0.5,4.0) -- (0.0,3.0);
\draw[ color = DarkBlue, ,](0.0,4.0) -- (0.0,3.0);
\draw[ color = DarkRed, very thin,](0.0,4.0) -- (-0.5,4.0);
\draw[ color = black, very thick,](0.5,3.0) -- (1.0,3.0);
\draw[ color = DarkBlue, ,](0.0,4.0) -- (1.0,3.0);
\draw[ color = DarkRed, very thin,](0.0,4.0) -- (0.5,3.0);
\draw[ color = black, very thick,](0.5,3.0) -- (0.0,3.0);
\draw[ color = DarkBlue, ,](0.0,4.0) -- (0.0,3.0);
\draw[ color = DarkRed, very thin,](0.0,4.0) -- (0.5,3.0);
\draw[ color = black, very thick,](-0.5,5.0) -- (0.0,5.0);
\draw[ color = DarkBlue, ,](0.0,4.0) -- (0.0,5.0);
\draw[ color = DarkRed, very thin,](0.0,4.0) -- (-0.5,5.0);
\draw[ color = black, very thick,](-0.5,5.0) -- (-1.0,5.0);
\draw[ color = DarkBlue, ,](0.0,4.0) -- (-1.0,5.0);
\draw[ color = DarkRed, very thin,](0.0,4.0) -- (-0.5,5.0);
\draw[->,  color = purple,] (0.8,1.8) -- (0.4,1.8);
\draw[->,  color = purple,] (0.4,1.8) -- (0.4,1.4);
\draw[->,  color = purple,] (0.4,1.4) -- (0.8,0.6);
\draw[->,  color = purple,] (0.8,0.6) -- (0.6,0.6);
\draw[->,  color = purple,] (0.6,0.6) -- (0.6,0.2);
\draw[->,  color = purple,] (0.6,0.2) -- (0.2,0.2);
\end{tikzpicture}} \\
    \dynkinCIIa & \dynkinDIIb & \dynkinAIVb\\[4mm]
    \widetilde C_2=C_2^{(1)}  & \widetilde C_2^\vee=D_3^{(2)} & \widetilde{BC}_2=A_4^{(2)}
  \end{array}$
  \end{bigcenter}
  \caption{The alcove pictures and Dynkin diagrams for the three realizations of the
    Coxeter group $\widetilde C_2 = C_2^{(1)}$ as an affine Weyl group, drawn
    in the coweight lattice. The sample alcove walk is the same as in
    Figure~\ref{figure.antisorting}. Notice that the pictures are
    identical up to a diagonal deformation.}
  \label{figure.rescaling}
\end{figure}

In this subsection, we show that the geometric picture is independent
of the chosen generalized Cartan matrix of $W$ (see Figure~\ref{figure.rescaling}).
In other words, this paper is
really about Coxeter groups which happen to have a realization as
affine Weyl groups, and not about Weyl groups. In particular, one
could always assume without loss of generality that the chosen
geometric representation comes from a realization of $\W$ as an
untwisted affine Weyl group.

Let $\W$ be any Coxeter group, and $M$ and $M'$ be two symmetrizable
generalized Cartan matrices for $\W$, and $D=(d_i)_{i\in I}$ be the
diagonal matrix such that $M'=D M D^{-1}$. We denote by
$\coweightspace$ and $\coweightspace'$ the corresponding geometric
realizations of $\W$, by $\clcoweightspace$ and ${\clcoweightspace}'$
the linear span of the coroots, etc.
Consider the isomorphism $d:{\rootspace}'\to\rootspace$ determined by
$d(\alpha'_i) := \frac1{d_i}\alpha_i$. Further fix an isomorphism
$\dc:{\coweightspace}'\to\coweightspace$ such that $\dc({\coroot_i}') :=
d_i\coroot_i$ ($\dc$ is a well-defined and unique isomorphism from
${\clcoweightspace}'$ to $\clcoweightspace$: given the relation
between $M'$ and $M$ linear relations between the ${\coroot_i}'$'s are
mapped to linear relations between the $\coroot_i$'s, and one can
extend it to $\coweightspace$).

Straightforward computations show that $\ip{\dc(\xc), d(y)} =
\ip{\xc,y}$, $s_i(\dc(\xc))=\dc(s_i(\xc))$ and
$s_i(d(y))=d(s_i(y))$, so that $\dc$ and $d$ are $\W$-morphisms.  It
follows that a root $\alpha'=w\alpha'_i$ of $\W$ in ${\rootspace}'$ is
mapped by $\dc$ to a positive scalar multiple of $\alpha=w\alpha_i$ in
$\rootspace$.  So, $\dc$ preserves the hyperplane $H_\alpha$ and the
half spaces $H^+_\alpha$ and $H^-_\alpha$. Therefore $\dc$ preserves
chambers and in particular the fundamental one, the Tits cone, the
bijection between chambers and elements of $\W$; furthermore $\dc$ is a
morphism for the action of the $\pi_i$'s.

Assume now that $\W$ can be realized as an affine Weyl group.  The
action of $\W$ on the level $0$-hyperplanes are isomorphic, and thus
$\clW'$ and $\clW$ form the same quotient of $\W$. Also, the level $0$
action of $\W$ and of the $0$-Hecke monoid on $\clW$, and therefore
the representation of the $q_1,q_2$-affine Hecke algebra on $\kclW$
match. The set of translations (elements of $\W$ acting trivially at
level $0$) are the same, and for $\lc$ in the coroot lattice of $\clW$
we get identical expressions for $t_\lc$ in terms of the $s_i$'s, and
for $Y^\lc$ in terms of the $T_i$'s.

Finally, $\dc$ can be chosen such as to further preserve the level and
therefore the full alcove picture.

\fi

\iflongversion
\subsection{Explicit (co)ambient space realizations for types $A_n,B_n,C_n,D_n$}
\label{subsection.ambient}

In the sequel, we use for types $A_n$, $B_n$, $C_n$, and $D_n$ the
following ambient space realizations of the finite coroot systems
which realize $\clW$ as groups of signed
permutations~\cite{Bjorner_Brenti.2005,eriksson_eriksson.1998}. For
type $A_n$, we take $\coweightspace=\QQ^{n+1}$ and for types
$B_n,C_n$, and $D_n$ $\coweightspace=\QQ^n$. Denoting by
$(\e^\vee_i)_i$ the canonical basis of $\QQ^{n+1}$ (resp. $\QQ^n$) and
identifying it with its dual basis $(\e_i)_i$, the simple roots are
given by
\begin{equation} \label{eq:simple roots}
\begin{split}
\text{Type $A_n$:} & \qquad
	\alpha_i = \begin{cases}
	\e_{n+1} - \e_1 & \text{for $i=0$,}\\
	\e_i - \e_{i+1} & \text{for $1\le i\le n$;}
	\end{cases}\\[2mm]
\text{Type $B_n$}: & \qquad
	\alpha_i = \begin{cases}
	- \e_1 - \e_2 & \text{for $i=0$,}\\
	\e_i - \e_{i+1} & \text{for $1\le i<n$,}\\
	\e_n & \text{for $i=n$;}
	\end{cases}\\[2mm]
\text{Type $C_n$}: & \qquad 
	\alpha_i = \begin{cases}
	-2 \e_1 & \text{for $i=0$,}\\
	\e_i - \e_{i+1} & \text{for $1\le i<n$,}\\
	2 \e_n & \text{for $i=n$;}
	\end{cases}\\[2mm]
\text{Type $D_n$}: & \qquad
	\alpha_i = \begin{cases}
	-\e_1 - \e_2 & \text{for $i=0$,}\\
	\e_i - \e_{i+1} & \text{for $1\le i<n$,}\\
	\e_{n-1} + \e_n & \text{for $i=n$.}
	\end{cases}
\end{split}
\end{equation}
With this, the action~\eqref{eq:pi} of $\pi_i$ on $\xc=(x_1,x_2,\dots)\in\coweightspace$ becomes
\begin{equation} \label{eq:pi on ambient}
\begin{split}
\text{Type $A_n$:} & \qquad
  	\pi_i(\xc) = \begin{cases}
		(x_{n+1},x_2,\ldots,x_n,x_1) &
		\text{if $i=0$ and $x_{n+1}>x_1$,}\\
		(x_1,\ldots,x_{i+1},x_i,\ldots,x_{n+1}) & 
  		\text{if $1\le i\le n$ and  $x_i>x_{i+1}$,}\\
		\x & \text{otherwise;}
		\end{cases}\\
\text{Type $B_n$:} & \qquad		
   	\pi_i(\xc) = \begin{cases}
		(-x_2,-x_1,x_3,\ldots,x_n) &
		\text{if $i=0$ and $x_1+x_2<0$,}\\
		(x_1,\ldots,x_{i+1},x_i,\ldots,x_n) & 
  		\text{if $1\le i<n$ and  $x_i>x_{i+1}$,}\\
		(x_1,\ldots,x_{n-1},-x_n) &
		\text{if $i=n$ and $x_n>0$,}\\
		\x & \text{otherwise;}
		\end{cases}\\
\text{Type $C_n$:} & \qquad
  	\pi_i(\xc) = \begin{cases}
		(-x_1,x_2,\ldots,x_n) &
		\text{if $i=0$ and $x_1<0$,}\\
		(x_1,\ldots,x_{i+1},x_i,\ldots,x_n) & 
  		\text{if $1\le i<n$ and  $x_i>x_{i+1}$,}\\
		(x_1,\ldots,x_{n-1},-x_n) &
		\text{if $i=n$ and $x_n>0$,}\\
		\x & \text{otherwise;}
		\end{cases}\\
\text{Type $D_n$:} & \qquad
  	\pi_i(\xc) = \begin{cases}
		(-x_2,-x_1,x_3,\ldots,x_n) &
		\text{if $i=0$ and $x_1+x_2<0$,}\\
		(x_1,\ldots,x_{i+1},x_i,\ldots,x_n) & 
  		\text{if $1\le i<n$ and  $x_i>x_{i+1}$,}\\
		(x_1,\ldots,x_{n-2},-x_n,-x_{n-1}) &
		\text{if $i=n$ and $x_{n-1}+x_n>0$,}\\
		\x & \text{otherwise.}
		\end{cases}
\end{split}
\end{equation}
We may pick $\rhoc:=(d,d-1,\ldots,1)$ (where $d$ is the
dimension of $\coweightspace$) as representative of the fundamental
chamber for $\clW$: $\ip{\rhoc,\alpha_i} > 0$, for all $i=1,\dots,n$.
Instead of $s_i$ and $\pi_i$ acting on the coambient space, they can
equivalently act on group elements themselves. The correspondence can
be realized by evaluating $w(\rhoc)$. Whereas the action on the coambient
space~\eqref{eq:pi on ambient} is an action from the left, the action on the group
itself is an action from the right.
\fi

\section{Transitivity of the level $0$ action of affine $0$-Hecke
  algebras}
\label{section.transitivity}

In this section we state and prove the core combinatorial
Theorem~\ref{theorem.connected} of this paper about transitivity of
the level $0$ action of affine $0$-Hecke algebras and mention some
applications to crystal graphs.

\iflongversion
\subsection{Transitivity}
\fi

We start with type $A_n$ to illustrate the results. Here, each $\pi_i$ can be interpreted as a
partial (anti)sort operator: it acts on a permutation (or word)
$w:=(w_1,\ldots,w_{n+1})$ by exchanging $w_i$ and $w_{i+1}$ if
$w_i<w_{i+1}$. By bubble sort, any permutation can be mapped via
$\pi_1,\dots,\pi_n$ to the maximal permutation $\Wmax$, but not
conversely. More precisely the (oriented) graph of the action is the
usual right permutohedron, which is acyclic with $1$ as minimal
element and $\Wmax$ as maximal element.

Consider now $w$ as written along a circle, and let $\pi_0$ act as
above with $i$ taken modulo $n+1$. As suggested by
Figure~\ref{figure.antisorting} for $n=2$, adding the $0$ edges makes
the graph of the action strongly connected.
\begin{proposition} \label{proposition.transitiveA}
  $\pi_0,\dots,\pi_n$ act transitively on permutations of $\{1,\dots,n+1\}$.
\end{proposition}
\iflongversion
\begin{proof}
\fi
  We start with any permutation $w$ and identify it with
  $w(\rhoc)=:\xc=(x_1,\ldots,x_{n+1})$%
  \iflongversion\else%
  , where $\rhoc=(n+1,\ldots,1)$%
  \fi.
  Then the $\pi_i$ act as in~\eqref{eq:pi on ambient}.  

  Suppose that the letter $z=n+1$ is at position $k$ in $\xc$. Then
  $\pi_0\pi_n \cdots \pi_{k+1} \pi_k(\xc)$ has letter $z$ in position
  $1$.  The operator $\tilde{\pi}_0=(\pi_0 \pi_n \cdots
  \pi_1)^{n-1}(\pi_0\pi_n)(\pi_0)(\pi_{n-1}\cdots \pi_1)$ acts in the
  same way as $\pi_0$, except only on the last $n$ positions:
    \begin{equation}
    \def\d{\ar@/^6em/[d]}
    \vcenter{
    \xymatrix@R=.2em{
      z \; x_1 \; x_2 \ldots x_{n-1} \; x_n \d^{\pi_{n-1}\cdots \pi_1}\\
      x_1 \; x_2 \ldots x_{n-1} \; z \; x_n \d^{\pi_0}\\
      x_1' \; x_2 \ldots x_{n-1} \; z \; x_n' \d^{\pi_0\pi_n}\\
      z \; x_2 \ldots x_{n-1} x_n' \; x_1' \d^{(\pi_0 \pi_n \cdots \pi_1)^{n-1}}\\
      z \; x_1' \; x_2 \ldots x_{n-1}\; x_n'
      }}
  \end{equation}
  where $x_1'=x_n$ and $x_n' = x_1$ if $x_n>x_1$ and $x'_1=x_1$ and
  $x'_n=x_n$ otherwise.  In the last step we have used that the
  operator $\pi_0\pi_n \cdots \pi_1$ rotates the last $n$ letters
  cyclically one step to the left, leaving the letter $z$ in position
  $1$ unchanged.
  The result follows by induction.
\iflongversion
\end{proof}
\fi

Let now $\clW$ be any finite Weyl group, and $\heckeW[\clW]{0}$ its
$0$-Hecke algebra. Via $\pi_1,\dots,\pi_n$ the identity of $\clW$ can
be mapped to any $w\in\clW$, but not back (the graph of the action
is just the Hasse diagram of the right weak Bruhat order). Now embed
$\clW$ in an affine Weyl group $\W$, and consider the extra generator
$\pi_0$ of its $0$-Hecke algebra acting on $\clW$. As the
dominant chamber of $\clW$ is on the negative side of $H_{\alpha_0}$,
$\pi_0$ tends to map elements of $\clW$ back to the identity (see
Figure~\ref{figure.antisorting}).
\begin{theorem}
  \label{theorem.connected}
  Let $\W$ be an affine Weyl group, $\clW$ the associated
  finite Weyl group, and
  $\pi_0,\pi_1,\dots,\pi_n$ the generators of the $0$-Hecke algebra of
  $\W$. Then, the level $0$ action of $\pi_0,\pi_1,\dots,\pi_n$
  on $\clW$ (or equivalently on the chambers of $\clW$) is transitive.
\end{theorem}

We prove Theorem~\ref{theorem.connected} by a type free geometric
argument using Lemma~\ref{lemma.positivePathDominant} below.
Figure~\ref{figure.antisorting} illustrates the proof,
\iflongversion%
and thanks to Section~\ref{section.independence} covers 
\else%
actually covering
\fi%
all the rank $2$ affine Weyl groups.

\begin{figure}
  \centerline{$\begin{array}{ccc}
    \begin{tikzpicture}[>=latex,join=bevel,scale=.5,baseline=(current bounding box.east)]
\tiny%
  \node (N_1) at (50bp,211bp) [draw=none] {$123$};
  \node (N_2) at (10bp,143bp) [draw=none] {$132$};
  \node (N_3) at (90bp,143bp) [draw=none] {$213$};
  \node (N_4) at (10bp, 75bp) [draw=none] {$312$};
  \node (N_5) at (90bp, 75bp) [draw=none] {$231$};
  \node (N_6) at (50bp,  7bp) [draw=none] {$321$};

  \draw [->,DarkRed]  (N_1) -- node [left]  {$2$} (N_2);
  \draw [->,DarkBlue] (N_1) -- node [right] {$1$} (N_3);

  \draw [->,DarkBlue] (N_2) -- node [left]  {$1$} (N_4);
  \draw [->,DarkRed]  (N_3) -- node [right] {$2$} (N_5);

  \draw [->,DarkRed]  (N_4) -- node [left]  {$2$} (N_6);
  \draw [->,DarkBlue] (N_5) -- node [right] {$1$} (N_6);

  \draw [->]      (N_5) -- node [left =1mm] {$0$} (N_2);
  \draw [->]      (N_4) -- node [right=1mm] {$0$} (N_3);

  \draw [->]      (N_6) -- 
                            (N_1);
\end{tikzpicture} &
    \begin{tikzpicture}[>=latex,join=bevel,scale=.5,baseline=(current bounding box.east)]
\tiny%
  \node (N_4) at (50bp,287bp) [draw,draw=none] {$12$};

  \node (N_6) at (10bp,218bp) [draw,draw=none] {$21$};
  \node (N_2) at (90bp,218bp) [draw,draw=none] {$1\underline{2}$};

  \node (N_8) at (10bp,148bp) [draw,draw=none] {$2\underline{1}$};
  \node (N_5) at (90bp,148bp) [draw,draw=none] {$\underline{2}1$};

  \node (N_3) at (10bp,78bp) [draw,draw=none] {$\underline{1}2$};
  \node (N_7) at (90bp,78bp) [draw,draw=none] {$\underline{2}\underline{1}$};

  \node (N_1) at (50bp,8bp) [draw,draw=none] {$\underline{1}\underline{2}$};

  \draw [->,DarkBlue] (N_4) -- node [left]  {$1$} (N_6);
  \draw [->,DarkRed]  (N_4) -- node [right] {$2$} (N_2);

  \draw [->,DarkRed]  (N_6) -- node [left]  {$2$} (N_8);
  \draw [->,DarkBlue] (N_2) -- node [right] {$1$} (N_5);

  \draw [->,DarkBlue] (N_8) -- node [left]  {$1$} (N_3);
  \draw [->,DarkRed]  (N_5) -- node [right] {$2$} (N_7);

  \draw [->]      (N_5) -- node [above] {$0$} (N_6);

  \draw [->,DarkRed]  (N_3) -- node [left]  {$2$} (N_1);
  \draw [->,DarkBlue] (N_7) -- node [right] {$1$} (N_1);
  \draw [->]      (N_3) -- node [left]  {$0$} (N_4);
  \draw [->]      (N_7) -- node [above] {$0$} (N_8);

  \draw [->]      (N_1) -- node [right] {$0$} (N_2);
\end{tikzpicture}&
    \scalebox{0.75}{\begin{tikzpicture}[>=latex,join=bevel,scale=.5,baseline=(current bounding box.east)]
\tiny%
  \node (N_1) at (50bp,376bp) [draw,draw=none] {$1$};

  \node (N_2) at (10bp,314bp) [draw,draw=none] {$$};
  \node (N_3) at (90bp,314bp) [draw,draw=none] {$$};

  \node (N_4) at (10bp,252bp) [draw,draw=none] {$$};
  \node (N_11) at (90bp,252bp) [draw,draw=none] {$$};

  \node (N_7) at (10bp,190bp) [draw,draw=none] {$$};
  \node (N_5) at (90bp,190bp) [draw,draw=none] {$$};

  \node (N_9) at (10bp,128bp) [draw,draw=none] {$$};
  \node (N_6) at (90bp,128bp) [draw,draw=none] {$$};

  \node (N_10) at (10bp,66bp) [draw,draw=none] {$$};
  \node (N_8) at (90bp,66bp) [draw,draw=none] {$$};

  \node (N_12) at (50bp,4bp) [draw,draw=none] {$\Wmax$};

  \draw [->,DarkRed] (N_1) -- node [left] {$2$} (N_2);
  \draw [->,DarkBlue] (N_1) -- node [right] {$1$} (N_3);

  \draw [->,DarkBlue] (N_2) -- node [left] {$1$} (N_4);
  \draw [->,DarkRed] (N_3) -- node [right] {$2$} (N_11);

  \draw [->,DarkRed] (N_4) -- node [left] {$2$} (N_7);
  \draw [->,DarkBlue] (N_11) -- node [right] {$1$} (N_5);

  \draw [->,DarkBlue] (N_7) -- node [left] {$1$} (N_9);
  \draw [->,DarkRed] (N_5) -- node [right] {$2$} (N_6);
  \draw [->] (N_5) -- node [above=2mm]{$0$} (N_4);

  \draw [->,DarkRed] (N_9) -- node [left] {$2$} (N_10);
  \draw [->,DarkBlue] (N_6) -- node [right] {$1$} (N_8);
  \draw [->] (N_9) -- node [right] {$0$} (N_3);
  \draw [->] (N_6) -- node [above=2mm] {$0$} (N_7);

  \draw [->,DarkBlue] (N_10) -- node [left] {$1$} (N_12);
  \draw [->,DarkRed] (N_8) -- node [right] {$2$} (N_12);
  \draw [->] (N_10) -- node [right]{$0$} (N_11);
  \draw [->] (N_8) -- node [above=3mm] {$\,0$} (N_1);

  \draw [->] (N_12) -- node [below=3mm] {$\!0$} (N_2);

\end{tikzpicture}}\\\\
    \iflongversion
    \scalebox{.9}{\begin{tikzpicture}[baseline=(current bounding box.east)]
\draw[->,  color = black,] (0.0,0.0) -- (0.0,-1.732050808) node[at end, auto=right] {$\alpha^\vee_{0}$};
\draw[->,  color = DarkBlue,] (0.0,0.0) -- (1.5,0.8660254038) node[at end, auto=right] {$\alpha^\vee_{1}$};
\draw[->,  color = DarkRed,] (0.0,0.0) -- (-1.5,0.8660254038) node[at end, auto=right] {$\alpha^\vee_{2}$};
\draw[ color = black, very thick,](0.5,0.8660254038) -- (-0.5,0.8660254038);
\draw[ color = DarkBlue, ,](0.0,0.0) -- (-0.5,0.8660254038);
\draw[ color = DarkRed, very thin,](0.0,0.0) -- (0.5,0.8660254038);
\draw[ color = black, very thick,](0.5,0.8660254038) -- (1.0,0.0);
\draw[ color = DarkBlue, ,](0.0,0.0) -- (1.0,0.0);
\draw[ color = DarkRed, very thin,](0.0,0.0) -- (0.5,0.8660254038);
\draw[ color = black, very thick,](0.5,-0.8660254038) -- (1.0,0.0);
\draw[ color = DarkBlue, ,](0.0,0.0) -- (1.0,0.0);
\draw[ color = DarkRed, very thin,](0.0,0.0) -- (0.5,-0.8660254038);
\draw[ color = black, very thick,](0.5,-0.8660254038) -- (-0.5,-0.8660254038);
\draw[ color = DarkBlue, ,](0.0,0.0) -- (-0.5,-0.8660254038);
\draw[ color = DarkRed, very thin,](0.0,0.0) -- (0.5,-0.8660254038);
\draw[ color = black, very thick,](-1.0,0.0) -- (-0.5,0.8660254038);
\draw[ color = DarkBlue, ,](0.0,0.0) -- (-0.5,0.8660254038);
\draw[ color = DarkRed, very thin,](0.0,0.0) -- (-1.0,0.0);
\draw[ color = black, very thick,](-1.0,0.0) -- (-0.5,-0.8660254038);
\draw[ color = DarkBlue, ,](0.0,0.0) -- (-0.5,-0.8660254038);
\draw[ color = DarkRed, very thin,](0.0,0.0) -- (-1.0,0.0);
\draw[ color = black, very thick,](0.5,2.598076211) -- (-0.5,2.598076211);
\draw[ color = DarkBlue, ,](0.0,1.732050808) -- (-0.5,2.598076211);
\draw[ color = DarkRed, very thin,](0.0,1.732050808) -- (0.5,2.598076211);
\draw[ color = black, very thick,](0.5,2.598076211) -- (1.0,1.732050808);
\draw[ color = DarkBlue, ,](0.0,1.732050808) -- (1.0,1.732050808);
\draw[ color = DarkRed, very thin,](0.0,1.732050808) -- (0.5,2.598076211);
\draw[ color = black, very thick,](0.5,0.8660254038) -- (1.0,1.732050808);
\draw[ color = DarkBlue, ,](0.0,1.732050808) -- (1.0,1.732050808);
\draw[ color = DarkRed, very thin,](0.0,1.732050808) -- (0.5,0.8660254038);
\draw[ color = black, very thick,](0.5,0.8660254038) -- (-0.5,0.8660254038);
\draw[ color = DarkBlue, ,](0.0,1.732050808) -- (-0.5,0.8660254038);
\draw[ color = DarkRed, very thin,](0.0,1.732050808) -- (0.5,0.8660254038);
\draw[ color = black, very thick,](-1.0,1.732050808) -- (-0.5,2.598076211);
\draw[ color = DarkBlue, ,](0.0,1.732050808) -- (-0.5,2.598076211);
\draw[ color = DarkRed, very thin,](0.0,1.732050808) -- (-1.0,1.732050808);
\draw[ color = black, very thick,](-1.0,1.732050808) -- (-0.5,0.8660254038);
\draw[ color = DarkBlue, ,](0.0,1.732050808) -- (-0.5,0.8660254038);
\draw[ color = DarkRed, very thin,](0.0,1.732050808) -- (-1.0,1.732050808);
\draw[->,  color = purple,thick] (0.0,1.154700538) -- (0.0,0.5773502692);
\end{tikzpicture}}&
    \scalebox{1.2}{\begin{tikzpicture}[baseline=(current bounding box.east)]
\draw[->,  color = black,] (0.0,0.0) -- (0.0,-1.0) node[at end, auto=right] {$\alpha^\vee_{0}$};
\draw[->,  color = DarkBlue,] (0.0,0.0) -- (-1.0,1.0) node[at end, auto=right] {$\alpha^\vee_{1}$};
\draw[->,  color = DarkRed,] (0.0,0.0) -- (1.0,0.0) node[at end, auto=right] {$\alpha^\vee_{2}$};
\draw[ color = black, very thick,](0.0,0.5) -- (0.5,0.5);
\draw[ color = DarkBlue, ,](0.0,0.0) -- (0.5,0.5);
\draw[ color = DarkRed, very thin,](0.0,0.0) -- (0.0,0.5);
\draw[ color = black, very thick,](0.0,0.5) -- (-0.5,0.5);
\draw[ color = DarkBlue, ,](0.0,0.0) -- (-0.5,0.5);
\draw[ color = DarkRed, very thin,](0.0,0.0) -- (0.0,0.5);
\draw[ color = black, very thick,](0.0,-0.5) -- (0.5,-0.5);
\draw[ color = DarkBlue, ,](0.0,0.0) -- (0.5,-0.5);
\draw[ color = DarkRed, very thin,](0.0,0.0) -- (0.0,-0.5);
\draw[ color = black, very thick,](0.0,-0.5) -- (-0.5,-0.5);
\draw[ color = DarkBlue, ,](0.0,0.0) -- (-0.5,-0.5);
\draw[ color = DarkRed, very thin,](0.0,0.0) -- (0.0,-0.5);
\draw[ color = black, very thick,](-0.5,0.0) -- (-0.5,0.5);
\draw[ color = DarkBlue, ,](0.0,0.0) -- (-0.5,0.5);
\draw[ color = DarkRed, very thin,](0.0,0.0) -- (-0.5,0.0);
\draw[ color = black, very thick,](-0.5,0.0) -- (-0.5,-0.5);
\draw[ color = DarkBlue, ,](0.0,0.0) -- (-0.5,-0.5);
\draw[ color = DarkRed, very thin,](0.0,0.0) -- (-0.5,0.0);
\draw[ color = black, very thick,](0.5,0.0) -- (0.5,0.5);
\draw[ color = DarkBlue, ,](0.0,0.0) -- (0.5,0.5);
\draw[ color = DarkRed, very thin,](0.0,0.0) -- (0.5,0.0);
\draw[ color = black, very thick,](0.5,0.0) -- (0.5,-0.5);
\draw[ color = DarkBlue, ,](0.0,0.0) -- (0.5,-0.5);
\draw[ color = DarkRed, very thin,](0.0,0.0) -- (0.5,0.0);
\draw[ color = black, very thick,](0.0,1.5) -- (0.5,1.5);
\draw[ color = DarkBlue, ,](0.0,1.0) -- (0.5,1.5);
\draw[ color = DarkRed, very thin,](0.0,1.0) -- (0.0,1.5);
\draw[ color = black, very thick,](0.0,1.5) -- (-0.5,1.5);
\draw[ color = DarkBlue, ,](0.0,1.0) -- (-0.5,1.5);
\draw[ color = DarkRed, very thin,](0.0,1.0) -- (0.0,1.5);
\draw[ color = black, very thick,](0.0,0.5) -- (0.5,0.5);
\draw[ color = DarkBlue, ,](0.0,1.0) -- (0.5,0.5);
\draw[ color = DarkRed, very thin,](0.0,1.0) -- (0.0,0.5);
\draw[ color = black, very thick,](0.0,0.5) -- (-0.5,0.5);
\draw[ color = DarkBlue, ,](0.0,1.0) -- (-0.5,0.5);
\draw[ color = DarkRed, very thin,](0.0,1.0) -- (0.0,0.5);
\draw[ color = black, very thick,](-0.5,1.0) -- (-0.5,1.5);
\draw[ color = DarkBlue, ,](0.0,1.0) -- (-0.5,1.5);
\draw[ color = DarkRed, very thin,](0.0,1.0) -- (-0.5,1.0);
\draw[ color = black, very thick,](-0.5,1.0) -- (-0.5,0.5);
\draw[ color = DarkBlue, ,](0.0,1.0) -- (-0.5,0.5);
\draw[ color = DarkRed, very thin,](0.0,1.0) -- (-0.5,1.0);
\draw[ color = black, very thick,](0.5,1.0) -- (0.5,1.5);
\draw[ color = DarkBlue, ,](0.0,1.0) -- (0.5,1.5);
\draw[ color = DarkRed, very thin,](0.0,1.0) -- (0.5,1.0);
\draw[ color = black, very thick,](0.5,1.0) -- (0.5,0.5);
\draw[ color = DarkBlue, ,](0.0,1.0) -- (0.5,0.5);
\draw[ color = DarkRed, very thin,](0.0,1.0) -- (0.5,1.0);
\draw[ color = black, very thick,](1.0,1.5) -- (1.5,1.5);
\draw[ color = DarkBlue, ,](1.0,1.0) -- (1.5,1.5);
\draw[ color = DarkRed, very thin,](1.0,1.0) -- (1.0,1.5);
\draw[ color = black, very thick,](1.0,1.5) -- (0.5,1.5);
\draw[ color = DarkBlue, ,](1.0,1.0) -- (0.5,1.5);
\draw[ color = DarkRed, very thin,](1.0,1.0) -- (1.0,1.5);
\draw[ color = black, very thick,](1.0,0.5) -- (1.5,0.5);
\draw[ color = DarkBlue, ,](1.0,1.0) -- (1.5,0.5);
\draw[ color = DarkRed, very thin,](1.0,1.0) -- (1.0,0.5);
\draw[ color = black, very thick,](1.0,0.5) -- (0.5,0.5);
\draw[ color = DarkBlue, ,](1.0,1.0) -- (0.5,0.5);
\draw[ color = DarkRed, very thin,](1.0,1.0) -- (1.0,0.5);
\draw[ color = black, very thick,](0.5,1.0) -- (0.5,1.5);
\draw[ color = DarkBlue, ,](1.0,1.0) -- (0.5,1.5);
\draw[ color = DarkRed, very thin,](1.0,1.0) -- (0.5,1.0);
\draw[ color = black, very thick,](0.5,1.0) -- (0.5,0.5);
\draw[ color = DarkBlue, ,](1.0,1.0) -- (0.5,0.5);
\draw[ color = DarkRed, very thin,](1.0,1.0) -- (0.5,1.0);
\draw[ color = black, very thick,](1.5,1.0) -- (1.5,1.5);
\draw[ color = DarkBlue, ,](1.0,1.0) -- (1.5,1.5);
\draw[ color = DarkRed, very thin,](1.0,1.0) -- (1.5,1.0);
\draw[ color = black, very thick,](1.5,1.0) -- (1.5,0.5);
\draw[ color = DarkBlue, ,](1.0,1.0) -- (1.5,0.5);
\draw[ color = DarkRed, very thin,](1.0,1.0) -- (1.5,1.0);
\draw[ color = black, very thick,](0.0,2.5) -- (0.5,2.5);
\draw[ color = DarkBlue, ,](0.0,2.0) -- (0.5,2.5);
\draw[ color = DarkRed, very thin,](0.0,2.0) -- (0.0,2.5);
\draw[ color = black, very thick,](0.0,2.5) -- (-0.5,2.5);
\draw[ color = DarkBlue, ,](0.0,2.0) -- (-0.5,2.5);
\draw[ color = DarkRed, very thin,](0.0,2.0) -- (0.0,2.5);
\draw[ color = black, very thick,](0.0,1.5) -- (0.5,1.5);
\draw[ color = DarkBlue, ,](0.0,2.0) -- (0.5,1.5);
\draw[ color = DarkRed, very thin,](0.0,2.0) -- (0.0,1.5);
\draw[ color = black, very thick,](0.0,1.5) -- (-0.5,1.5);
\draw[ color = DarkBlue, ,](0.0,2.0) -- (-0.5,1.5);
\draw[ color = DarkRed, very thin,](0.0,2.0) -- (0.0,1.5);
\draw[ color = black, very thick,](-0.5,2.0) -- (-0.5,2.5);
\draw[ color = DarkBlue, ,](0.0,2.0) -- (-0.5,2.5);
\draw[ color = DarkRed, very thin,](0.0,2.0) -- (-0.5,2.0);
\draw[ color = black, very thick,](-0.5,2.0) -- (-0.5,1.5);
\draw[ color = DarkBlue, ,](0.0,2.0) -- (-0.5,1.5);
\draw[ color = DarkRed, very thin,](0.0,2.0) -- (-0.5,2.0);
\draw[ color = black, very thick,](0.5,2.0) -- (0.5,2.5);
\draw[ color = DarkBlue, ,](0.0,2.0) -- (0.5,2.5);
\draw[ color = DarkRed, very thin,](0.0,2.0) -- (0.5,2.0);
\draw[ color = black, very thick,](0.5,2.0) -- (0.5,1.5);
\draw[ color = DarkBlue, ,](0.0,2.0) -- (0.5,1.5);
\draw[ color = DarkRed, very thin,](0.0,2.0) -- (0.5,2.0);
\draw[ color = black, very thick,](1.0,2.5) -- (1.5,2.5);
\draw[ color = DarkBlue, ,](1.0,2.0) -- (1.5,2.5);
\draw[ color = DarkRed, very thin,](1.0,2.0) -- (1.0,2.5);
\draw[ color = black, very thick,](1.0,2.5) -- (0.5,2.5);
\draw[ color = DarkBlue, ,](1.0,2.0) -- (0.5,2.5);
\draw[ color = DarkRed, very thin,](1.0,2.0) -- (1.0,2.5);
\draw[ color = black, very thick,](1.0,1.5) -- (1.5,1.5);
\draw[ color = DarkBlue, ,](1.0,2.0) -- (1.5,1.5);
\draw[ color = DarkRed, very thin,](1.0,2.0) -- (1.0,1.5);
\draw[ color = black, very thick,](1.0,1.5) -- (0.5,1.5);
\draw[ color = DarkBlue, ,](1.0,2.0) -- (0.5,1.5);
\draw[ color = DarkRed, very thin,](1.0,2.0) -- (1.0,1.5);
\draw[ color = black, very thick,](0.5,2.0) -- (0.5,2.5);
\draw[ color = DarkBlue, ,](1.0,2.0) -- (0.5,2.5);
\draw[ color = DarkRed, very thin,](1.0,2.0) -- (0.5,2.0);
\draw[ color = black, very thick,](0.5,2.0) -- (0.5,1.5);
\draw[ color = DarkBlue, ,](1.0,2.0) -- (0.5,1.5);
\draw[ color = DarkRed, very thin,](1.0,2.0) -- (0.5,2.0);
\draw[ color = black, very thick,](1.5,2.0) -- (1.5,2.5);
\draw[ color = DarkBlue, ,](1.0,2.0) -- (1.5,2.5);
\draw[ color = DarkRed, very thin,](1.0,2.0) -- (1.5,2.0);
\draw[ color = black, very thick,](1.5,2.0) -- (1.5,1.5);
\draw[ color = DarkBlue, ,](1.0,2.0) -- (1.5,1.5);
\draw[ color = DarkRed, very thin,](1.0,2.0) -- (1.5,2.0);
\draw[ color = black, very thick,](2.0,2.5) -- (2.5,2.5);
\draw[ color = DarkBlue, ,](2.0,2.0) -- (2.5,2.5);
\draw[ color = DarkRed, very thin,](2.0,2.0) -- (2.0,2.5);
\draw[ color = black, very thick,](2.0,2.5) -- (1.5,2.5);
\draw[ color = DarkBlue, ,](2.0,2.0) -- (1.5,2.5);
\draw[ color = DarkRed, very thin,](2.0,2.0) -- (2.0,2.5);
\draw[ color = black, very thick,](2.0,1.5) -- (2.5,1.5);
\draw[ color = DarkBlue, ,](2.0,2.0) -- (2.5,1.5);
\draw[ color = DarkRed, very thin,](2.0,2.0) -- (2.0,1.5);
\draw[ color = black, very thick,](2.0,1.5) -- (1.5,1.5);
\draw[ color = DarkBlue, ,](2.0,2.0) -- (1.5,1.5);
\draw[ color = DarkRed, very thin,](2.0,2.0) -- (2.0,1.5);
\draw[ color = black, very thick,](1.5,2.0) -- (1.5,2.5);
\draw[ color = DarkBlue, ,](2.0,2.0) -- (1.5,2.5);
\draw[ color = DarkRed, very thin,](2.0,2.0) -- (1.5,2.0);
\draw[ color = black, very thick,](1.5,2.0) -- (1.5,1.5);
\draw[ color = DarkBlue, ,](2.0,2.0) -- (1.5,1.5);
\draw[ color = DarkRed, very thin,](2.0,2.0) -- (1.5,2.0);
\draw[ color = black, very thick,](2.5,2.0) -- (2.5,2.5);
\draw[ color = DarkBlue, ,](2.0,2.0) -- (2.5,2.5);
\draw[ color = DarkRed, very thin,](2.0,2.0) -- (2.5,2.0);
\draw[ color = black, very thick,](2.5,2.0) -- (2.5,1.5);
\draw[ color = DarkBlue, ,](2.0,2.0) -- (2.5,1.5);
\draw[ color = DarkRed, very thin,](2.0,2.0) -- (2.5,2.0);
\draw[->,  color = purple,] (0.875,1.625) -- (0.875,1.375);
\draw[->,  color = purple,] (0.875,1.375) -- (0.625,1.125);
\draw[->,  color = purple,] (0.625,1.125) -- (0.375,1.125);
\draw[->,  color = purple,] (0.375,1.125) -- (0.375,0.875);
\draw[->,  color = purple,] (0.375,0.875) -- (0.125,0.625);
\draw[->,  color = purple,] (0.125,0.625) -- (0.125,0.375);
\end{tikzpicture}}&
    \scalebox{.9}{\input{minimalAntisortingAlcoveWalk-ambient-G2}}\\
    \fi
    \scalebox{.9}{\begin{tikzpicture}[baseline=(current bounding box.east)]
\draw[->,  color = black,] (0.0,0.0) -- (0.0,-1.732050808) node[at end, auto=right] {$\alpha^\vee_{0}$};
\draw[->,  color = DarkBlue,] (0.0,0.0) -- (1.5,0.8660254038) node[at end, auto=right] {$\alpha^\vee_{1}$};
\draw[->,  color = DarkRed,] (0.0,0.0) -- (-1.5,0.8660254038) node[at end, auto=right] {$\alpha^\vee_{2}$};
\draw[ color = black, very thick,](0.5,0.8660254038) -- (-0.5,0.8660254038);
\draw[ color = DarkBlue, ,](0.0,0.0) -- (-0.5,0.8660254038);
\draw[ color = DarkRed, very thin,](0.0,0.0) -- (0.5,0.8660254038);
\draw[ color = black, very thick,](0.5,0.8660254038) -- (1.0,0.0);
\draw[ color = DarkBlue, ,](0.0,0.0) -- (1.0,0.0);
\draw[ color = DarkRed, very thin,](0.0,0.0) -- (0.5,0.8660254038);
\draw[ color = black, very thick,](0.5,-0.8660254038) -- (1.0,0.0);
\draw[ color = DarkBlue, ,](0.0,0.0) -- (1.0,0.0);
\draw[ color = DarkRed, very thin,](0.0,0.0) -- (0.5,-0.8660254038);
\draw[ color = black, very thick,](0.5,-0.8660254038) -- (-0.5,-0.8660254038);
\draw[ color = DarkBlue, ,](0.0,0.0) -- (-0.5,-0.8660254038);
\draw[ color = DarkRed, very thin,](0.0,0.0) -- (0.5,-0.8660254038);
\draw[ color = black, very thick,](-1.0,0.0) -- (-0.5,0.8660254038);
\draw[ color = DarkBlue, ,](0.0,0.0) -- (-0.5,0.8660254038);
\draw[ color = DarkRed, very thin,](0.0,0.0) -- (-1.0,0.0);
\draw[ color = black, very thick,](-1.0,0.0) -- (-0.5,-0.8660254038);
\draw[ color = DarkBlue, ,](0.0,0.0) -- (-0.5,-0.8660254038);
\draw[ color = DarkRed, very thin,](0.0,0.0) -- (-1.0,0.0);
\draw[->,  color = purple,thick] (0.0,-0.5773502692) -- (0.0,0.5773502692);
\end{tikzpicture}}&
    \scalebox{1.2}{\begin{tikzpicture}[baseline=(current bounding box.east)]
\draw[->,  color = black,] (0.0,0.0) -- (0.0,-1.0) node[at end, auto=right] {$\alpha^\vee_{0}$};
\draw[->,  color = DarkBlue,] (0.0,0.0) -- (-1.0,1.0) node[at end, auto=right] {$\alpha^\vee_{1}$};
\draw[->,  color = DarkRed,] (0.0,0.0) -- (1.0,0.0) node[at end, auto=right] {$\alpha^\vee_{2}$};
\draw[ color = black, very thick,](0.0,0.5) -- (0.5,0.5);
\draw[ color = DarkBlue, ,](0.0,0.0) -- (0.5,0.5);
\draw[ color = DarkRed, very thin,](0.0,0.0) -- (0.0,0.5);
\draw[ color = black, very thick,](0.0,0.5) -- (-0.5,0.5);
\draw[ color = DarkBlue, ,](0.0,0.0) -- (-0.5,0.5);
\draw[ color = DarkRed, very thin,](0.0,0.0) -- (0.0,0.5);
\draw[ color = black, very thick,](0.0,-0.5) -- (0.5,-0.5);
\draw[ color = DarkBlue, ,](0.0,0.0) -- (0.5,-0.5);
\draw[ color = DarkRed, very thin,](0.0,0.0) -- (0.0,-0.5);
\draw[ color = black, very thick,](0.0,-0.5) -- (-0.5,-0.5);
\draw[ color = DarkBlue, ,](0.0,0.0) -- (-0.5,-0.5);
\draw[ color = DarkRed, very thin,](0.0,0.0) -- (0.0,-0.5);
\draw[ color = black, very thick,](-0.5,0.0) -- (-0.5,0.5);
\draw[ color = DarkBlue, ,](0.0,0.0) -- (-0.5,0.5);
\draw[ color = DarkRed, very thin,](0.0,0.0) -- (-0.5,0.0);
\draw[ color = black, very thick,](-0.5,0.0) -- (-0.5,-0.5);
\draw[ color = DarkBlue, ,](0.0,0.0) -- (-0.5,-0.5);
\draw[ color = DarkRed, very thin,](0.0,0.0) -- (-0.5,0.0);
\draw[ color = black, very thick,](0.5,0.0) -- (0.5,0.5);
\draw[ color = DarkBlue, ,](0.0,0.0) -- (0.5,0.5);
\draw[ color = DarkRed, very thin,](0.0,0.0) -- (0.5,0.0);
\draw[ color = black, very thick,](0.5,0.0) -- (0.5,-0.5);
\draw[ color = DarkBlue, ,](0.0,0.0) -- (0.5,-0.5);
\draw[ color = DarkRed, very thin,](0.0,0.0) -- (0.5,0.0);
\draw[->,  color = purple,] (-0.125,-0.375) -- (-0.125,0.375);
\draw[->,  color = purple,] (-0.125,0.375) -- (-0.375,0.125);
\draw[->,  color = purple,] (-0.375,0.125) -- (0.375,0.125);
\draw[->,  color = purple,] (0.375,0.125) -- (0.375,-0.125);
\draw[->,  color = purple,] (0.375,-0.125) -- (0.125,-0.375);
\draw[->,  color = purple,] (0.125,-0.375) -- (0.125,0.375);
\end{tikzpicture}}&
    \scalebox{.9}{\begin{tikzpicture}[baseline=(current bounding box.east)]
\draw[->,  color = black,] (0.0,0.0) -- (0.0,-1.732050808) node[at end, auto=right] {$\alpha^\vee_{0}$};
\draw[->,  color = DarkBlue,] (0.0,0.0) -- (-1.5,0.8660254038) node[at end, auto=right] {$\alpha^\vee_{1}$};
\draw[->,  color = DarkRed,] (0.0,0.0) -- (3.0,0.0) node[at end, auto=right] {$\alpha^\vee_{2}$};
\draw[ color = black, very thick,](0.0,0.8660254038) -- (0.5,0.8660254038);
\draw[ color = DarkBlue, ,](0.0,0.0) -- (0.5,0.8660254038);
\draw[ color = DarkRed, very thin,](0.0,0.0) -- (0.0,0.8660254038);
\draw[ color = black, very thick,](0.0,0.8660254038) -- (-0.5,0.8660254038);
\draw[ color = DarkBlue, ,](0.0,0.0) -- (-0.5,0.8660254038);
\draw[ color = DarkRed, very thin,](0.0,0.0) -- (0.0,0.8660254038);
\draw[ color = black, very thick,](0.0,-0.8660254038) -- (0.5,-0.8660254038);
\draw[ color = DarkBlue, ,](0.0,0.0) -- (0.5,-0.8660254038);
\draw[ color = DarkRed, very thin,](0.0,0.0) -- (0.0,-0.8660254038);
\draw[ color = black, very thick,](0.0,-0.8660254038) -- (-0.5,-0.8660254038);
\draw[ color = DarkBlue, ,](0.0,0.0) -- (-0.5,-0.8660254038);
\draw[ color = DarkRed, very thin,](0.0,0.0) -- (0.0,-0.8660254038);
\draw[ color = black, very thick,](0.75,0.4330127019) -- (0.5,0.8660254038);
\draw[ color = DarkBlue, ,](0.0,0.0) -- (0.5,0.8660254038);
\draw[ color = DarkRed, very thin,](0.0,0.0) -- (0.75,0.4330127019);
\draw[ color = black, very thick,](0.75,0.4330127019) -- (1.0,0.0);
\draw[ color = DarkBlue, ,](0.0,0.0) -- (1.0,0.0);
\draw[ color = DarkRed, very thin,](0.0,0.0) -- (0.75,0.4330127019);
\draw[ color = black, very thick,](-0.75,0.4330127019) -- (-0.5,0.8660254038);
\draw[ color = DarkBlue, ,](0.0,0.0) -- (-0.5,0.8660254038);
\draw[ color = DarkRed, very thin,](0.0,0.0) -- (-0.75,0.4330127019);
\draw[ color = black, very thick,](-0.75,0.4330127019) -- (-1.0,0.0);
\draw[ color = DarkBlue, ,](0.0,0.0) -- (-1.0,0.0);
\draw[ color = DarkRed, very thin,](0.0,0.0) -- (-0.75,0.4330127019);
\draw[ color = black, very thick,](-0.75,-0.4330127019) -- (-1.0,0.0);
\draw[ color = DarkBlue, ,](0.0,0.0) -- (-1.0,0.0);
\draw[ color = DarkRed, very thin,](0.0,0.0) -- (-0.75,-0.4330127019);
\draw[ color = black, very thick,](-0.75,-0.4330127019) -- (-0.5,-0.8660254038);
\draw[ color = DarkBlue, ,](0.0,0.0) -- (-0.5,-0.8660254038);
\draw[ color = DarkRed, very thin,](0.0,0.0) -- (-0.75,-0.4330127019);
\draw[ color = black, very thick,](0.75,-0.4330127019) -- (1.0,0.0);
\draw[ color = DarkBlue, ,](0.0,0.0) -- (1.0,0.0);
\draw[ color = DarkRed, very thin,](0.0,0.0) -- (0.75,-0.4330127019);
\draw[ color = black, very thick,](0.75,-0.4330127019) -- (0.5,-0.8660254038);
\draw[ color = DarkBlue, ,](0.0,0.0) -- (0.5,-0.8660254038);
\draw[ color = DarkRed, very thin,](0.0,0.0) -- (0.75,-0.4330127019);
\draw[->,  color = purple,] (-0.25,-0.7216878365) -- (-0.25,0.7216878365);
\draw[->,  color = purple,] (-0.25,0.7216878365) -- (-0.5,0.5773502692);
\draw[->,  color = purple,] (-0.5,0.5773502692) -- (-0.75,0.1443375673);
\draw[->,  color = purple,] (-0.75,0.1443375673) -- (-0.75,-0.1443375673);
\draw[->,  color = purple,] (-0.75,-0.1443375673) -- (-0.5,-0.5773502692);
\draw[->,  color = purple,] (-0.5,-0.5773502692) -- (0.75,0.1443375673);
\draw[->,  color = purple,] (0.75,0.1443375673) -- (0.75,-0.1443375673);
\draw[->,  color = purple,] (0.75,-0.1443375673) -- (0.5,-0.5773502692);
\draw[->,  color = purple,] (0.5,-0.5773502692) -- (0.25,-0.7216878365);
\draw[->,  color = purple,] (0.25,-0.7216878365) -- (0.25,0.7216878365);
\end{tikzpicture}}\\
    \iflongversion
    \dynkinAIIa & \dynkinCIIa & \dynkinGIIa\\[4mm]
    \fi
    \widetilde A_2 = A_2^{(1)}
    & \widetilde C_2 = C_2^{(1)} & \widetilde G_2 = G_2^{(1)}
  \end{array}$}
  \caption{
    \textbf{Top:} Graph of the action of
    $\pi_0,\pi_1,\dots,\pi_n$ on the finite Weyl group $\clW$,
    using (signed) permutation notation. 
    %
    \newline
    \textbf{Center:} 
    The alcove picture in the ambient space, with a shortest alcove 
    walk from an alcove $w(A)$
    in the dominant chamber such that $\cl(w)=\Wmax$ down to the
    fundamental alcove $A$. An $i$-crossing is negative if it goes down or
    straight to the left.
    \newline
    \textbf{Bottom:} The top graph can be realized geometrically in
    the \emph{Steinberg torus}, quotient of
    the alcove picture by the translations, or equivalently by
    identification of the opposite edges of the fundamental polygon. An $i$-arrow in the
    graph corresponds to a negative $i$-crossing.
    The alcove walk of the center figure
    then becomes a path from the antifundamental
    chamber $\Wmax$(A) back the fundamental chamber $A$.
  }    
  \FIXME{newlines or not?}
  \label{figure.antisorting}
\end{figure}

\begin{lemma}[Cf. Remark 3.5 of~\cite{Ram.2006}]
  \label{lemma.positivePathDominant}
  Let $w(A)$ be an alcove in the dominant chamber of $\clW$, and
  consider a shortest alcove walk $i_1,\dots,i_r$ from $A$ to
  $w(A)$. Then, each crossing is positive. In particular, $i_k$ is a
  descent of $cl(s_{i_1}\!\cdots s_{i_{k-1}})$.
\end{lemma}

\begin{proof}
  If $w(A)$ is the fundamental alcove $A$, the path is empty, and we
  are done. Otherwise, let $H_{\alpha,m}$ be the wall separating
  $w(A)$ from the previous alcove $s_{i_1}\!\cdots
  s_{i_{r-1}}(A)$. Assume that $w(A)$ is in $H^-_{\alpha,m}$.  Taking
  some point $\xc$ in $w(A)$,
  \begin{equation}
    0 > \ip{\xc,\alpha-\delta m} = \ip{\xc,\alpha} - \ell m\,.
  \end{equation}
  Then, using that $w(A)$ is in the fundamental chamber, $m >
  \frac1\ell\ip{\xc,\alpha} > 0$.
  On the other hand, since the alcove walk is shortest, $H_{\alpha,m}$
  separates $w(A)$ and $A$, so $A\in H^+_{\alpha,m}$. Since $0^\ell$
  is in the closure of $A$, $0\leq\ip{0^\ell,\alpha-\delta m} = 0
  -\ell m$. It follows that $m\leq 0$, a contradiction.
\end{proof}

\begin{proof}[Proof of Theorem~\ref{theorem.connected}]
  Take $w\in \clW$, and $w(A)$ the corresponding alcove. One can choose
  a long enough strictly dominant element $\lc$ of the coroot lattice
  so that $t_\lc(w(A))$ lies in the dominant chamber of $\clW$. Consider some
  shortest alcove walk $i_1,\dots,i_r$ from $t_\lc(w(A))$ back to the
  fundamental alcove $A$ (see Figure~\ref{figure.antisorting}). Then, in $\clW$, $w\cl(s_{i_1})\cdots
  \cl(s_{i_r})=1$. Furthermore, by
  Lemma~\ref{lemma.positivePathDominant}, at each step $i_k$ is not a
  descent of $w\cl(s_{i_1})\cdots\cl(s_{i_{k-1}})$. Therefore,
  $w.\pi_{i_1}\dots\pi_{i_r}=w \cl(s_{i_1})\cdots \cl(s_{i_r}) = 1$, as desired.
\end{proof}

We now exhibit a recursive sorting algorithm for type $B_n$,
\iflongversion
where the operators $\pi_i$ act on the coambient space as outlined in
Section~\ref{subsection.ambient},
\else
using the usual signed permutation representation~\cite{Bjorner_Brenti.2005,eriksson_eriksson.1998},
\fi
similar to the recursive sorting algorithm for type $A$ at the beginning of this section.
This is an explicit algorithm which achieves the results of Theorem~\ref{theorem.connected}
(but not necessarily in the most efficient way). 
This sorting algorithm actually contains all the ingredients for type
$C_n$ and $D_n$, since the Dynkin diagram of type $B_n$ contains both
kinds of endings. We have also verified by computer that explicit
recursive sorting algorithms exist for the exceptional types; the base
cases $B_2$, $B_3$, $C_2$, and $D_3$ can be worked out
explicitly. Details are available upon request.
  
Let $w$ be a permutation of type $B_n$ for $n\geq 4$. As before we identify $w$ with
$w(\rhoc)=\xc=(x_1,\ldots,x_n)$.  We can bring the maximal letter $z=n$ to any 
position, as $z$ or $-z$:
\begin{equation}
    \def\d{\ar@/^7em/[d]}
    \vcenter{
    \xymatrix@R=.2em{
      x_1 \ldots x_{k-1} \; z \; x_k \ldots x_{n-1} \; \d^{\pi_{n-1} \cdots \pi_k}\\
      x_1 \; x_2 \ldots x_{n-1} \; z  \qquad \quad \; \;  \; \d^{\pi_n}\\
      x_1 \; x_2 \ldots x_{n-1} \; -z  \qquad \;  \; \d^{\pi_2 \cdots \pi_{n-1}}\\
      x_1 \; -z \; x_1 \ldots x_{n-1}  \quad \; \; \; \; \d^{\pi_0}\\
      z \; -x_1 \; x_2 \ldots x_{n-1} \qquad \;\d^{\;\pi_{k-1}\cdots \pi_1}\\
      -x_1 \ldots x_{k-1} \; z \; x_k \ldots x_{n-1}
    }}
\end{equation}
In  particular, we can move $z$ to the left of $y=n-1$ (or $-z$ to the right of $-y$).
The pair $zy$ (or $-y-z$) can move around in a circle to any position by similar arguments
as above without disturbing any of the other letters, noting that if $zy$ are in the last two positions 
of $\xc$, then $\pi_n \pi_{n-1}\pi_n(\xc)$ contains $-y-z$ in the last two positions, and if $-y-z$ is in the 
first two positions of $\xc$, then $\pi_0(\xc)$ contains $zy$ in the first two positions.
  
Next suppose that $zy$ occupy the first two positions of $\xc$. We construct $\tilde{\pi}_0$ on such
$\xc$, which acts the same way as $\pi_0$, but on the last $n-2$ letters:
\begin{equation}
    \def\d{\ar@/^5.5em/[d]}
    \vcenter{
    \xymatrix@R=.2em{
      z \; y \; x_1 \; x_2 \cdots x_{n-2} \d^{\pi_2\pi_1\pi_3\pi_2}\\
      x_1 \; x_2 \; z \; y \cdots x_{n-2} \d^{\pi_0}\\
      x_1' \; x_2' \; z \; y \cdots x_{n-2}
    }}
\end{equation}
followed by the above circling to move $zy$ back to position 1 and 2.

\begin{problem}
  We had first proved a variant of
  Proposition~\ref{proposition.transitiveA} with the cycle
  $(1,\dots,n)$ and $\pi_1,\dots,\pi_n$ as operators. There, the
  sorting of a permutation $\sigma$ involves decomposing it
  recursively in terms of the following strong generating set of
  $\sg$ (as a permutation group):
  \begin{equation}
    \left( ((1,\dots,i)^k)_{k=0,\dots,i-1}\right) _{i=1,\dots n}.
  \end{equation}
  The sequence $(k_n,\dots,k_1)$ describing which power $k_i$ of
  $(1,\dots,i)$ is used for each base point $i$ is (essentially) the
  flag code of $\sigma$, as defined in~\cite{Adin_Roichman.2001}.

  Similar flag codes have been defined for types $B_n$, $C_n$, $D_n$,
  and even for general reflection
  groups~\cite{Adin_Brenti_Roichman.2005,Biagioli_Caselli.2004,Bagno_Biagioli.2007}.
  Do there exist related recursive sorting algorithms?
\end{problem}

\iflongversion
\subsection{Strong connectivity of crystals}
\label{subsection.crystal}

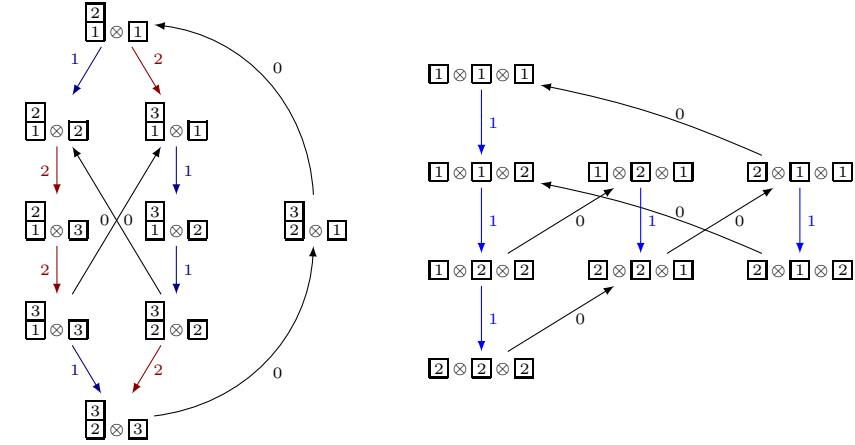
\begin{figure}
  \centerline{
    \begin{tikzpicture}[>=latex,join=bevel,scale=.5,baseline=(current bounding box.east)]
\tiny%
\node (N_1) at (100bp,300bp) [draw=none] {${\def\lr#1#2#3{\multicolumn{1}{#1@{\hspace{.6ex}}c@{\hspace{.6ex}}#2}{\raisebox{-.3ex}{$#3$}}}\raisebox{-.6ex}{$\begin{array}[b]{c}\cline{1-1}\lr{|}{|}{2}\\\cline{1-1}\lr{|}{|}{1}\\\cline{1-1}\end{array}$}} \otimes {\def\lr#1#2#3{\multicolumn{1}{#1@{\hspace{.6ex}}c@{\hspace{.6ex}}#2}{\raisebox{-.3ex}{$#3$}}}\raisebox{-.6ex}{$\begin{array}[b]{c}\cline{1-1}\lr{|}{|}{1}\\\cline{1-1}\end{array}$}}$};
  \node (N_2) at (55bp,225bp) [draw=none] {${\def\lr#1#2#3{\multicolumn{1}{#1@{\hspace{.6ex}}c@{\hspace{.6ex}}#2}{\raisebox{-.3ex}{$#3$}}}\raisebox{-.6ex}{$\begin{array}[b]{c}\cline{1-1}\lr{|}{|}{2}\\\cline{1-1}\lr{|}{|}{1}\\\cline{1-1}\end{array}$}} \otimes {\def\lr#1#2#3{\multicolumn{1}{#1@{\hspace{.6ex}}c@{\hspace{.6ex}}#2}{\raisebox{-.3ex}{$#3$}}}\raisebox{-.6ex}{$\begin{array}[b]{c}\cline{1-1}\lr{|}{|}{2}\\\cline{1-1}\end{array}$}}$};
  \node (N_3) at (145bp,225bp) [draw=none] {${\def\lr#1#2#3{\multicolumn{1}{#1@{\hspace{.6ex}}c@{\hspace{.6ex}}#2}{\raisebox{-.3ex}{$#3$}}}\raisebox{-.6ex}{$\begin{array}[b]{c}\cline{1-1}\lr{|}{|}{3}\\\cline{1-1}\lr{|}{|}{1}\\\cline{1-1}\end{array}$}} \otimes {\def\lr#1#2#3{\multicolumn{1}{#1@{\hspace{.6ex}}c@{\hspace{.6ex}}#2}{\raisebox{-.3ex}{$#3$}}}\raisebox{-.6ex}{$\begin{array}[b]{c}\cline{1-1}\lr{|}{|}{1}\\\cline{1-1}\end{array}$}}$};
  \node (N_4) at (55bp,150bp) [draw=none] {${\def\lr#1#2#3{\multicolumn{1}{#1@{\hspace{.6ex}}c@{\hspace{.6ex}}#2}{\raisebox{-.3ex}{$#3$}}}\raisebox{-.6ex}{$\begin{array}[b]{c}\cline{1-1}\lr{|}{|}{2}\\\cline{1-1}\lr{|}{|}{1}\\\cline{1-1}\end{array}$}} \otimes {\def\lr#1#2#3{\multicolumn{1}{#1@{\hspace{.6ex}}c@{\hspace{.6ex}}#2}{\raisebox{-.3ex}{$#3$}}}\raisebox{-.6ex}{$\begin{array}[b]{c}\cline{1-1}\lr{|}{|}{3}\\\cline{1-1}\end{array}$}}$};
  \node (N_5) at (145bp,150bp) [draw=none] {${\def\lr#1#2#3{\multicolumn{1}{#1@{\hspace{.6ex}}c@{\hspace{.6ex}}#2}{\raisebox{-.3ex}{$#3$}}}\raisebox{-.6ex}{$\begin{array}[b]{c}\cline{1-1}\lr{|}{|}{3}\\\cline{1-1}\lr{|}{|}{1}\\\cline{1-1}\end{array}$}} \otimes {\def\lr#1#2#3{\multicolumn{1}{#1@{\hspace{.6ex}}c@{\hspace{.6ex}}#2}{\raisebox{-.3ex}{$#3$}}}\raisebox{-.6ex}{$\begin{array}[b]{c}\cline{1-1}\lr{|}{|}{2}\\\cline{1-1}\end{array}$}}$};
  \node (N_6) at (55bp,75bp) [draw=none] {${\def\lr#1#2#3{\multicolumn{1}{#1@{\hspace{.6ex}}c@{\hspace{.6ex}}#2}{\raisebox{-.3ex}{$#3$}}}\raisebox{-.6ex}{$\begin{array}[b]{c}\cline{1-1}\lr{|}{|}{3}\\\cline{1-1}\lr{|}{|}{1}\\\cline{1-1}\end{array}$}} \otimes {\def\lr#1#2#3{\multicolumn{1}{#1@{\hspace{.6ex}}c@{\hspace{.6ex}}#2}{\raisebox{-.3ex}{$#3$}}}\raisebox{-.6ex}{$\begin{array}[b]{c}\cline{1-1}\lr{|}{|}{3}\\\cline{1-1}\end{array}$}}$};
  \node (N_7) at (145bp,75bp) [draw=none] {${\def\lr#1#2#3{\multicolumn{1}{#1@{\hspace{.6ex}}c@{\hspace{.6ex}}#2}{\raisebox{-.3ex}{$#3$}}}\raisebox{-.6ex}{$\begin{array}[b]{c}\cline{1-1}\lr{|}{|}{3}\\\cline{1-1}\lr{|}{|}{2}\\\cline{1-1}\end{array}$}} \otimes {\def\lr#1#2#3{\multicolumn{1}{#1@{\hspace{.6ex}}c@{\hspace{.6ex}}#2}{\raisebox{-.3ex}{$#3$}}}\raisebox{-.6ex}{$\begin{array}[b]{c}\cline{1-1}\lr{|}{|}{2}\\\cline{1-1}\end{array}$}}$};
  \node (N_8) at (250bp,150bp) [draw=none] {${\def\lr#1#2#3{\multicolumn{1}{#1@{\hspace{.6ex}}c@{\hspace{.6ex}}#2}{\raisebox{-.3ex}{$#3$}}}\raisebox{-.6ex}{$\begin{array}[b]{c}\cline{1-1}\lr{|}{|}{3}\\\cline{1-1}\lr{|}{|}{2}\\\cline{1-1}\end{array}$}} \otimes {\def\lr#1#2#3{\multicolumn{1}{#1@{\hspace{.6ex}}c@{\hspace{.6ex}}#2}{\raisebox{-.3ex}{$#3$}}}\raisebox{-.6ex}{$\begin{array}[b]{c}\cline{1-1}\lr{|}{|}{1}\\\cline{1-1}\end{array}$}}$};
  \node (N_9) at (100bp,0bp) [draw=none] {${\def\lr#1#2#3{\multicolumn{1}{#1@{\hspace{.6ex}}c@{\hspace{.6ex}}#2}{\raisebox{-.3ex}{$#3$}}}\raisebox{-.6ex}{$\begin{array}[b]{c}\cline{1-1}\lr{|}{|}{3}\\\cline{1-1}\lr{|}{|}{2}\\\cline{1-1}\end{array}$}} \otimes {\def\lr#1#2#3{\multicolumn{1}{#1@{\hspace{.6ex}}c@{\hspace{.6ex}}#2}{\raisebox{-.3ex}{$#3$}}}\raisebox{-.6ex}{$\begin{array}[b]{c}\cline{1-1}\lr{|}{|}{3}\\\cline{1-1}\end{array}$}}$};
  \draw [->,DarkBlue] (N_1) -- node [above left] {$1$} (N_2);
  \draw [->,DarkRed]  (N_1) -- node [above right]{$2$} (N_3);

  \draw [->,DarkRed]  (N_2) -- node [left] {$2$} (N_4);
  \draw [->,DarkBlue] (N_3) -- node [auto] {$1$} (N_5);

  \draw [->,DarkRed]  (N_4) -- node [left] {$2$} (N_6);
  \draw [->,DarkBlue] (N_5) -- node [auto] {$1$} (N_7);

  \draw [<-,    ] (N_3) -- node [left] {$0$} (N_6);
  \draw [<-,    ] (N_2) -- node [right] {$0$} (N_7);

  \draw [->,DarkBlue] (N_6) -- node [left]{$1$} (N_9);
  \draw [->,DarkRed]  (N_7) -- node [right] {$2$} (N_9);

  \draw [<-,    ] (N_1) to [bend left=40] node [auto] {$0$}  (N_8);
  \draw [<-,    ] (N_8) to [bend left=40] node [auto] {$0$} (N_9);
\end{tikzpicture} \qquad
    \begin{tikzpicture}[>=latex,join=bevel,scale=.5,baseline=(current bounding box.east)]
\tiny%
\node (N_1) at (80bp,232bp) [draw,draw=none] {${\def\lr#1#2#3{\multicolumn{1}{#1@{\hspace{.6ex}}c@{\hspace{.6ex}}#2}{\raisebox{-.3ex}{$#3$}}}\raisebox{-.6ex}{$\begin{array}[b]{c}\cline{1-1}\lr{|}{|}{1}\\\cline{1-1}\end{array}$}} \otimes {\def\lr#1#2#3{\multicolumn{1}{#1@{\hspace{.6ex}}c@{\hspace{.6ex}}#2}{\raisebox{-.3ex}{$#3$}}}\raisebox{-.6ex}{$\begin{array}[b]{c}\cline{1-1}\lr{|}{|}{1}\\\cline{1-1}\end{array}$}} \otimes {\def\lr#1#2#3{\multicolumn{1}{#1@{\hspace{.6ex}}c@{\hspace{.6ex}}#2}{\raisebox{-.3ex}{$#3$}}}\raisebox{-.6ex}{$\begin{array}[b]{c}\cline{1-1}\lr{|}{|}{1}\\\cline{1-1}\end{array}$}}$};
  \node (N_2) at (80bp,158bp) [draw,draw=none] {${\def\lr#1#2#3{\multicolumn{1}{#1@{\hspace{.6ex}}c@{\hspace{.6ex}}#2}{\raisebox{-.3ex}{$#3$}}}\raisebox{-.6ex}{$\begin{array}[b]{c}\cline{1-1}\lr{|}{|}{1}\\\cline{1-1}\end{array}$}} \otimes {\def\lr#1#2#3{\multicolumn{1}{#1@{\hspace{.6ex}}c@{\hspace{.6ex}}#2}{\raisebox{-.3ex}{$#3$}}}\raisebox{-.6ex}{$\begin{array}[b]{c}\cline{1-1}\lr{|}{|}{1}\\\cline{1-1}\end{array}$}} \otimes {\def\lr#1#2#3{\multicolumn{1}{#1@{\hspace{.6ex}}c@{\hspace{.6ex}}#2}{\raisebox{-.3ex}{$#3$}}}\raisebox{-.6ex}{$\begin{array}[b]{c}\cline{1-1}\lr{|}{|}{2}\\\cline{1-1}\end{array}$}}$};
  \node (N_3) at (80bp,84bp) [draw,draw=none] {${\def\lr#1#2#3{\multicolumn{1}{#1@{\hspace{.6ex}}c@{\hspace{.6ex}}#2}{\raisebox{-.3ex}{$#3$}}}\raisebox{-.6ex}{$\begin{array}[b]{c}\cline{1-1}\lr{|}{|}{1}\\\cline{1-1}\end{array}$}} \otimes {\def\lr#1#2#3{\multicolumn{1}{#1@{\hspace{.6ex}}c@{\hspace{.6ex}}#2}{\raisebox{-.3ex}{$#3$}}}\raisebox{-.6ex}{$\begin{array}[b]{c}\cline{1-1}\lr{|}{|}{2}\\\cline{1-1}\end{array}$}} \otimes {\def\lr#1#2#3{\multicolumn{1}{#1@{\hspace{.6ex}}c@{\hspace{.6ex}}#2}{\raisebox{-.3ex}{$#3$}}}\raisebox{-.6ex}{$\begin{array}[b]{c}\cline{1-1}\lr{|}{|}{2}\\\cline{1-1}\end{array}$}}$};
  \node (N_8) at (80bp,10bp) [draw,draw=none] {${\def\lr#1#2#3{\multicolumn{1}{#1@{\hspace{.6ex}}c@{\hspace{.6ex}}#2}{\raisebox{-.3ex}{$#3$}}}\raisebox{-.6ex}{$\begin{array}[b]{c}\cline{1-1}\lr{|}{|}{2}\\\cline{1-1}\end{array}$}} \otimes {\def\lr#1#2#3{\multicolumn{1}{#1@{\hspace{.6ex}}c@{\hspace{.6ex}}#2}{\raisebox{-.3ex}{$#3$}}}\raisebox{-.6ex}{$\begin{array}[b]{c}\cline{1-1}\lr{|}{|}{2}\\\cline{1-1}\end{array}$}} \otimes {\def\lr#1#2#3{\multicolumn{1}{#1@{\hspace{.6ex}}c@{\hspace{.6ex}}#2}{\raisebox{-.3ex}{$#3$}}}\raisebox{-.6ex}{$\begin{array}[b]{c}\cline{1-1}\lr{|}{|}{2}\\\cline{1-1}\end{array}$}}$};

  \node (N_4) at (200bp,158bp) [draw,draw=none] {${\def\lr#1#2#3{\multicolumn{1}{#1@{\hspace{.6ex}}c@{\hspace{.6ex}}#2}{\raisebox{-.3ex}{$#3$}}}\raisebox{-.6ex}{$\begin{array}[b]{c}\cline{1-1}\lr{|}{|}{1}\\\cline{1-1}\end{array}$}} \otimes {\def\lr#1#2#3{\multicolumn{1}{#1@{\hspace{.6ex}}c@{\hspace{.6ex}}#2}{\raisebox{-.3ex}{$#3$}}}\raisebox{-.6ex}{$\begin{array}[b]{c}\cline{1-1}\lr{|}{|}{2}\\\cline{1-1}\end{array}$}} \otimes {\def\lr#1#2#3{\multicolumn{1}{#1@{\hspace{.6ex}}c@{\hspace{.6ex}}#2}{\raisebox{-.3ex}{$#3$}}}\raisebox{-.6ex}{$\begin{array}[b]{c}\cline{1-1}\lr{|}{|}{1}\\\cline{1-1}\end{array}$}}$};
  \node (N_5) at (200bp,84bp) [draw,draw=none] {${\def\lr#1#2#3{\multicolumn{1}{#1@{\hspace{.6ex}}c@{\hspace{.6ex}}#2}{\raisebox{-.3ex}{$#3$}}}\raisebox{-.6ex}{$\begin{array}[b]{c}\cline{1-1}\lr{|}{|}{2}\\\cline{1-1}\end{array}$}} \otimes {\def\lr#1#2#3{\multicolumn{1}{#1@{\hspace{.6ex}}c@{\hspace{.6ex}}#2}{\raisebox{-.3ex}{$#3$}}}\raisebox{-.6ex}{$\begin{array}[b]{c}\cline{1-1}\lr{|}{|}{2}\\\cline{1-1}\end{array}$}} \otimes {\def\lr#1#2#3{\multicolumn{1}{#1@{\hspace{.6ex}}c@{\hspace{.6ex}}#2}{\raisebox{-.3ex}{$#3$}}}\raisebox{-.6ex}{$\begin{array}[b]{c}\cline{1-1}\lr{|}{|}{1}\\\cline{1-1}\end{array}$}}$};

  \node (N_6) at (320bp,158bp) [draw,draw=none] {${\def\lr#1#2#3{\multicolumn{1}{#1@{\hspace{.6ex}}c@{\hspace{.6ex}}#2}{\raisebox{-.3ex}{$#3$}}}\raisebox{-.6ex}{$\begin{array}[b]{c}\cline{1-1}\lr{|}{|}{2}\\\cline{1-1}\end{array}$}} \otimes {\def\lr#1#2#3{\multicolumn{1}{#1@{\hspace{.6ex}}c@{\hspace{.6ex}}#2}{\raisebox{-.3ex}{$#3$}}}\raisebox{-.6ex}{$\begin{array}[b]{c}\cline{1-1}\lr{|}{|}{1}\\\cline{1-1}\end{array}$}} \otimes {\def\lr#1#2#3{\multicolumn{1}{#1@{\hspace{.6ex}}c@{\hspace{.6ex}}#2}{\raisebox{-.3ex}{$#3$}}}\raisebox{-.6ex}{$\begin{array}[b]{c}\cline{1-1}\lr{|}{|}{1}\\\cline{1-1}\end{array}$}}$};
  \node (N_7) at (320bp,84bp) [draw,draw=none] {${\def\lr#1#2#3{\multicolumn{1}{#1@{\hspace{.6ex}}c@{\hspace{.6ex}}#2}{\raisebox{-.3ex}{$#3$}}}\raisebox{-.6ex}{$\begin{array}[b]{c}\cline{1-1}\lr{|}{|}{2}\\\cline{1-1}\end{array}$}} \otimes {\def\lr#1#2#3{\multicolumn{1}{#1@{\hspace{.6ex}}c@{\hspace{.6ex}}#2}{\raisebox{-.3ex}{$#3$}}}\raisebox{-.6ex}{$\begin{array}[b]{c}\cline{1-1}\lr{|}{|}{1}\\\cline{1-1}\end{array}$}} \otimes {\def\lr#1#2#3{\multicolumn{1}{#1@{\hspace{.6ex}}c@{\hspace{.6ex}}#2}{\raisebox{-.3ex}{$#3$}}}\raisebox{-.6ex}{$\begin{array}[b]{c}\cline{1-1}\lr{|}{|}{2}\\\cline{1-1}\end{array}$}}$};

  \draw [->,blue] (N_1) -- node [auto] {$1$} (N_2);
  \draw [->,blue] (N_2) -- node [auto] {$1$} (N_3);
  \draw [->,blue] (N_3) -- node [auto] {$1$} (N_8);

  \draw [->,blue] (N_4) -- node [auto] {$1$} (N_5);

  \draw [->,blue] (N_6) -- node [auto] {$1$} (N_7);

  \draw [<-]      (N_1) to [bend left=6] node [right=2mm] {$0$} (N_6);
  \draw [<-]      (N_2) to [bend left=6] node [right=2mm] {$0$} (N_7);

  \draw [<-]      (N_4) -- node [right=1mm] {$0$} (N_3);
  \draw [<-]      (N_6) -- node [right=1mm] {$0$} (N_5);
  \draw [<-]      (N_5) -- node [right=1mm] {$0$} (N_8);
\end{tikzpicture}}
   \caption{
     \textbf{Left}: Crystal $B^{2,1}\otimes B^{1,1}$ of type $A_2^{(1)}$. By
     contraction of all $i$-strings to a single edge $i$, one recovers
     the left most graph of Figure~\ref{figure.antisorting}.
     \textbf{Right}: Crystal $(B^{1,1})^{\otimes 3}$ of type $A_1^{(1)}$.}
  \label{figure.crystal}
\end{figure}

Crystal bases are combinatorial bases of modules of quantum algebras $U_q(\geh)$ as the
parameter $q$ tends to zero. They consist of a non-empty set $B$ together with 
raising and lowering operators $e_i$ and $f_i$ for $i\in I$ from $B$ to $B\cup \{0\}$ 
and a weight function $\wt:B\to \weightspace$. For more information on crystal theory 
see~\cite{Hong_Kang.2002}. Of particular interest are crystals coming from finite-dimensional
affine $U_q(\geh)$-modules, where $\geh$ is an affine Kac-Moody algebra.
These crystals are not highest weight. In this section we deduce from the transitivity
of the level 0 action of the 0-Hecke algebra on $\clW$ of Theorem~\ref{theorem.connected}
that these finite-dimensional affine crystals are strongly connected;
that is, any two elements $b,b'\in B$ can be connected via a sequence of operators $f_i$:
$b'=f_{i_1}\cdots f_{i_r}(b)$ for $i_j\in I$.

There is an action of the Weyl group on any finite affine crystal $B$ defined by
\begin{equation}
  s_i(b) = \begin{cases}
    f_i^{\ip{\alpha_i^\vee,\wt(b)}}(b) & \text{if $\ip{\alpha_i^\vee,\wt(b)} > 0$,}\\
    e_i^{-\ip{\alpha_i^\vee,\wt(b)}}(b) & \text{if $\ip{\alpha_i^\vee,\wt(b)} \le 0$,}
    \end{cases}
\end{equation}
where $b\in B$ and $i\in I$. This action is compatible with the weights, that is,
$s_i(\wt(b)) = \wt(s_i(b))$. In particular we also have $\wt(\pi_i(b)) = \pi_i(\wt(b))$, where 
\begin{equation}
\label{equation.piOnCrystal}
  \pi_i(b) := \begin{cases}
    f_i^{\ip{\alpha_i^\vee,\wt(b)}}(b) & \text{if $\ip{\alpha_i^\vee,\wt(b)} > 0$,}\\
    b & \text{if $\ip{\alpha_i^\vee,\wt(b)} \le 0$.}
    \end{cases}
\end{equation}

\begin{remark}
\label{remark.piOnCrystal}
Comparing~\eqref{eq:pi} and~\eqref{equation.piOnCrystal}, it is clear that 
if a sequence $i_1,\dots,i_r$ is such that at each step in
$\pisequence(\wt(b))$ the operator $\pi_i$ acts as $s_i$, then the same
holds in $\pisequence(b)$.
\end{remark}

\begin{theorem}
\label{theorem.strongCrystal}
  Let $B$ be a finite connected affine crystal. Then $B$ is strongly
  connected.
\end{theorem}

\begin{proof}
  It is sufficient to prove that if $x$ and $y$ in $B$ are in the same
  $i_0$-string with $y=f_{i_0}^a(x)$ for some $i_0\in I$ and $a>0$, then
  there is an $f$-path from $y$ to $x$. Using finiteness, we may
  further assume without loss of generality that
  $y=s_{i_0}(x)=\pi_{i_0}(x)$ (moving for example $x$ and $y$ to
  respectively to top and bottom of the string).

  By Theorem~\ref{theorem.connected}, there exists a sequence
  $i_1,\dots,i_r$ such that
  $\pisequence(\wt(y))=\wt(x)$. Choose such a sequence
  of minimal length, so that each $\pi_{i_j}$ above acts as
  $s_{i_j}$. Consider $p:=\pisequence \pi_{i_0}$, and
  $\s:=s_{i_r} \cdots s_{i_0}$. Then,
  $p(\wt(x))=\s(\wt(x))=\wt(x)$. Now, $p(x)$ might not be $x$,
  but by Remark~\ref{remark.piOnCrystal} we may apply $p$ repeatedly and still have
  $p^k(x)=\s^k(x)$. Since the crystal is finite,
  eventually we will have $p^k(x)=\s^k(x)=x$. Since any application
  of $\pi_i$ results from a sequence of applications of $f_i$,
  this proves the existence of an $f$-path from $y$ back to $x$.
\end{proof}

Theorem~\ref{theorem.strongCrystal} is equivalent to~\cite[Theorem 3.37]{Kashiwara.2008}.

\begin{remark}
As noted in the proof of Theorem~\ref{theorem.strongCrystal}, the action of the affine Weyl group 
on a crystal is \emph{not} necessarily the level $0$ action: only a power of $p$ maps a given crystal
element $x$ to itself $p^k(x)=x$. Take for example $x=\begin{array}{|c|} \hline 1\\ \hline \end{array}
\otimes \begin{array}{|c|} \hline 1\\ \hline \end{array} \otimes \begin{array}{|c|} \hline 2\\ \hline 
\end{array}$ in $(B^{1,1})^{\otimes 3}$ of type $A_1^{(1)}$, where $B^{r,s}$ denotes
a Kirillov--Reshetikhin crystal. Then for $p=s_0 s_1$ we have 
$p(\wt(x))=\wt(x)$, but only $p^3(x)=x$ as can be seen from Figure~\ref{figure.crystal}.
\end{remark}

\begin{remark}
  Interpreting the $\pi_i$'s as Demazure operators,
  Theorem~\ref{theorem.connected} is related to properties of affine
  crystals. Let $\geh$ be an affine Kac--Moody algebra, $\W$ the corresponding 
  affine Weyl group, and $B^{r,s}$ a Kirillov--Reshetikhin crystal of type 
  $\geh$~\cite{Hatayama_Kuniba_Okado_Takagi_Tsuboi.2001,okado_schilling.2008}.
  Consider the affine crystal $B := B^{n,1} \otimes B^{n-1,1} \otimes
  \cdots \otimes B^{1,1}$, and define the Demazure operators on $b\in
  B$ as in~\cite{Kashiwara.1993}:
  \begin{equation}
  	\D_i(b) = \begin{cases}
        \phantom{-}
        \sum_{0\le k \le \phantom{-}\ip{\alpha_i^\vee,\wt(b)}} f_i^k(b) &  
	\text{if  $\ip{\alpha_i^\vee,\wt(b)} \ge 0$\,,}\\
 	-\sum_{1\le k \le -\ip{\alpha_i^\vee,\wt(b)}} e_i^k(b) &  
	\text{if  $\ip{\alpha_i^\vee,\wt(b)}< 0$\,.}
	\end{cases}
  \end{equation}
  Let $\Lambda_i$ be the fundamental weights of $\geh$, and take for $u_i$ the
  unique element in $B^{i,1}$ of weight $\Lambda_i-\Lambda_0$.  Then,
  the transitivity of the action of $\heckeW{0}$ on $\clW$ is
  closely related to the strong connectivity of the graph
  generated by $\D_0,\ldots,\D_n$ acting on $u_n \otimes \cdots
  \otimes u_1$~\cite{kashiwara.2002,fourier_schilling_shimozono.2007}, see 
  Figure~\ref{figure.crystal}.
\end{remark}
\fi

\section{Hecke group algebras as quotients of affine Hecke algebras}
\label{section.quotient}

We are now in the position to state the main theorem of this paper.
Let $\W$ be an affine Weyl group and $\heckeW{q_1,q_2}$ its Hecke algebra. Let 
$\clW$ be the associated
finite Weyl group, and $\heckeWW[\clW]$ its Hecke group algebra.
Then the level $0$-representation
\begin{equation}
\label{eq:cl}
  \cl: 
  \begin{cases}
    \heckeW{q_1,q_2} & \to \End(\kclW)\\
    T_i    & \mapsto (q_1+q_2)\pi_i - q_1 s_i
  \end{cases}
\end{equation}
actually defines a morphism from $\heckeW{q_1,q_2}$ to
$\heckeWW[\clW]$. (Note that $\pi_\alpha$ and in particular $\pi_0$ is indeed an element of 
$\heckeWW[\clW]$: it can be written as $\pi_\alpha = \s\pi_i\s^{-1}$ where $\s$ is an
element of $\clW$ conjugating $\alpha$ to some simple root $\alpha_i$.)
When the Dynkin diagram has special automorphisms $\Omega$, this morphism can
be extended to the extended affine Hecke algebra by sending the
special Dynkin diagram automorphisms to the corresponding element of
the finite Weyl group $\clW$.

\begin{theorem}
  \label{theorem.quotient}
  Let $\W$ be an affine Weyl group.  Except when $q_1+q_2=0$
  (and possibly when $q:=\q$ is a $k$th root of unity
  with $k\leq 2\height(\theta^\vee)$), the morphism $\cl:
  \heckeW{q_1,q_2}\to \heckeWW[\clW]$ is surjective and makes the
  Hecke group algebra $\heckeWW[\clW]$ into a quotient of the
  affine Hecke algebra $\heckeW{q_1,q_2}$.
\end{theorem}

\begin{proof}
  Here we outline the proof which relies on material in the next two sections.

  When $q_1+q_2=0$, the image of $\cl$ is obviously $\clWa$ (or
  just $\{0\}$ if $q_1=q_2=0$); so the morphism is not surjective.

  If $q_1=0$ and $q_2\ne0$, this is exactly
  Corollary~\ref{corollary.piGenerators} below. If $q_2=0$ and $q_1\ne0$,
  then $\cl(\oT_i) = q_1 \opi_i$, and by symmetry, we can also use
  Corollary~\ref{corollary.piGenerators}.
  The theorem follows right away for all values of $q=-q_1/q_2$ but a finite
  number using a standard specialization argument: take $q$ formal,
  and consider the family $B_q$ obtained from $B$ by replacing each
  $\pi_i$ by $(1-q) \pi_i + q s_i$. This family has polynomial
  coefficients when expressed in terms of the basis $\{ \s\pi_\t
  \suchthat \Des(\s) \cap \Rec(\t) = \emptyset\}$ of $\heckeWW[\clW]$. Its
  determinant is a polynomial in $q$ with a non-zero constant since
  $B_0$ is a basis.  Thus it vanishes for at most a finite number of
  values of $q$.
  
  Theorem~\ref{theorem.quotientQualibrated} below allows to further reduce
  the possible inappropriate values of $q$ to $k$th roots of unity
  with $k$ small. Note however that
  Theorem~\ref{theorem.quotientQualibrated} does not apply at $q_1=0$
  or $q_2=0$.
\end{proof}

Theorem~\ref{theorem.quotient} raises immediately the following
problem, currently under investigation together with Nicolas Borie.
\begin{problem}
  Determine for which roots of unity $q$ the morphism $\cl$ is not
  surjective.
\end{problem}

\section{Alternative generators for Hecke group algebras}
\label{section.generators}

In this section we show
that the Hecke group algebra can be entirely generated by $\pi_0,\pi_1,\dots,\pi_n$.

\begin{proposition}
  \label{proposition.piGenerators}
  Let $\clW$ be a finite Coxeter group, and $\setofroots$ be a set of
  roots of $\clW$ such that the associated projections $\{\pi_\alpha
  \suchthat \alpha\in\setofroots\}$ act transitively on $\clW$. Then,
  the Hecke group algebra $\heckeWW[\clW]$ is generated as an algebra by
  $\{\pi_\alpha \suchthat \alpha\in\setofroots\}$.
\end{proposition}

\begin{proof}
  First note that $\pi_\alpha$ is indeed an element of $\heckeWW[\clW]$: it
  can be written as $\pi_\alpha = \s\pi_i\s^{-1}$ where $\s$ is an
  element of $\clW$ conjugating $\alpha$ to some simple root $\alpha_i$.
  In Proposition~\ref{proposition.piBasis} below, we exhibit a
  sufficiently large family of operators which are linearly
  independent, because they display the same triangularity property as
  the basis $\{ \s\pi_\t \suchthat \Des(\s) \cap \Rec(\t) =
  \emptyset\}$ of $\heckeWW[\clW]$ (see Lemma~3.8
  of~\cite{Hivert_Thiery.HeckeGroup.2007}).
\end{proof}

\begin{corollary}
  \label{corollary.piGenerators}
  Let $\W$ be an affine Weyl group, $\clW$ be the associated
  finite Weyl group, and $\pi_0,\dots,\pi_n$ be the projections
  associated to the roots $\cl(\alpha_0),\dots,\cl(\alpha_n)$ of the finite Weyl group. 
  Then, the Hecke group algebra
  $\heckeWW[\clW]$ is generated as an algebra by $\pi_0,\dots,\pi_n$.

  Alternatively, $\pi_0$ may be replaced by any $\Omega\in W$ mapping
  $\alpha_0$ to some simple root, typically one induced by some
  special Dynkin diagram automorphism.
\end{corollary}

Let $\s\in \clW$. An $\setofroots$-\textit{reduced word} for $\s$ is a
word $i_1,\dots,i_r$ of minimal length such that $i_j\in\setofroots$
and $\s^{-1}.\pi_{i_1}\dots\pi_{i_r} = 1$. Since the $\{\pi_\alpha
\suchthat \alpha\in\setofroots\}$ acts transitively on $\clW$, there
always exists such an $\setofroots$-reduced word, and we choose once
for all one of them for each $\s$. More generally, for a right coset
$\s \clW_J$, we choose an $\setofroots$-reduced word $i_1,\dots,i_r$ of
minimal length such that there exists $\nu\in \clW_J\s^{-1}$ and $\mu\in
\clW_J$ with $\nu.\pi_{i_1}\dots\pi_{i_r} = \mu$. 

\begin{example}
  In type $C_2$, the word $0,1,2,0,1,0$ is $\setofroots$-reduced for
  $\Wmax=\Wmax^{-1}=(\underline1,\underline2)$, where we write 
  $\underline{1}$ and $\underline{2}$ for $-1$ and $-2$ (see
  Figure~\ref{figure.antisorting}).

  In type $A_3$ the word $1,0$ is $\setofroots$-reduced for
  $4123\clW_{\{1,3\}} $.  Here $w=4123$, $\nu=w^{-1}=2341$, and
  $\mu=1243$. Looking at $\clW_J$ left-cosets is the Coxeter equivalent
  to looking at words with repetitions: we may think of left
  $\clW_{\{1,3\}}$-cosets as identifying the values $1,2$ and $3,4$, and
  represent $\clW_{\{1,3\}}\s^{-1}=2341$ by the word $1331$; this word
  gets sorted by $\pi_1\pi_0$ to $1133$ which represents
  $\clW_{\{1,3\}}$.
\end{example}
Setting $\pim_i := (\pi_i-1)$, define the operator $\pim_{\s
  \clW_J}:=\pim_{i_1}\cdots\pim_{i_r}$ where $i_1,\dots,i_r$ is the
chosen $\setofroots$-reduced word. The operator may actually depend on
the choice of the $\setofroots$-reduced word, but this is irrelevant
for our purpose.
\begin{proposition}
  \label{proposition.piBasis}
  The following family forms a basis for $\heckeWW[\clW]$:
  \begin{equation}
    B := \{ \pim_{\s \clW_{\Rec(\t)}} \pi_\t \suchthat \Rec(\s) \cap \Des(\t) =
    \emptyset \} \,.
  \end{equation}
\end{proposition}

\begin{proof}
The number of elements of $B$ is the same as the dimension of the Hecke group
algebra by Section~\ref{subsection.heckeGroupAlgebras}. Corollary~\ref{corollary.triangularity}
shows that the elements in $B$ are linearly independent. This proves the claim.
\end{proof}

\iflongversion
\begin{lemma}
  \label{lemma.sortingOnCosets}%
  Let $\s \clW_J$ be a right coset in $\clW$, and $i_1,\dots,i_r$ be the
  corresponding $\setofroots$-reduced word. Set $\s'=s_{i_1}\cdots s_{i_r}$. Then,
  $\pi_{i_1}\cdots\pi_{i_r}$ restricted to $\clW_J\s^{-1}$ acts by
  right multiplication by $\s'$. In particular, it induces a bijection
  from $\clW_J\s^{-1}$ to $\clW_J$.
\end{lemma}
\begin{proof}
  Take $\nu$ in $\clW_J\s^{-1}$ such that $\nu.\pi_{i_1}\dots\pi_{i_r}\in
  \clW_J$. By minimality of the $\setofroots$-reduced word, no $\pi_i$ acts trivially,
  so $\nu.\pi_{i_1}\dots\pi_{i_r}=\nu\s'$. Furthermore,
  $\pi_{i_1}\cdots\pi_{i_r}$ is in $\heckeWW[\clW]$ and thus preserves
  left-antisymmetries. Taking $i\in I$, this implies that
  $(s_i\nu).\pi_{i_1}\dots\pi_{i_r}$ is either $s_i\nu\s'$ or
  $\nu\s'$. By minimality of the $\setofroots$-reduced word, the latter case is
  impossible: indeed if any of the $\pi_{i_j}$ acts trivially we get a
  strictly shorter $\setofroots$-reduced word from $s_i\nu\in \clW_J\s^{-1}$ to
  $\nu\s'\in \clW_J$.  Applying transitivity, we get that
  $\pi_{i_1}\cdots\pi_{i_r}$ acts by multiplication by $\s'$ on
  $\clW_J\s^{-1}$.
\end{proof}

Let $<$ be any linear extension of the right Bruhat order on $\clW$.
Given an endomorphism $f$ of $\kclW$, we order the rows and columns of
its matrix $M_f:=[f_{\mu\nu}:= f(\nu)_{|\mu}]$ according to $<$
(beware that, the action being on the right, $M_{fg} = M_g M_f$).
Denote by
$\init(f) := \min\{ \mu \suchthat \exists \nu, f_{\mu\nu}\ne 0\}$ the
index of the first non-zero row of $M_f$.
\begin{lemma}
  \label{lemma.coninjectivity}
  Let $f:=\pim_{\s \clW_J}$. Then, for any $\mu\in \clW$, there exists a unique
  $\nu\in \clW$ such that the coefficient $f_{\mu\nu}$ is non-zero; this
  coefficient is either $1$ or $-1$ (in other words, $f$ is the
  transpose of a signed-monoidal application).

  In particular, if $\mu\in \clW_J$ then $\nu$ belongs to $\clW_J\s^{-1}$,
  and $f_{\mu\nu} = 1$.
\end{lemma}
\begin{proof}
  This is clear if $f=\pim_J$; here is for example the matrix of
  $\pim_1$ in type $A_1$:
  \begin{equation}
    \begin{pmatrix}
      -1 & 0\\
      1  & 0
    \end{pmatrix}\,.
  \end{equation}
  By products, this extends to any $f$.

  Take now $\mu\in \clW_J$. Using
  Lemma~\ref{lemma.sortingOnCosets}, let $\nu$ be 
  the unique element in $\clW_J\s^{-1}$ such that
  $\nu.\pi_{i_1}\dots\pi_{i_r}=\mu$.  By minimality of the
  $\setofroots$-reduced word, $\mu$ cannot occur in any other term of the
  expansion of
  \begin{equation}
    \nu.\pim_{i_1}\dots\pim_{i_r} = \nu.(\pi_{i_1}-1)\dots(\pi_{i_r}-1)\,.
  \end{equation}
  Therefore, $f_{\mu\nu}=1$, and $f_{\mu\nu'}=0$ for $\nu'\ne \nu$.
\end{proof}

We get as a corollary that the basis $B$ is triangular.

\begin{corollary}
  \label{corollary.triangularity}
  Let $f:=\pim_{\s\clW_{\Rec(\t)}} \pi_\t$ in $B$. Then, $\init(f)=\t$, and
  \begin{equation}
    f_{\t\nu} =
    \begin{cases}
      1 & \text{ if } \nu\in \clW_{\Rec(\t)} \s^{-1}\,,\\
      0 & \text{ otherwise.}
    \end{cases}
  \end{equation}
\end{corollary}
\fi

\section{Hecke group algebras and principal series representations
  of affine Hecke algebras}
\label{section.calibrated}

Let $t: Y \to \CC^*$ be a character of the multiplicative group
$Y$ (or equivalently of the additive group $\corootlattice$). It
induces a representation $M(t) := t \uparrow_{\CC[Y]}^{\heckeW{q_1,q_2}}$ 
called \emph{principal series representation} of the affine Hecke algebra
$\heckeW{q_1,q_2}$. Since
$\heckeW{q_1,q_2}=\heckeW[\clW]{q_1,q_2}\otimes \CC[Y]$, this
representation is of dimension $|\clW|$. When $t$ is regular, the
representation is \emph{calibrated}: it admits a basis $(E_w)_{w\in \clW}$ 
which diagonalizes the action of $Y$ with a distinct character
$wt$ on each $E_w$. This basis can be constructed explicitly by means
of intertwining operators $\tau_i$ which skew commute with the
elements of $Y$. We refer to~\cite[Section 2.5]{Ram.2003} for details. Note also
that the construction of the $\tau_i$ operators by deformation of the
$T_i$ is reminiscent of Yang-Baxter
graphs~\cite[\S~10.7]{Lascoux.2003.CBMS}, in which $t$ corresponds
to a choice of \emph{spectral parameters}.

The main result of this section is that for $q_1,q_2\ne 0$ and $q$
not a root of unity, there exists a suitable character $t$, such that
the level $0$ representation of the affine Hecke algebra is isomorphic
to the principal series representation $M(t)$
(Theorem~\ref{theorem.isomorphismCalibrated}), and to deduce that the
morphism $\cl:\heckeW{q_1,q_2}\mapsto \heckeWW[\clW]$ is
surjective (Theorem~\ref{theorem.quotientQualibrated}).

\begin{theorem}
  \label{theorem.isomorphismCalibrated}
  Assume $q_1,q_2$ are such that $q_1, q_2\ne 0$ and
  $q:=\q$ is not a $k$th root of unity with $k\leq
  2\height(\theta^\vee)$. Then, the level $0$ representation of the
  affine Hecke algebra $\heckeW{q_1,q_2}$ is isomorphic to the
  principal series representation $M(t)$ for the character
  $t:Y^{\lc}\mapsto q^{-\height(\lc)}$.
\end{theorem}
Note that $t(Y^{\coroot_i})=q^{-1}$ for any simple coroot. By a result
of Kato~\cite[Theorem 2.2]{Kato.1981} (see also~\cite[Theorem 2.12
(c)]{Ram.2003}) one sees right away that $M(t)$ is not
irreducible. Note also that this is, up to inversion, the same
character as for the action of $\CC[Y]$ on the constant Macdonald
polynomial $\mathbf{1}$ \cite[Equation (3.4)]{Ram.2008}.

\begin{proof}
  In the upcoming Lemma~\ref{lemma.w0eigenvector}, we prove that $\Wmax$
  is an eigenvector for the character $t$, and check that $t$ is
  regular (that is the orbit $\clW t$ of $t$ is of size
  $|\clW|$). We then mimic~\cite{Ram.2003} and use the
  intertwining operators to explicitly diagonalize the action of $Y$
  on $\kclW$ in
  Proposition~\ref{proposition.eigenvectors}. Although this is more
  than strictly necessary to prove the desired isomorphism, the
  results will be useful for the subsequent
  Theorem~\ref{theorem.quotientQualibrated}.
\end{proof}

\iflongversion
\begin{lemma}
  \label{lemma.alcoveWalkTiOnW0}
  Let $i_1,\dots,i_r$ be an alcove walk from the fundamental alcove,
  and $\epsilon_1,\dots,\epsilon_r$ as defined in
  Section~\ref{subsection.affine.weyl.groups}. Then,
  \begin{equation}
    \Wmax.T_{i_1}^{\epsilon_1}\cdots T_{i_r}^{\epsilon_r} =
    q_2^{\epsilon_1+\dots+\epsilon_r}\Wmax s_{i_1}\cdots s_{i_r}\;.
  \end{equation}
\end{lemma}
\begin{proof}
  Take $w\in \clW$, and $i\in \{0,\dots,n\}$. If $i$ is not a
  descent of $w$, then, using~\eqref{eq:cl}:
  \begin{equation}
    w.T_i = w.\quotientmap(T_i) = w.\left((q_1+q_2)\pi_i -q_1 s_i\right)
    = w ((q_1+q_2) s_i - q_1 s_i) = q_2 w s_i\,.
  \end{equation}
  Inverting this equation yields that, when $i$ is a descent of
  $w$, $w.T_i^{-1} = q_2^{-1} w s_i$.
  
  We conclude by induction since $\epsilon_k=1$ if and only if $i_k$
  is a descent of $w_{k-1} = s_{i_1}\cdots s_{i_{k-1}}$
  (cf. Remark~\ref{remark.positiveCrossing}), that is not a descent of
  $\Wmax s_{i_1}\cdots s_{i_{k-1}}$.
\end{proof}
\fi

\begin{proposition}
  \label{proposition.eigenvectors}
  Assuming the same conditions as in
  Theorem~\ref{theorem.isomorphismCalibrated}, there exists a basis
  $(E_w)_{w\in \clW}$ of $\kclW$ which diagonalizes
  simultaneously all $Y^\lavee$:
  \begin{equation}
    E_w. Y^\lavee = (wt)(Y^\lavee) E_w\,,
  \end{equation}
  where $(wt)(Y^\lavee):=q^{-\height(w(\lavee))}$.
\iflongversion
  In particular, the eigenvalue for $Y^{\lavee}$ on $E_w$ is $q^{-1}$ if
  and only if $w(\lc)$ is a simple coroot.
\fi
\end{proposition}
\iflongversion
Note that acting with $\oT_i$'s instead, or equivalently defining
$T_i$'s in term of the operators $\opi_i$'s would allow to revert the
picture and use $1$ as initial eigenvector instead of $\Wmax$. We also
get the following side result on the Hecke group algebra.
\fi

\begin{proof}
  First note that $t$ is regular; indeed, $\clrho$ is regular, and $q$
  is not a $k$th root of unity with $k$ too small, so one can use
  \begin{equation}
    (wt)(Y^{\alpha_i^\vee})=q^{-\height(w(\alpha_i^\vee))} =
    q^{-\ip{\alpha_i^\vee,w^{-1}(\clrho)}}
  \end{equation}
  to recover the coordinates of $w(\clrho)$ on each $i$th fundamental
  weight. For the same reason, $(wt)(Y^{\alpha_i^\vee})$ is never $1$.

  We first prove in Lemma~\ref{lemma.w0eigenvector} that $E_{1}=\Wmax$
  is an appropriate eigenvector, and then define intertwining
  operators $\tau_i$ to construct the other $E_w$
  (Lemma~\ref{lemma.eigenvectors}).
\end{proof}

\FIXME{Should we move this later, just before theorem 7.7?}
\begin{corollary}
  \label{corollary.decompositionOfIdentity}
  Each choice of $q_1$ and $q_2$ as in
  Theorem~\ref{theorem.isomorphismCalibrated} determines in
  $\heckeWW[\clW]$ a maximal decomposition of the identity into
  idempotents, namely, $1=\sum_{w\in W} p_w$, where $p_w$ is the
  projection onto $E_w$, orthogonal to all $E_{w'}$, $w'\ne w$.
\end{corollary}
\begin{proof}
  Since $t$ is regular, one can construct each $p_w$ from
  $\cl(Y^{\coroot_1}),\ldots,\cl(Y^{\coroot_n})\in \heckeWW[\clW]$ by
  multivariate Lagrange interpolation. Therefore $p_w$ belongs to
  $\heckeWW[\clW]$.
\end{proof}

\iflongversion
\begin{lemma}
  \label{lemma.w0eigenvector}
  Let $\Wmax$ be the maximal element of $\clW$ in $\kclW$, and $\lavee$
  an element of the (finite) coroot lattice. Then $\Wmax$ is an
  eigenvector for $Y^\lavee$ with eigenvalue $q^{-\height(\lavee)}$.
\end{lemma}
\begin{proof}
  Let $i_1,\dots,i_r$ be an alcove walk for the translation
  $t_\lc$. Then, $s_{i_1}\cdots s_{i_r}$ acts trivially on the finite
  Weyl group: $\Wmax s_{i_1}\cdots s_{i_r}=\Wmax$.  Therefore,
  \begin{equation}
    \begin{split}
      \Wmax. Y^\lavee&
      = \Wmax.(-q_1q_2)^{-\height(\lc)} T_{i_1}^{\epsilon_1}\cdots T_{i_r}^{\epsilon_r} \\&
      = (-q_1q_2)^{-\frac12(\epsilon_1+\dots+\epsilon_r)} q_2^{(\epsilon_1+\dots+\epsilon_r)} \Wmax s_{i_1}\dots s_{i_r} \\&
      = \left(-\frac{q_1}{q_2}\right)^{-\frac12(\epsilon_1+\dots+\epsilon_r)} \Wmax
      = q^{-\height(\lc)}\Wmax\,,
    \end{split}
  \end{equation}
  using Equation~\eqref{equation.Y},
  Lemma~\ref{lemma.alcoveWalkTiOnW0}, and Remark~\ref{remark.height}.
\end{proof}
\fi

\iflongversion
As in \cite{Ram.2008}, define
$\tau_i:=T_i-\frac{q_1+q_2}{1-Y^{-\alpha_i^\vee}} \in \End(\kclW)$ for
$i=1,\dots,n$. Note that this operator is a priori only defined for
eigenvectors of $Y^{-\alpha_i^\vee}$ for an eigenvalue $\ne
1$. Whenever they are well-defined, they satisfy the braid relations,
as well as the following skew-commutation relation: $\tau_i Y^\lavee =
Y^{s_i(\lavee)} \tau_i$. Therefore, $\tau_i$ sends an $Y$-weight space
for the character $wt$ to an $Y$-weight space for the character
$ws_i t$.

For $w\in \clW$, define $E_w:=\Wmax.\tau_{i_1}\cdots \tau_{i_r}$ where
$i_1,\dots,i_r$ is any reduced word for $w$. 
\fi

\iflongversion
\begin{lemma}
  \label{lemma.eigenvectors}
  The $(E_w)_{w\in \clW}$ are well-defined, and triangular with respect to the
  canonical basis of $\kclW$:
  \begin{equation}
    E_w = (-q_1)^{\len(w)}  \Wmax w + \sum_{w'>\Wmax w} c_{w,w'} w'\,,
  \end{equation}
  for some coefficients $c_{w,w'}\in \CC$. In particular, the $E_w$ are
  all non-zero.
\end{lemma}
\begin{proof}
  The definition of $E_w$ does not depend on the choice of the reduced
  word thanks to the braid relations. Furthermore, at each step the
  application of $\tau_i$ on $E_{ws_i}$ is well-defined because
  $(ws_it){Y^{-\alpha_i^\vee}}\ne 1$.

  The triangularity is easily proved by induction: when $i$ is not a
  descent of $w$:
  \begin{equation}
    E_{ws_i} = E_w.\tau_i = E_w. \left( (q_1+q_2)\pi_i - q_1 s_i
      -\frac{q_1+q_2}{1-Y^{-\alpha_i^\vee}} \right)\,,
  \end{equation}
  and only the second term can contribute to the coefficient of $\Wmax
  w$.
\end{proof}
\fi

\begin{theorem}
  \label{theorem.quotientQualibrated}
  The morphism $\cl$ from the affine Hecke algebra $\heckeW{q_1,q_2}$
  to the Hecke group algebra $\heckeWW[\clW]$ is surjective for
  $q_1,q_2$ as in Theorem~\ref{theorem.isomorphismCalibrated}.
\end{theorem}
\begin{proof}
  Consider the decomposition $1=\sum_{w\in W} p_w$ of the identity of
  $\heckeWW[\clW]$ given in
  Corollary~\ref{corollary.decompositionOfIdentity}.

  Writing $(1-Y^{-\alpha_i^\vee})^{-1} = \sum_{w\in \clW} p_w
  (1-Y^{-\alpha_i^\vee})^{-1} = \sum_{w\in \clW} p_w
  (1-(wt)(Y^{-\alpha_i^\vee}))^{-1}$ shows that $1-Y^{-\alpha_i^\vee}$
  is invertible not only in $\End(\kclW)$ but even inside
  $\cl(\heckeW{q_1,q_2})$. Therefore
  $\tau_i=T_i-\frac{q_1+q_2}{1-Y^{-\alpha_i^\vee}}$ also belongs to
  $\cl(\heckeW{q_1,q_2})$.

  Consider the operator $p_w\tau_i$ which kills all eigenspaces
  $\CC.{E_{w'}}, w \ne w'$, and sends the eigenspace $\CC.E_w$ to
  $\CC.E_{ws_i}$. 

  The \emph{calibration graph} is the graph on $\W$ with an arrow from
  $w$ to $ws_i$ if $p_w\tau_i\ne 0$, or equivalently if
  $E_w.\tau_i\ne0$. We claim that this is the case if and only if
  $\Rec(w)\subset \Rec(ws_i)$. Take indeed $w\in W$ with a non-descent
  at position $i$. Then, $\Rec(w)\subset \Rec(ws_i)$ and by
  Lemma~\ref{lemma.eigenvectors}, $E_w.\tau_i = E_{ws_i} \ne 0$.
  Next, there is no arrow back from $ws_i$ to $w$ if and only if
  $E_w.(\tau_i)^2=0$. Using the quadratic relation satisfied by
  $\tau_i$,
  \begin{equation}
    \tau_i^2 = \frac{(q_1+q_2Y^{\alpha_i^\vee})(q_1+q_2 Y^{-\alpha_i^\vee})}
    {(1-Y^{\alpha_i^\vee}) (1-Y^{-\alpha_i^\vee}) }\,,
  \end{equation}
  this is the case if $\height(w(\coroot_i))=\pm1$. Since $i$ is not a
  descent of $w$, this is equivalent to $w(\coroot_i)=-\coroot_j$ for
  some simple coroot $\coroot_j$, that is $ws_i=s_jw$. In turn, this is
  equivalent to $\Rec(ws_i)=\Rec(s_jw)\supsetneq \Rec(w)$, which
  concludes the claim.

  For each $w$ and $w'$ with $\Rec(w) \subset \Rec(w')$ there exists a
  path $i_1,\dots,i_r$ from $w$ to $w'$ in the calibration graph;
  choose one, and set $\tau_{w,w'} = \tau_{i_1}\cdots\tau_{i_r}$.  The
  following family
  \begin{equation}
    \left\{ p_w\tau_{w, w'} \suchthat \Rec(w) \subset \Rec(w') \right\}
  \end{equation}
  is linearly independent, and by dimension comparison with $\heckeWW[\clW]$
  forms a basis $\cl(\heckeW{q_1,q_2})$. Therefore,
  $\cl(\heckeW{q_1,q_2})=\heckeWW[\clW]$.
\end{proof}

\section*{Acknowledgements}
\label{section.acknowledgements}

We would like to thank Jean-Yves Thibon for suggesting the
investigation of the connection between affine Hecke algebras and
Hecke group algebras. We are also very grateful to Masaki Kashiwara
for sharing his private notes on finite-dimensional representations of
quantized affine algebras with us, to Arun Ram for pointing out the
link with calibrated representations, to Francesco Brenti, Christophe
Holweg, Mark Shimozono, John Stembridge, and Monica Vazirani for
fruitful discussions, and to the anonymous referees for very helpful
suggestions.

This research was partially supported by NSF grants DMS-0501101,
DMS-0652641, and DMS-0652652. It started during a visit of the authors
at the University of California, San Diego in 2006, under the kind
invitation of Adriano Garsia and Richard and Isabelle Kauffmann. It
was completed while the third author was visiting the University of
California at Davis and during the inspiring 2008 MSRI Combinatorial
Representation Theory program.

The research was driven by computer exploration using the open-source
algebraic combinatorics package
\texttt{MuPAD-Combinat}~\cite{MuPAD-Combinat}. The pictures have been
produced (semi)-automatically, using \texttt{MuPAD-Combinat},
\texttt{graphviz}, \texttt{dot2tex}, and \texttt{pgf/tikz}.

\bibliographystyle{alpha}
\bibliography{main}
\end{document}